\documentclass[reqno]{amsart}

\usepackage{amssymb}

\usepackage[dvips,bookmarks,colorlinks=true,linkcolor=blue,urlcolor=red]{hyperref}
\newcounter{numcount}
\newcommand{\labelnummer}{\mbox{\normalfont (\roman{numcount})}}%

\makeatletter

   {\let\curlabelspeicher\@currentlabel%
     \begin{list}{\labelnummer}%
       {\usecounter{numcount}\leftmargin0pt%
         \topsep0.5ex\partopsep2ex\parsep0pt\itemsep0ex\@plus1\p@%
         \labelwidth2.5em\itemindent3.5em\labelsep1em%
       }%
     \let\saveitem\item%
     \def\item{\saveitem%
  \def\@currentlabel{\curlabelspeicher\hskip0.5pt\labelnummer}}%
     \let\savelabel\label%
     \def\label##1{\savelabel{##1}%
       \@bsphack%
         \ifmmode\else%
           \protected@write\@auxout{}%
           {\string\newlabel{##1item}{{\labelnummer}{\thepage}}}%
         \fi%
       \@esphack%
     }%
   }{\end{list}}%

\numberwithin{equation}{section}

\theoremstyle{plain}
\newtheorem{Th}{Theorem}[section]
\newtheorem{Le}{Lemma}[section]
\newtheorem{Pro}{Proposition}[section]

\theoremstyle{definition}
\newtheorem{Rem}{Remark}[section]
\newtheorem{Def}{Definition}[section]

\DeclareMathOperator{\supp}{supp}
\DeclareMathOperator{\tr}{tr}

\DeclareMathOperator{\dist}{dist}

\newcommand\R{\mathbb R}
\newcommand\N{\mathbb N}

\newcommand\Z{\mathbb Z}

\renewcommand\P{\mathbb P}
\newcommand\E{\mathbb E}

\newcommand{\esp}{\mathbb{E}}
\newcommand{\pro}{\mathbb{P}}

\newcommand{\vers}{\operatornamewithlimits{\to}}
\newcommand{\car}{\mathbf{1}}

\newcommand\K{\mathcal{K}}

\newcommand\I{\mathcal{I}}

\newcommand\e{\mathrm{e}}

\newcommand\eps{\varepsilon}

\newcommand{\set}[1]{\left\{ #1 \right\}}
\newcommand{\pa}[1]{\left( #1 \right)}

\newcommand\beq{\begin{equation}}
\newcommand\eeq{\end{equation}}

\begin{document}


\title[Spectral statistics in the localized regime] {Spectral
  statistics for random Schr{\"o}dinger operators in the localized regime}

\author{Fran\c cois Germinet} \address[Fran\c cois
Germinet]{Universit{\'e} de Cergy-Pontoise, CNRS UMR 8088, IUF,
  D{\'e}partement de math{\'e}matiques, F-95000 Cergy-Pontoise, France}
\email{\href{mailto:francois.germinet@u-cergy.fr}{francois.germinet@u-cergy.fr}}

\author{Fr{\'e}d{\'e}ric Klopp} 
\address[Fr{\'e}d{\'e}ric Klopp]{IMJ, UMR CNRS 7586 Universit{\'e} Pierre et Marie
  Curie Case 186, 4 place Jussieu F-75252 Paris cedex 05 France}
\email{\href{mailto:klopp@math.jussieu.fr}{klopp@math.jussieu.fr}}


\thanks{The authors are supported by the grant
  ANR-08-BLAN-0261-01. The authors would also like to thank the Banff
  International Research Station (Alberta, Canada) and the Centre
  Interfacultaire Bernoulli (EPFL, Lausanne), where this work was
  partially written, for their hospitality.\\
  It is pleasure to thank C. Shirley and N. veniaminov for pointing
  out various mistakes in previous versions of the paper.}


\begin{abstract}
  We study various statistics related to the eigenvalues and
  eigenfunctions of random Hamiltonians in the localized
  regime. Consider a random Hamiltonian at an energy $E$ in the
  localized phase. Assume the density of states function is not
  too flat near $E$. Restrict it to some large cube $\Lambda$. Consider
  now $I_\Lambda$, a small energy interval centered at $E$ that
  asymptotically contains infintely many eigenvalues when the volume
  of the cube $\Lambda$ grows to infinity. We prove that, with
  probability one in the large volume limit, the eigenvalues of the
  random Hamiltonian restricted to the cube inside the interval are
  given by independent identically distributed random variables, up to
  an error of size an arbitrary  power of the volume of the cube.\\
  As a consequence, we derive
  \begin{itemize}
  \item uniform Poisson behavior of the locally unfolded
    eigenvalues,
  \item a.s. Poisson behavior of the joint distibutions of the
    unfolded energies and unfolded localization centers in a large
    range of scales.
  \item the distribution of the unfolded level spacings, locally and
    globally,
  \item the distribution of the unfolded localization centers, locally
    and globally.
  \end{itemize}
  \vskip.5cm\noindent \textsc{R{\'e}sum{\'e}.}  Nous {\'e}tudions diff{\'e}rentes
  statistiques associ{\'e}es aux valeurs propres et vecteurs propres d'un
  op{\'e}rateur al{\'e}atoire dans le r{\'e}gime localis{\'e}. Consid{\'e}rons un
  op{\'e}rateur al{\'e}atoire au voisinage d'une {\'e}nergie $E$ suppos{\'e}e se
  trouver dans le r{\'e}gime localis{\'e}. On consid{\`e}re la restirction de
  l'op{\'e}rateur {\`a} un grand cube $\Lambda$. Soit $I_\Lambda$, un petit
  intervalle d'{\'e}nergie centr{\'e} en $E$ qui, asymptotiquement, contient
  une infinit{\'e} de valeurs propres quand $\Lambda$ grandit. Nous
  d{\'e}montrons qu'asymptotiquement presque s{\^u}rement, les valeurs propres
  de l'op{\'e}rateur al{\'e}atoire restreint au cube contenu dans l'intervalle
  sont tr{\`e}s bien approch{\'e}es par des variables al{\'e}atoires ind{\'e}pendantes
  identiquement distribu{\'e}es.\\
  De cette caract{\'e}risation, nous d{\'e}duisons
  \begin{itemize}
  \item que localement, les valeurs propres renormalis{\'e}es comportement
    uniform{\'e}ment poissonien,
  \item que les distibutions jointes des valeurs propres et centres de
    localisation renormalis{\'e}s sont asymptotiquement distribu{\'e}es selon
    une loi de Poisson,
  \item la distribution des espacements de niveaux renormalis{\'e}s,
    localement et globalement,
  \item la distribution des espacements des centres de localisation
    renormalis{\'e}s, localement et globalement.
  \end{itemize}
\end{abstract}

\setcounter{section}{-1}
\maketitle

\section{Introduction}
\label{sec:introduction}
To introduce our results, we first restrict our exposition to the
celebrated the random Anderson model, that is, on $\ell^2(\Z^d)$, we
consider the operator
\begin{equation*}
  H_{\omega}=-\Delta+V_\omega
\end{equation*}
where $-\Delta$ is the free discrete Laplace operator
\begin{equation*}
  (-\Delta u)_n=\sum_{|m-n|=1}u_m\quad\text{ for }
  u=(u_n)_{n\in\Z^d}\in\ell^2(\Z^d)
\end{equation*}
and $V_\omega$ is the random potential
\begin{equation*}
  (V_\omega u)_n=\omega_n u_n\quad\text{ for }u=(u_n)_{n\in\Z^d}
  \in\ell^2(\Z^d).
\end{equation*}
We assume that the random variables $(\omega_n)_{n\in\Z^d}$ are
independent identically distributed and that their distribution admits
a compactly supported smooth density, say $g$.\\
It is well known (see e.g.~\cite{MR2509110}) that, $\omega$ almost
surely, the spectrum of $H_\omega$ is equal to a fixed closed set,
say, $\Sigma$. Moreover, there exists a Lebesgue almost everywhere
bounded density of states, say $\lambda\mapsto\nu(\lambda)$, such
that, for any continuous function $\varphi:\ \R\to\R$, one has
\begin{equation*}
  \int_\R\varphi(\lambda)\nu(\lambda)d\lambda=
  \mathbb{E}(\langle\delta_0,\varphi(H_\omega)\delta_0\rangle).
\end{equation*}
The function $\nu$ is the density of a probability measure on
$\Sigma$.
\par For $L>1$, consider $\Lambda=[-L,L]^d\cap\Z^d$, a cube on the lattice
and let $H_\omega(\Lambda)$ be the random Hamiltonian $H_\omega$
restricted to $\Lambda$ with periodic boundary conditions. It is a
finite dimensional symmetric matrix; let us denote its eigenvalues
ordered increasingly and repeated according to multiplicity by
$E_1(\omega,\Lambda)\leq E_2(\omega,\Lambda)\leq \cdots\leq
E_N(\omega,\Lambda)$. For $x\geq0$, define the empirical level spacings
distribution of $H_\omega(\Lambda)$ as
\begin{equation*}
  DLS(x;\omega,\Lambda)=\frac{\#\{j;\
    (E_{j+1}(\omega,\Lambda)-E_j(\omega,\Lambda))|\Lambda|\geq x\} 
  }{|\Lambda|}.
\end{equation*}
A result that is typical of the results we prove in the present paper
is
\begin{Th}
  \label{thr:17}
  There exists $\lambda_0>0$ such that, for $|\lambda|>\lambda_0$,
  with probability 1, as $|\Lambda|\to+\infty$, $DLS(x;\omega,\Lambda)$
  converges uniformly to the distribution $x\mapsto g(x)$ where
  \begin{equation}
    g(x)=\int_{\Sigma}\e^{-\nu(\lambda)x}\nu(\lambda)d\lambda,
  \end{equation} 
  that is, $\omega$ almost surely,
  \begin{equation*}
    \sup_{x\geq0}\left|DLS(x;\omega,\Lambda)
      -g(x)\right|\vers_{|\Lambda|\to+\infty}0.
  \end{equation*}
\end{Th}
\noindent This result shows that, for the discrete Anderson
Hamiltonian with smoothly distributed random potential at sufficiently
large coupling, the limit of the level spacings distribution is that of
i.i.d. random variables distributed according to the density of states
of the random Hamiltonian.
\vskip.2cm\noindent To the best of our knowledge, this is the first
rigorous study of the level spacings distribution of random Schr{\"o}dinger operators.
\vskip.2cm\noindent The purpose of the paper is to study spectral
statistics for random Hamiltonians in the localized regime. The large
coupling Anderson Hamiltonian described above is the typical
example. Spectral statistics have been studied in various works,
mainly for discrete or continuous Anderson models (see
e.g.~\cite{MR84e:34081,MR97d:82046,MR2299191,MR2002i:82054}) but up to
now, to the best of our knowledge, studies have only described the
local spectral statistics. For a random Hamiltonian restricted to a
cube $\Lambda$, the existence of the density of states, that is, of a
limit for the number of eigenlevels per unit of volume
(see~\eqref{eq:83}) implies that the average distance between levels
is of order $|\Lambda|^{-1}$. Local spectral statistics are the
statistics in energy intervals $I$ of the size $|\Lambda|^{-1}$. Thus,
such intervals contain typically a number of eigenvalues that is
bounded uniformly in the volume of the cube. This, in particular,
prevents the study of the empirical distribution of level spacings.\\
In the present paper, we go beyond this. Therefore, we introduce a new
method of study of the eigenlevels and localization centers that is
quite close to the physical heuristics (see
e.g.~\cite{citeulike:693492,RevModPhys.57.287,citeulike:3832118}). The
method consists in approximating the eigenvalues of the true random
operator (restricted to some finite cube) by a family of i.i.d. random
variables that are constructed as eigenvalues of the random operator
restricted to smaller cubes. That this is possible is a consequence of
localization: due to their exponential falloff, eigenfunctions only
see the random potential surrounding them. This construction is only
feasible under some restrictions on the relative size of the region
where one wants to study eigenvalues and the size of the cube on which
one restricts the random operator. If one wants to control, with a
good probability, all the eigenvalues in some interval $I$, then, one
roughly needs $I$ to be of size $|\Lambda|^{-\alpha}$, the inverse of
the volume of the cube to some power $\alpha$ smaller than but close
to $1$ (see Theorem~\ref{thr:vsmall1}). If one wants to enlarge the
interval $I$, one can go up to sizes that are of order a negative
power of $\log|\Lambda|$ at the expense of being able to describe only
most of the eigenvalues (see Theorem~\ref{thr:vbig1}).\\
The basic tools that we use to control the eigenvalues are the
so-called ``Wegner'' and ``Minami'' estimates (see (W) and (M) in
section~\ref{sec:random-model}). \vskip.1cm
Using the approximation described above, we obtain a large deviation
estimate for the number of eigenvalues inside a possibly shrinking
interval of a random operator restricted to some large cube (see
Theorem~\ref{thr:16}). This bound shows that, for intervals $I$ that
are not too small, with good probability, the number of eigenvalues in
$I$ is given by the weight that the integrated density of states gives
to $I$ times the volume of the cube up to an error of smaller order. \\
Then, we derive the almost sure level spacings statistics near fixed
energies as well as inside non trivial compact intervals (see
Theorems~\ref{thr:1}, \ref{thr:9} and~\ref{thr:DCS}). We also compute
the localization center spacings statistics (see
Theorem~\ref{thr:DCS}).\vskip.1cm\noindent
This is the first time that these statistics are obtained for random
Schr{\"o}dinger operators.\vskip.1cm\noindent
The next result is the uniform local statistic for the eigenvalues and
localization centers when they are rescaled covariantly i.e. the
scaling in energy is of order the scaling in space to the power $-d$
(see Theorem~\ref{thr:2}); we prove that the covariantly scaled local
statistics are independent of the scale (if they are not too small)
i.e. one always obtains Poisson statistics. In the non-covariant
scaling case, we obtain almost sure results on the counting function
(see Theorem~\ref{thr:7}). In the case of the standard scale
i.e. energies are scaled by $|\Lambda|^{-1}$ on a cube of volume
$\Lambda$, we also study the asymptotic independence of these local
processes (see Theorems~\ref{thr:4} and~\ref{thr:8}). This extends
known results of \cite{MR84e:34081,MR97d:82046,MR2299191}.\\
We point out that our analysis goes beyond the previous results also
in the sense that locally the images of the eigenvalues by the IDS are
shown to exhibit a Poissonian behavior. When the DOS exists and is non
zero, we recover the previous known statements, but a vanishing
derivative of the IDS is allowed in the present
work.\vskip.1cm
We also consider the problem from a different point of view. The usual
procedure consists in restricting the random Hamiltonian to some
finite cube and study the statistics for this operator in the limit
when the cube grows to be the whole space. One can also consider the
Hamiltonian in the whole space. In a compact interval in the localized
region, say $I$, the Hamiltonian admits countably many eigenvalues. We
enumerate them using the localization center attached to an associated
eigenfunction (see Proposition~\ref{pro:1}). I.e. we consider the
eigenvalues in $I$ having localization in some finite cube. We derive
the almost sure statistics of the level spacings distributions (see
Theorem~\ref{thr:10}); they are the same as the ones obtained for the
Hamiltonian restricted to a cube.  \vskip.1cm
Finally let us us conclude this introduction by saying that a number
of the results obtained in the present paper for general random
Schr{\"o}dinger operators were already described for discrete random
operators in~\cite{Ge-Kl:10b}.



\section{The main results}
\label{sec:results}
After this short illustration of what can be obtained from our method,
we now turn to the description of the main results of this
paper. Results will be given for general random Schr{\"o}dinger operators
under a number of assumptions that are known to hold for, e.g. the
Anderson model and some continuous Anderson models.\\
We shall use the following standard notations: $a\lesssim b$ means
there exists $c<\infty$ so that $a\le cb$; $a\asymp b$ means
$a\lesssim b$ and $b \lesssim a$; $\langle x\rangle =
(1+|x|^2)^\frac12$.
\subsection{The random model}
\label{sec:random-model}
Consider $H_\omega=H_0+V_\omega$, a $\Z^d$-ergodic random
Schr{\"o}\-dinger operator on $\mathcal{H}=L^2(\R^d)$ or $\ell^2(\Z^d)$
(see e.g.~\cite{MR94h:47068,MR1935594}). Typically, the background
potential $H_0$ is the Laplacian $-\Delta$, possibly perturbed by a
periodic potential. Magnetic fields can be considered as well; in
particular, the Landau Hamiltonian is also admissible as a background
Hamiltonian. For the sake of simplicity, we assume that $V_\omega$ is
almost surely bounded; hence, almost surely, $H_\omega$ share the same
domain $H^2(\R^d)$ or $\ell^2(\Z^d)$.\\
For $\Lambda$, a cube in either $\R^d$ or $\Z^d$, we let
$H_\omega(\Lambda)$ be the operator $H_\omega$ restricted to $\Lambda$
with periodic boundary conditions. Our analysis stays valid for
Dirichlet boundary conditions.\\
Furthermore, we shall denote by $\car_J(H)$ the spectral projector of
the operator $H$ on the energy interval $J$. $\E(\cdot)$ denotes the
expectation with respect to $\omega$;
Our first assumption will be an independence assumption on the local
Hamiltonian that is
\begin{description}
\item[(IAD)] Independence At a Distance: there exists $R_0>0$ such
  that for dist$(\Lambda,\Lambda')> R_0$, the random Hamiltonians
  $H_\omega(\Lambda)$ and $H_\omega(\Lambda')$ are independent.
\end{description}
Such an assumption is clearly satisfied by standard models like the
Anderson model, the Poisson model or the random displacement model if
the single site potential is compactly supported (see
e.g.~\cite{MR94h:47068,MR1935594}).
\begin{Rem}
  \label{rem:8}
  As will be seen in the course of the proofs, it can be weakened to
  assuming that the correlations between the local Hamiltonians decay
  polynomially at a rate that is faster than the $-d$-th power of the
  distance separating the cubes on which the local Hamiltonians are
  considered.
\end{Rem}
Let $\Sigma$ be the almost sure spectrum of $H_\omega$. Pick $I$ a
relatively compact open subset of $\Sigma$. Assume the following
holds:
\begin{description}
\item[(W)] a Wegner estimate holds in $I$, i.e. given $I$ there exists
  $C>0$ such that, for $J\subset I$, and $\Lambda$, a cube in $\R^d$
  or $\Z^d$, one has
  \begin{equation}
    \label{eq:1}
    \E\left[\text{tr}(\car_J(H_\omega(\Lambda)))
    \right]\leq C |J|\,|\Lambda| .
  \end{equation}
 
\item[(M)] a Minami estimate holds in $I$, i.e. given $I$ there exists
  $C>0$ and $\rho>0$ such that, for $J\subset I$, and $\Lambda$, a
  cube in $\R^d$ or $\Z^d$, one has
  \begin{equation}
    \label{eq:2}
    \E\left[\text{tr}(\car_J(H_\omega(\Lambda)))
      \cdot[\text{tr}(\car_J(H_\omega(\Lambda)))-1]\right]\leq
    C (|J|\,|\Lambda|)^{1+\rho}.
  \end{equation}
\end{description}
\begin{Rem}
  \label{rem:1}
  The Wegner estimate has been proved for many random Schr{\"o}dinger
  models both discrete and continuous Anderson models under rather
  general conditions on the single site potential and on the
  randomness (see e.g.~\cite{MR2509108,MR2307751,MR2378428}). The
  right hand side in~\eqref{eq:1} can be lower bounded by the
  probability to have at least one eigenvalue in $J$ (for $J$
  small).\\
  As the proofs will show, one can weaken assumption (W) and replace
  the right hand side with $C |J|^\alpha\,|\Lambda|^\beta$ for
  arbitrary positive $\alpha$ and $\beta$. Such Wegner estimates are
  known to hold also for some non monotonous models (see
  e.g.~\cite{MR95m:82080,MR1934351,MR2423576}). \vskip.1cm
  On the Minami estimate, much less is known: it holds for the
  discrete Anderson model with $I=\Sigma$
  (see~\cite{MR97d:82046,MR2290333,MR2360226,MR2505733}). These proofs yield an optimal exponent $\rho=1$. \\
  In dimension 1, it has been proved recently (see~\cite{Kl:11a}),
  that, for general random models, the Minami estimate (for any
  $\rho\in(0,1)$) is a consequence of the Wegner estimate within the
  localization region (see section~\ref{sec:localized-regime}).\\
  In higher dimensions, for continuous Anderson models, proving a Minami estimate is still a challenging open problem.
  The right hand side in~\eqref{eq:2} can be lower bounded by the
  probability to have at least two eigenvalues in $J$. For $\rho=1$,
  it behaves as the square of the probability to have one eigenvalue
  in $J$. So, roughly speaking, close by eigenvalues behave as
  independent random variables.
\end{Rem}
The integrated density of states is defined as
\begin{equation}
  \label{eq:83}
  N(E):=\lim_{|\Lambda|\to+\infty}\frac{\#\{\text{e.v. of
    }H_\omega(\Lambda)\text{ less than E}\}}{|\Lambda|} .
\end{equation}
By (W), $N(E)$ is the distribution function of a measure that is
absolutely continuous with respect to the Lebesgue measure on
$\R$. Let $\nu$ be the density of state of $H_\omega$ i.e. the
distributional derivative of $N$. In the sequel, for a set $I$, we
will often write $N(I)$ for the mass the measure $\nu(E)dE$ puts on
$I$ i.e.
\begin{equation}
  \label{defIDS}
  N(I)=\int_I\nu(E)dE.  
\end{equation}
\subsection{The localized regime}
\label{sec:localized-regime}
For $L\geq1$, $\Lambda_L$ denotes the cube $[-L/2,L/2]^d$ in either
$\R^d$ or $\Z^d$. Let $\mathcal{H}_\Lambda$ be
$\ell^2(\Lambda\cap\Z^d)$ in the discrete case and $L^2(\Lambda)$
in the continuous one. For a vector $\varphi\in\mathcal{H}$, we define
\begin{equation}
  \label{eq:22}
  \|\varphi\|_x=
  \begin{cases}
    \|\car_{\Lambda(x)}\varphi\|_2\text{ where
    }\Lambda(x)=\{y;|y-x|\leq1/2\}&\quad\text{
    if }\mathcal{H}=L^2(\R^d),\\
  \quad\quad|\varphi(x)|&\quad\text{ if }\mathcal{H}=\ell^2(\Z^d).
  \end{cases}
\end{equation}
In the discrete case, the definition is that given in
section~\ref{sec:local-cent-stat}.\\
Let $I$ be a compact interval. We assume that $I$ lies in the region
of complete localization (see e.g.~\cite{MR2078370,MR2203782}) for
which we use the following finite volume version:
\begin{description}
\item[(Loc)] for all $\xi\in(0,1)$, one has
  \begin{equation}
    \label{eq:84}
    \sup_{L>0}\,\sup_{\substack{\text{supp} f\subset I \\ |f|\leq1}}\,
      \esp\left(\sum_{\gamma\in\Z^d} \e^{|\gamma|^\xi}\,
        \|\car_{\Lambda(0)}f(H_\omega(\Lambda_L))
        \car_{\Lambda(\gamma)}\|_2\right)<+\infty.
  \end{equation}
\end{description}
Whenever the fractional moment method is available, one may replace
the factor $\e^{|\gamma|^\xi}$ by an exponential one
$\e^{\eta|\gamma|}$, where $\eta>0$. 
\begin{Rem}
  \label{rem:11}
  Assumption (Loc) may be relaxed into asking~\eqref{eq:84} for a
  single $\xi$. This will not change the subsequent results in an
  essential way, but only modify some constants.
\end{Rem}
We note that the assumption (Loc) implies in particular that the
spectrum of $H_\omega$ is pure point in $I$ (see
e.g.~\cite{MR2203782,MR2509110}).  We refer to
Appendix~\ref{sec:append-local-contr} where, in Theorem~\ref{thmFVE},
we provide equivalent finite volume properties of the region of
complete localization, and show it coincides with the infinite volume
one.  For the sake of the exposition, from Theorem~\ref{thmFVE}, we
extract the following lemma that we shall use intensively in this
paper.
\begin{Le}
  \label{le:3}
  Assume (W) and (Loc).
  \\
  (I) For any $p>0$ and $\xi\in(0,1)$, for $L\geq1$ large enough,
  there exists a set of configuration $\mathcal{U}_{\Lambda_L}$ such
  that $\pro(\mathcal{U}_{\Lambda_L})\geq1-L^{-p}$ and for
  $\omega\in\mathcal{U}_{\Lambda_L}$, if
  \begin{enumerate}
  \item $\varphi_j(\omega, \Lambda_L)$ is a normalized eigenvector of
    $H_{\omega}(\Lambda_L)$ associated to $E_j(\omega,\Lambda_L)\in I$,
  \item $x_j(\omega, \Lambda_L)\in \Lambda_L$ is a maximum of
    $x\mapsto\|\varphi_j(\omega, \Lambda_L)\|_x$ in $\Lambda_L$,
  \end{enumerate}
  then, for $x\in\Lambda_L$, one has
  \begin{equation}
    \label{eq:19}
    \|\varphi_j(\omega, \Lambda_L)\|_x\leq L^{p+d} \e^{-|x-x_j(\omega,
      \Lambda_L)|^\xi}. 
  \end{equation}
  (II)  For any $\nu,\xi\in(0,1)$, $\nu<\xi$, for $L\geq1$ large enough, there
  exists a set of configuration $\mathcal{V}_{\Lambda_L}$ such that
  $\pro(\mathcal{V}_{\Lambda_L})\geq 1-e^{-L^\nu}$,  and for
  $\omega\in\mathcal{V}_{\Lambda_L}$, if
  \begin{enumerate}
  \item $\varphi_j(\omega, \Lambda_L)$ is a normalized eigenvector of
    $H_{\omega}(\Lambda_L)$ associated to $E_j(\omega, \Lambda_L)\in I$,
  \item $x_j(\omega, \Lambda_L)\in \Lambda_L$ is a maximum of
    $x\mapsto\|\varphi_j(\omega, \Lambda_L)\|_x$ in $\Lambda_L$,
  \end{enumerate}
  then, for $x\in\Lambda_L$, one has
  \begin{equation}
    \label{eq:95}
    \|\varphi_j(\omega, \Lambda_L)\|_x\leq e^{2L^\nu}
    e^{-|x-x_j(\omega, \Lambda_L)|^{\xi}}. 
  \end{equation}
\end{Le}
\begin{Rem}
  Both Part (I) and (II) of Lemma~\ref{le:3} are consequences of the
  localization hypothesis. We shall use both of these
  characterizations of localization. Part (I) is relevant for large
  scales (typically, powers of the volume of the reference box) and
  yields a smaller constant in front of the exponential, while part
  (II) will be used for smaller scales (typically powers of $\log$ of
  the volume of the box) and yields a better probability.
\end{Rem}
\noindent Such a result can essentially be found in~\cite{MR2203782}
for the continuous case and in~\cite{Kl:10} for the discrete case. \\
Clearly, the function $x\mapsto\|\varphi_j(\omega, \Lambda_L)\|_x$
need not have a unique maximum in $\Lambda_L$. But, as, for any
$x\in\Lambda_L$, one has
\begin{equation*}
  \sum_{\gamma\in\Lambda_L\cap\Z^d}\|\varphi_j(\omega,
  \Lambda_L)\|^2_{x+\gamma}
  =\|\varphi_j(\omega, \Lambda_L)\|^2=1,
\end{equation*}
if $x_j(\omega, \Lambda_L)$ is a maximum, then $\|\varphi_j(\omega,
\Lambda_L)\|^2_{x_n(\omega)}\geq (2L+1)^{-d}$. Hence, if $x_j(\omega,
\Lambda_L)$ and $x'_j(\omega, \Lambda_L)$ are two maxima, then (Loc),
through Lemma~\ref{le:3}(I), implies that, for any $p$, there exists
$C_p>0$ such that, with a probability larger than $1-L^{-p}$, we have
\begin{equation*}
  |x_j(\omega, \Lambda_L)-x'_j(\omega, \Lambda_L)|\leq C_p(\log L)^{1/\xi}.
\end{equation*}
For $\varphi\in\mathcal{H}_\Lambda$, define the set of localization
centers for $\varphi$ as
\begin{equation}
  \label{eq:12}
  C(\varphi)=\{x\in\Lambda;\
 \|\varphi\|_x=\max_{\gamma\in\Lambda}  \|\varphi\|_{\gamma}\}.
\end{equation}
As a consequence of Lemma~\ref{le:3}, one has
\begin{Le}
  \label{le:1}
  Pick $I$ in the localized regime for $H_\omega$.  For any $p>0$ ,
  there exists $C_p>0$, such that, with probability larger than
  $1-L^{-p}$, if $E_j(\omega,\Lambda_L)\in I$ then the diameter of
  $C(\varphi_j(\omega,\Lambda_L))$ is less than
  $C_p\log^{1/\xi}|\Lambda|$.
\end{Le}
\noindent From now on, a localization center for a function $\varphi$
will denote any point in the set of localization centers $C(\varphi)$
and let $x_j(\omega,\Lambda_L)$ be a localization center for
$\varphi_j(\omega,\Lambda_L)$. One can e.g. order them
lexicographically and pick the one with largest coefficients.

\subsection{The asymptotic description of the eigenvalues}
\label{sec:asympt-descr-eigenv}
We now state our main results. They are also the main technical
results on which we base all our studies of the spectral statistics.\\
They consist in a precise approximation of the eigenvalues of
$H_\omega(\Lambda)$ in $I_\Lambda$ by independent random variables,
that follows in a rather straightforward way from the standard
properties of random Schr{\"o}dinger operators recalled above ((IAD), (W),
(M) and localization). This approximation is at the heart of the
proofs of the new statistical results we present in this paper: from
it, we derive estimations on the number of eigenvalues in small
intervals, we provide the first computation of the level spacings
distribution, and we extend the known results about the convergence to
Poisson of rescaled (or unfolded) eigenvalues. \\
We will give two different descriptions depending on the size of
$N(I_\Lambda)$. When this quantity is sufficiently small with
respect to $|\Lambda|^{-1}$, our procedure enables us to control all
the eigenvalues. If it is not, we only control most of the
eigenvalues. \\
Recall that our cube of reference is $\Lambda=\Lambda_L$, with center
$0$ and sidelength $L$.
\subsubsection{Controlling all the eigenvalues}
\label{sec:contr-all-eigenv}
To start, pick $\tilde\rho$ such that
\begin{equation}
  \label{condrho}
  0\leq\tilde{\rho}<\frac{\rho}{1+d\rho}.
\end{equation}
Assume $E_0$ is such that~\eqref{eq:60} holds. Now, pick $I_\Lambda$
centered at $E_0$ such that $N(I_\Lambda)\asymp|\Lambda|^{-\alpha}$
for $\alpha\in(\alpha_{d,\rho,\tilde\rho},1)$ where
$\alpha_{d,\rho,\tilde\rho}$ is defined as
\begin{equation}
  \label{eq:59}
  \alpha_{d,\rho,\tilde\rho}:=(1+\tilde\rho)\frac{d\rho+1}{d\rho+1+\rho},
\end{equation}
where $\tilde\rho>0$ and $\rho$ is defined in the Minami estimate
(M). Assumption~\eqref{condrho} clearly implies that
$\alpha_{d,\rho,\tilde\rho}<1$.\\
Our restriction will enable us to control all the eigenvalues of
$H_\omega(\Lambda)$ in $I_\Lambda$.
\begin{Th}
  \label{thr:vsmall1}
  Assume $E_0$ is such that~\eqref{eq:60} holds for some
  $\tilde\rho\in[0,\rho/(1+d\rho))$. Recall that
  $\alpha_{d,\rho,\tilde\rho}$ is defined in~\eqref{eq:59} and pick
  $\alpha\in(\alpha_{d,\rho,\tilde\rho},1)$. Pick $I_\Lambda$ centered
  at $E_0$ such that $N(I_\Lambda)\asymp|\Lambda|^{-\alpha}$. There
  exists $\beta>0$ and $\beta'\in(0,\beta)$ small so that
  $1+\beta\rho<\alpha\frac{1+\rho}{1+\tilde\rho}$ and, for $\ell\asymp
  L^\beta$ and $\ell'\asymp L^{\beta'}$, there exists a decomposition
  of $\Lambda$ into disjoint cubes of the form
  $\Lambda_\ell(\gamma_j):=\gamma_j+[0,\ell]^d$ satisfying:
  \begin{itemize}
  \item $\cup_j\Lambda_\ell(\gamma_j)\subset\Lambda$,
  \item $\dist (\Lambda_\ell(\gamma_j),\Lambda_\ell(\gamma_k))\ge
    \ell'$ if $j\not=k$,
  \item $\dist (\Lambda_\ell(\gamma_j),\partial\Lambda)\ge \ell'$
  \item $|\Lambda\setminus\cup_j\Lambda_\ell(\gamma_j)|\lesssim |
    \Lambda|\ell'/\ell$,
  \end{itemize}
  and such that, for $L$ sufficiently large, there exists a set of
  configurations $\mathcal{Z}_ \Lambda $ s.t.:
  \begin{itemize}
  \item $\P(\mathcal{Z}_\Lambda)\geq 1 - |
    \Lambda|^{-(\alpha-\alpha_{d,\rho,\tilde\rho})}$,
  \item for $\omega\in\mathcal{Z}_\Lambda $, each centers of
    localization associated to $H_\omega(\Lambda)$ belong to some
    $\Lambda_\ell(\gamma_j)$ and each box $\Lambda_\ell(\gamma_j)$
    satisfies:
    \begin{enumerate}
    \item the Hamiltonian $H_\omega(\Lambda_\ell(\gamma_j))$ has at
      most one eigenvalue in $I_\Lambda $, say,
      $E_j(\omega,\Lambda_\ell(\gamma_j))$;
    \item $\Lambda_\ell(\gamma_j)$ contains at most one center of
      localization, say $x_{k_j}(\omega,\Lambda)$, of an eigenvalue of
      $H_\omega(\Lambda)$ in $I_\Lambda $, say $E_{k_j}(\omega,\Lambda)$;
    \item $\Lambda_\ell(\gamma_j)$ contains a center
      $x_{k_j}(\omega,\Lambda)$ if and only if
      $\sigma(H_\omega(\Lambda_\ell(\gamma_j)))\cap
      I_\Lambda\not=\emptyset$; in which case, one has
      \begin{equation}
        \label{error}
        |E_{k_j}(\omega,\Lambda)-E_j(\omega,\Lambda_\ell(\gamma_j))| \leq  
        e^{-(\ell')^{\xi}}\text{ and }\mathrm{dist}(x_{k_j}(\omega,\Lambda),
        \Lambda \setminus \Lambda_\ell(\gamma_j))\geq \ell'.
      \end{equation}
    \end{enumerate}
  \end{itemize}
  In particular if $\omega\in\mathcal{Z}_\Lambda$, all the eigenvalues
  of $H_\omega(\Lambda)$ are described by~\eqref{error}.
\end{Th}
\noindent With a probability tending to $1$, Theorem~\ref{thr:vsmall1}
describes all the eigenvalues of $H_\omega(\Lambda)$ inside a
sufficiently small interval $I_\Lambda$ as i.i.d. random variables
defined as the unique eigenvalue of a copy of the random Hamiltonian
$H_\omega(\Lambda_\ell(0))$ inside $I_\Lambda$.\\
As one can easily imagine, this description yields the local
statistics for both eigenlevels and localization center. Actually as
the intervals under consideration are larger than $|\Lambda|^{-1}$, it
yields uniform local statistics (see
sections~\ref{sec:strong-poiss-conv}
and~\ref{sec:strong-poiss-conv-1}). Theorem~\ref{thr:vsmall1} is at
the heart of the proof the proofs of
Theorems~\ref{thr:3},~\ref{thr:2},~\ref{thr:11},~\ref{thr:5}
and~\ref{thr:7} found in section~\ref{sec:loc-level-statistics}
and~\ref{sec:local-cent-stat}. Moreover, under additional the
decorrelations estimates (see assumptions (GM) and (D) in
section~\ref{sec:asympt-indep-local}), Theorem~\ref{thr:vsmall1} will
be sufficient to prove the mutual independence of the local processes
at distinct energies when they are sufficiently far apart, that is,
when they are separated by a distance that is asymptotically infinite
with respect to $|\Lambda|^{-1}$ (see
section~\ref{sec:asympt-indep-local}).
\subsubsection{Controlling most eigenvalues}
\label{sec:contr-most-eigenv}
Then, we now state a result that works on intervals $I_\Lambda$ such
that $N(I_\Lambda)$ be of size $(\log|\Lambda|)^{-d/\xi}$ but gives
the control only on most of the eigenvalues. This is enough to control
the levelspacings on such sets. This is the main tool to obtain
Theorems~\ref{thr:1},~\ref{thr:9},~\ref{thr:DCS} and~\ref{thr:10}. For
that purpose we state it in a more axiomatic way than what we did for
Theorems~\ref{thr:vsmall1}.
\begin{Def}
  \label{def:1}
  Pick $\xi\in(0,1)$, $R>1$ large and $\rho'\in(0,\rho)$ where $\rho$
  is defined in (M). For a cube $\Lambda$, consider an interval
  $I_\Lambda=[a_\Lambda,b_\Lambda]\subset I$. Set $\ell'_\Lambda=
  (R\log |\Lambda|)^{\frac 1\xi}$. We say that the sequence
  $(I_\Lambda)_\Lambda$ is $(\xi,R,\rho')$-admissible if, for any
  $\Lambda$, one has
  \begin{equation}
    \label{eq:93}
    |\Lambda| N(I_\Lambda)\geq1,\quad
    N(I_\Lambda)|I_\Lambda|^{-(1+\rho')}\geq1, \quad 
    N(I_\Lambda)^{\frac 1{1+\rho'}} (\ell'_\Lambda)^d\leq 1.
  \end{equation}
\end{Def}
\noindent One has
\begin{Th}
  \label{thr:vbig1}
  Assume (IAD), (W), (M) and (Loc) hold. Pick
  $\rho'\in[0,\rho/(1+(\rho+1)d))$ where $\rho$ is defined in (M). For
  any $q>0$, for $L$ sufficiently large, depending only on $\xi,
  R,\rho',p$, for any sequence of intervals $(I_\Lambda)_\Lambda$ that
  is $(\xi,R,\rho')$-admissible, and any sequence of scales
  $\tilde\ell_\Lambda$ so that $\ell'_\Lambda\ll\tilde\ell_\Lambda\ll
  L$ and
  \begin{equation}
    \label{eq:102}
    N(I_\Lambda)^{\frac
      1{1+\rho'}}\tilde\ell_\Lambda^d\vers_{|\Lambda|\to\infty}0, 
  \end{equation}
  there exists
  \begin{itemize}
  \item a decomposition of $\Lambda$ into disjoint cubes of the form
    $\Lambda_{\ell_\Lambda}(\gamma_j):=\gamma_j+[0,\ell_\Lambda]^d$,
    where $\displaystyle\ell_\Lambda=\tilde\ell_\Lambda(1+\mathcal{O}
    (\tilde\ell_\Lambda / |\Lambda|))=\tilde\ell_\Lambda(1+o(1))$ such
    that
    \begin{itemize}
    \item $\cup_j\Lambda_{\ell_\Lambda}(\gamma_j)\subset\Lambda$,
    \item $\dist
      (\Lambda_{\ell_\Lambda}(\gamma_j),\Lambda_{\ell_\Lambda}(\gamma_k))
      \ge \ell'_\Lambda$ if $j\not=k$,
    \item $\dist (\Lambda_{\ell_\Lambda}(\gamma_j),\partial\Lambda)\ge
      \ell'_\Lambda$
    \item
      $|\Lambda\setminus\cup_j\Lambda_{\ell_\Lambda}(\gamma_j)|\lesssim
      | \Lambda| \ell'_\Lambda/\ell_\Lambda$,
    \end{itemize}
  \item a set of configurations $\mathcal{Z}_\Lambda$ such that
    \begin{itemize}
    \item $\mathcal{Z}_\Lambda$ is large, namely,
      \begin{equation}
        \label{eq:94}
        \pro(\mathcal{Z}_\Lambda)\geq 1 - |\Lambda|^{-q}- 
        \exp\left(-c|I_\Lambda|^{1+\rho}|\Lambda|
          \ell_\Lambda^{d\rho}) \right)\\-\exp\left( -c
          |\Lambda||I_\Lambda|\ell'_\Lambda\ell_\Lambda^{-1}\right)
      \end{equation}
    \end{itemize}
  \end{itemize}
  so that
  \begin{itemize}
  \item for $\omega\in\mathcal{Z}_\Lambda$, there exists at least
    $\displaystyle\frac{|\Lambda|}{\ell_\Lambda^d}\left(1+
      o\left(N(I_\Lambda)^{\frac1{1+\rho'}}
          \ell_\Lambda^d\right)\right)$ disjoint boxes
      $\Lambda_{\ell_\Lambda}(\gamma_j)$ satisfying the properties
      (1), (2) and (3) described in Theorem~\ref{thr:vsmall1} with
      $\ell'_\Lambda=(R\log |\Lambda|)^{1/\xi}$ in~\eqref{error}; we
      note that $N(I_\Lambda)\ell^{d-1}_\Lambda\ell'_\Lambda=o(1)$ as
      $|\Lambda|\to+\infty$;
  \item the number of eigenvalues of $H_\omega(\Lambda)$ that are not
    described above is bounded by 
    \begin{equation}\label{controlev}
      C N(I_\Lambda)|\Lambda|
      \left(N(I_\Lambda)^{\frac{\rho-\rho'}{1+\rho'} }
        \ell_\Lambda^{d(1+\rho)}+
        N(I_\Lambda)^{-\frac{\rho'}{1+\rho'}}
        (\ell'_\Lambda)^{d+1}\ell_\Lambda^{-1}\right);
    \end{equation}
    this number is $o(N(I_\Lambda) |\Lambda|)$ provided
    \begin{equation}
      \label{condell} 
      N(I_\Lambda)^{-\frac{\rho'}{1+\rho'}}
      (\ell'_\Lambda)^{d+1} \ll \ell_\Lambda \ll
      N(I_\Lambda)^{- \frac{\rho-\rho'}{d(1+\rho)(1+\rho')}} .
    \end{equation}
  \end{itemize}
\end{Th}
\noindent Before turning to the description of the statistics global
and local, let us make two remarks about the choice of parameters and
length scales for which Theorem~\ref{thr:vbig1} is useful. First, if
$N(I_\Lambda)^{-1}$ is of an order much larger that of
$(\ell_\Lambda')^d$, then, the condition~\eqref{condell} essentially
imposes that $\displaystyle\frac{\rho-\rho'}{(1+\rho')(1+\rho)}>
\frac{d\rho'}{1+\rho'}$ which is satisfied if
$\rho'\in[0,\rho/(1+(\rho+1)d))$. Condition~\eqref{condell} then
guarantees that the condition~\eqref{eq:102} is met: indeed,
$\tilde\ell_\Lambda\sim\ell_\Lambda$ and, as $\rho-\rho'<1+\rho$ and
$N(I_\Lambda)\to0$, one has
\begin{equation*}
  N(I_\Lambda)^{-1/(1+\rho')}\ll
  N(I_\Lambda)^{- \frac{\rho-\rho'}{(1+\rho)(1+\rho')}}.
\end{equation*}
We shall use this in the proof of the large deviation
estimate~\eqref{eq:58} in conjunction with a choice of $\ell_\Lambda$
and $N(I_\Lambda)^{-1}$ as large powers of $\ell'_\Lambda$, that is,
large powers of $\log|\Lambda|$ (see section~\ref{sec:local-count-funct}).\\
Note that, if $\rho'\in[0,\rho/(1+d(\rho+1)))$, then,~\eqref{condell}
and~\eqref{eq:102} are satisfied for some choice of $\alpha\in(0,1)$
and $\nu\in(0,1/d)$ if one sets
$\ell'_\Lambda\asymp(\log|\Lambda|)^{1/\xi}$, $\ell_\Lambda\asymp
N(I_\Lambda)^{-\nu}$ and $N(I_\Lambda)\asymp |\Lambda|^{-\alpha}$ (see
section~\ref{sec:some-prel-cons}). This is
the choice of parameters we shall use in our study of level spacings.\\
In \cite{Ge-Kl:10a}, for a less general class of models, by improving
on the assumptions (W) and (M), we shall improve on the
bound~\eqref{controlev}, which
\begin{itemize}
\item will enable us to relax the second condition in~\eqref{eq:93} so
  as to admit $N(I_\Lambda)$ of size $e^{-|I_\Lambda|^{-\alpha}}$ for
  $\alpha\in(0,1)$;
\item will provide a better large deviation estimate for the
number of eigenvalues of $H_\omega(\Lambda)$ in the interval
$I_\Lambda$ than the first rough estimate given in
Theorem~\ref{thr:16}.
\end{itemize}
\subsubsection{Comparing various cube sizes}
\label{sec:comp-vari-cube}
We end this section with  a related result we shall use in the sequel.
We prove
\begin{Pro}
  \label{pro:2}
  Assume (IAD) (Loc) and (W) in $J$. Fix $I\subset J$ a compact
  interval and $C$, $C'$, compact cubes in $\R^d$ such that
  $\overline{C}\subset\overset{\circ}{C'}$. Fix $p>0$ and a sequence
  $(\ell_\Lambda)_\Lambda$ satisfying~\eqref{eq:14}. Then, there
  exists a set of configurations $\mathcal{Z}_\Lambda$ so that
  $\P(\mathcal{Z}_\Lambda)\to1$ as $|\Lambda|\to+\infty$ and, for
  $\omega\in\mathcal{Z}_\Lambda$ and $|\Lambda|$ sufficiently large,
  \begin{itemize}
  \item to every eigenvalue of $H_\omega(\Lambda)$ in
    $E_0+\ell_\Lambda^{-d} I$ associated to a localization center in
    $\ell_\Lambda C$, say $E_j(\omega,\Lambda)$, one can associate an
    eigenvalue of $H_\omega(\ell_\Lambda C')$, say
    $E_j(\omega,\ell_\Lambda C')$; moreover, these eigenvalues satisfy
    \begin{equation}
      \label{eq:103}
      |E_j(\omega,\ell_\Lambda C')-E_j(\omega,\Lambda)|\leq |\Lambda|^{-p}.
    \end{equation}
  \item to every eigenvalue of $H_\omega(\ell_\Lambda C')$ in
    $E_0+\ell_\Lambda^{-d} I$ associated to a localization center in
    $\ell_\Lambda C$, say $E_j(\omega,\ell_\Lambda C')$, one can
    associate an eigenvalue of $H_\omega(\Lambda)$, say
    $E_j(\omega,\Lambda)$ with localization center in $\ell_\Lambda
    C'$; moreover, these eigenvalues satisfy
    \begin{equation}
      \label{eq:21}
      |E_j(\omega,\ell_\Lambda C')-E_j(\omega,\Lambda)|\leq |\Lambda|^{-p}.
    \end{equation}
  \end{itemize}
\end{Pro}
\noindent Proposition~\ref{pro:2} is a consequence of Lemma~\ref{lemcenter}.
\begin{Rem}
  \label{rem:9}
  As the proofs will show, the sizes of the intervals where the
  control of the eigenvalues is possible and the probability of the
  event where this control is possible both depend very much on the
  forms of the Wegner and Minami estimates, (W) and (M). In
  particular, if one replaces (W) by what is suggested in
  Remark~\ref{rem:8}, the constants appearing in
  Theorems~\ref{thr:vsmall1} and~\ref{thr:vbig1} and
  Proposition~\ref{pro:2} have to be modified.
\end{Rem}
\subsection{The level spacings statistics}
\label{sec:level-spa-statistics}
\noindent Our goal is now to understand the level spacings statistics
for eigenvalues near $E_0\in I$. Pick $I_\Lambda$ a compact interval
containing $E_0$ such that its density of states measure
$N(I_\Lambda)$ stays bounded. We note that, by the existence of the
density of states and also Theorem~\ref{thr:3}, the spacing between
the image of the eigenvalues of $H_\omega(\Lambda)$ through $N$ near
$E_0$ is of size $|\Lambda|^{-1}$. Hence, to study the empirical
statistics of level spacings in $I_\Lambda$, $N(I_\Lambda)$ should
contain asymptotically infinitely many images of energy levels of
$H_{\omega}(\Lambda)$. Let us first study this number.
\subsubsection{A large deviation principle for the eigenvalue counting
  function}
\label{sec:large-devi-princ}
Define the random numbers
\begin{equation}
  \label{eq:26}
  N(I_\Lambda,\Lambda,\omega):=\#\{j;\ E_j(\omega,\Lambda)\in
  I_\Lambda\}.
\end{equation}
Write $I_\Lambda=[a_\Lambda,b_\Lambda]$ and recall that
$N(I_\Lambda)=N(b_\Lambda)-N(a_\Lambda)$ where $N$ is the integrated
density of states. We show that $N(I_\Lambda,\Lambda,\omega)$
satisfies a large deviation principle, namely,
\begin{Th}
  \label{thr:16}
  Assume (IAD), (W), (M) and (Loc) hold. Fix
  $\tilde\rho\in(0,\rho/(1+d(\rho+1)))$ where $\rho$ is given by (M),
  and $\nu\in(0,1)$. Then, there exists $\delta>0$ small such that, if
  $(I_\Lambda)_\Lambda$ is a sequence of compact intervals in the
  localization region $I$ satisfying
  \begin{itemize}
  \item $N(I_\Lambda)\,(\log|\Lambda|)^{1/\delta}\to0$ as
    $|\Lambda|\to+\infty$
  \item $N(I_\Lambda)\,|\Lambda|^{1-\nu}\to+\infty$ as
    $|\Lambda|\to+\infty$
  \item $N(I_\Lambda)\,|I_\Lambda|^{-1-\tilde\rho}\to+\infty$ as
    $|\Lambda|\to+\infty$,
  \end{itemize}
  then, for any $p>0$, for $|\Lambda|$ sufficiently large (depending
  on $\rho'$ and $\nu$ but not on the specific sequence
  $(I_\Lambda)_\Lambda$), one has
  \begin{equation}
    \label{eq:58}
    \pro\left(\left|N(I_\Lambda,\Lambda,\omega)-N(I_\Lambda)|\Lambda|
      \right|\geq N(I_\Lambda)|\Lambda|(\log|\Lambda|)^{-\delta}\right)\leq
    |\Lambda|^{-p}.
  \end{equation}
\end{Th}
\noindent We note that we do not need that intervals
$(I_\Lambda)_\Lambda$ lie near points $E_0$ where~\eqref{eq:60} is
satisfied; the density of states may vanish near $E_0$ though not
faster than the rate fixed by the condition
$N(I_\Lambda)|I_\Lambda|^{-1-\rho}\to+\infty$. The large deviation
principle~\eqref{eq:58} is meaningful only if
$N(I_\Lambda)|\Lambda|\to+\infty$; as $N$ is Lipschitz continuous as a
consequence of (W), this implies that
\begin{equation}
  \label{eq:23}
  |\Lambda|\cdot|I_\Lambda|\vers
  +\infty\quad\text{ when }\quad|\Lambda|\to+\infty .
\end{equation}
For the discrete Anderson model, we improve upon~\eqref{eq:58} in
\cite{Ge-Kl:10a} by relaxing
$N(I_\Lambda)|I_\Lambda|^{-1-\rho}\to+\infty$ into
$N(I_\Lambda)|I_\Lambda|^{-\nu}\to+\infty$ for some arbitrarily large
$\nu>0$ and obtain precise estimates of the term
$o(N(I_\Lambda)|\Lambda|)$, exploiting an improved Wegner and Minami
estimate.
\subsubsection{The level spacing statistics near a given energy}
\label{sec:loc-level-spa-stat}
Fix $E_0\in I$. Pick $I_\Lambda=[a_\Lambda,b_\Lambda]$ so that
$|a_\Lambda|+|b_\Lambda|\to0$. Consider the unfolded eigenvalue
spacings, for $1\leq j\leq N$,
\begin{equation}
  \label{eq:3}
  \delta N_j(\omega,\Lambda)=|\Lambda|(N(E_{j+1}(\omega,\Lambda))-
  N(E_j(\omega,\Lambda)))\geq0.
\end{equation}
Define the empirical distribution of these spacings to be the random
numbers, for $x\geq0$
\begin{equation}
  \label{eq:4}
  DLS(x;E_0+I_\Lambda,\omega,\Lambda)=\frac{\#\{j;\
    E_j(\omega,\Lambda)\in E_0+I_\Lambda,\ \delta
    N_j(\omega,\Lambda)\geq x\}
  }{N(I_\Lambda,\Lambda,\omega)}.
\end{equation}
We will now study the spacings distributions of energies inside
intervals that shrink to a point but that asymptotically contain
infinitely many eigenvalues.\\
We prove
\begin{Th}
  \label{thr:1}
  Assume (IAD), (W), (M) and (Loc) hold.  Fix $E_0\in I$ such that,
  for some $\tilde\rho\in[0,\rho/(1+d(\rho+1)))$, there exists a
  neighborhood of $E_0$, say, $U$ such that
  \begin{equation}
    \label{eq:86}
    \forall (x,y)\in U^2,\quad |N(x)-N(y)|\geq |x-y|^{1+\tilde\rho}.
  \end{equation}
  Fix $(I_\Lambda)_\Lambda$ a decreasing sequence of intervals such
  that $\displaystyle\sup_{E\in I_\Lambda}|E|
  \vers_{|\Lambda|\to+\infty}0$.\\
  Assume that, for some $\delta>0$, one has
  \begin{equation}
    \label{eq:43}
    |\Lambda|^{1-\delta}\cdot N(E_0+I_\Lambda)
    \vers_{|\Lambda|\to+\infty}    +\infty
    \ \text{ and }\ \text{if }\ell'=o(L)\text{ then
    }\frac{N(E_0+I_{\Lambda_{L+\ell'}})}{N(E_0+I_{\Lambda_L})}
    \vers_{|\Lambda|\to+\infty}1.
  \end{equation}
  Then, with probability $1$, as $|\Lambda|\to+\infty$,
  $DLS(x;E_0+I_\Lambda,\omega,\Lambda)$ converges uniformly to the
  distribution $x\mapsto e^{-x}$, that is, with probability $1$,
  \begin{equation}
    \label{eq:5}
    \sup_{x\geq0}\left|DLS(x;E_0+I_\Lambda,\omega,\Lambda)
          -e^{-x}\right|\vers_{|\Lambda|\to+\infty}0.
  \end{equation}
\end{Th}
\noindent Hence, the unfolded level spacings behave as if the images
of the eigenvalues by the IDS were i.i.d. uniformly distributed random
variables (see \cite{MR0070874} or section 7 of~\cite{MR0216622}). The
exponential distribution of the level spacings is the one predicted by
physical heuristics in the localized regime
(\cite{citeulike:693492,RevModPhys.57.287,citeulike:3832118,Th:74}). It
is also in accordance with Theorem~\ref{thr:3}.
In~\cite{MR84e:34081,MR97d:82046}, the domains in energy where
the statistics were studied were much smaller than the ones
considered in Theorem~\ref{thr:1}. Indeed, in these works, the energy
interval is of order $|\Lambda|^{-1}$ whereas, here, it is assumed to
tend to $0$ but to be asymptotically infinite when compared to
$|\Lambda|^{-1}$.
\begin{Rem}
  \label{rem:7}
  The first condition in~\eqref{eq:43} ensures that $I_\Lambda$
  contains sufficiently many eigenvalues of $H_\omega(\Lambda)$. The
  second condition in~\eqref{eq:43} is a regularity condition of the
  decay of $|I_\Lambda|$.\\
  If, in~\eqref{eq:43}, one replaces the first condition by $|\Lambda|
  N(E_0+I_\Lambda)\to+\infty$ or omits the second or does both, one
  still gets convergence in probability of
  $DLS(x;E_0+I_\Lambda,\omega,\Lambda)$ to $e^{-x}$ (see
  Remark~\ref{rem:2}) i.e.
  \begin{equation*}
    \P\left(\sup_{x\geq0}\left|DLS(x;E_0+I_\Lambda,\omega,\Lambda)
        -e^{-x}\right|\geq\varepsilon\right)\vers_{|\Lambda|\to+\infty}0.
  \end{equation*}  
\end{Rem}
\noindent Condition~\eqref{eq:86} is slightly stronger
than~\eqref{eq:60}; it requires some uniformity in the lower
bound. Theorem~\ref{thr:1} can be applied to obtain the levelspacing
distribution near regular points of the IDS. Define $\mathcal{E}$ to
be the set of energies $E$ such that $\nu(E)=N'(E)$ exists and
\begin{equation}\label{defE}
  \lim_{|x|+|y|\to0}\frac{N(E+x)-N(E+y)}{x-y}=\nu(E).
\end{equation}
Obviously the set $\mathcal{E}$ contains the continuity points of
$\nu(E)$. We prove in Appendix 
\begin{Le}
  \label{le:5}
  The set $\mathcal{E}$ is of full Lebesgue measure.
\end{Le}
\noindent For $E_0\in\mathcal{E}$ such that $\nu(E_0)>0$,
assumption~\eqref{eq:86} holds with $\tilde\rho=0$. Moreover, one
has
\begin{equation*}
  \delta E_j(\omega,\Lambda)=\nu(E_0)|\Lambda|(E_{j+1}(\omega,\Lambda)-
  E_j(\omega,\Lambda))(1+o(1)).
\end{equation*}
Then as a corollary of Theorem~\ref{thr:1}, we immediately obtain
\begin{Th}
  \label{thr:12}
  Assume (IAD), (W), (M) and (Loc) hold.  Fix $E_0\in
  I\cap\mathcal{E}$ such that $\nu(E_0)>0$. Fix $(I_\Lambda)_\Lambda$
  a sequence of intervals such that
  $\displaystyle\sup_{I_\Lambda}|x|\vers_{|\Lambda|\to+\infty}0$.\\
  Assume that, for some $\delta>0$, one has
  \begin{equation*}
    |\Lambda|^{1-\delta}\cdot |I_\Lambda|
    \vers_{|\Lambda|\to+\infty}    +\infty
    \ \text{ and }\ \text{if }\ell'=o(L)\text{ then
    }\frac{|I_{\Lambda_{L+\ell'}}|}{|I_{\Lambda_L}|}
    \vers_{|\Lambda|\to+\infty}1.
  \end{equation*}
  Then, with probability $1$, as $|\Lambda|\to+\infty$, the empirical
  distribution function
  \begin{equation*}
    \frac{\#\{j;\ E_j(\omega,\Lambda)\in
      I_\Lambda,\ \nu(E_0)\,
      |\Lambda|(E_{j+1}(\omega,\Lambda)-E_j(\omega,\Lambda))\geq x\}
    }{N(I_\Lambda,\Lambda,\omega)}
  \end{equation*}
  converges uniformly to the distribution $x\mapsto e^{-x}$.
\end{Th}
\subsubsection{The level spacings statistics on macroscopic energy
  intervals}
\label{sec:level-spac-stat-1}
Theorem~\ref{thr:12} seems optimal as the density of states at $E_0$
enters into the correct rescaling to obtain a universal result. Hence,
the distribution of level spacings on larger intervals needs to take
into account the variations of the density of states on these
intervals. Indeed, on intervals of non vanishing size, under
additional regularity assumption on $\nu$, one can compute the
asymptotic distribution of the level spacings when one omits the local
density of states in the spacing and obtain the
\begin{Th}
  \label{thr:9}
  Assume (IAD), (W), (M) and (Loc) hold.  Pick $J\subset I$ a compact
  interval such $\lambda\mapsto\nu(\lambda)$ be continuous on $J$ and
  $N(J):=\int_J\nu(\lambda)d\lambda>0$. Define the unfolded eigenvalue
  spacings, for $1\leq j\leq N$,
  \begin{equation}
    \label{eq:24}
    \delta_J E_j(\omega,\Lambda)= \frac{N(J)}{|J|} |\Lambda|
    (E_{j+1}(\omega,\Lambda)-E_j(\omega,\Lambda))\geq 0 ,
  \end{equation}
  and the empirical distribution of these spacings to be the random
  numbers, for $x\geq0$,
  \begin{equation}
    \label{eq:25}
    DLS'(x;J,\omega,\Lambda)=\frac{\#\{j;\ E_j(\omega,\Lambda)\in
      J,\ \delta_J E_j(\omega,\Lambda)\geq x\}
    }{N(J,\omega,\Lambda)}.
  \end{equation}
  Then, with probability 1, as $|\Lambda|\to+\infty$,
  $DLS'(x;J,\omega,\Lambda)$ converges uniformly to the distribution
  $x\mapsto g_{\nu,J}(x)$ where
  \begin{equation}
    \label{eq:31}
    g_{\nu,J}(x)=\int_{J}\e^{-\nu_J(\lambda)|J|x}\nu_J(\lambda)
    d\lambda\text{ where }\nu_J=\frac{1}{N(J)}\nu.
  \end{equation}
\end{Th}
\noindent We see that, in the large volume limit, the unfolded
levelspacings behave as if the eigenvalues were i.i.d. random
variables distributed according to the density
$\frac{1}{N(J)}\nu(\lambda)$ i.e. to the density of states
renormalized to be a probability measure on $J$ (see section 7
of~\cite{MR0216622}).\vskip.1cm
\noindent Theorem~\ref{thr:17} is then an immediate consequence of
Theorem~\ref{thr:9} and the results on the regularity for the density
of states of the discrete Anderson model at large disorder obtained
in~\cite{MR89c:82050}. We point out that for random Hamiltonians in
the continuum, the continuity of the density of states is still an
open problem.
\subsection{The localization center spacings statistics}
\label{sec:center-statistics}
Pick $E_0\in I$. Inside a large cube $\Lambda$, the number of centers
that corresponds to energies in $I_\Lambda$,
$N(E_0+I_\Lambda,\Lambda,\omega)$ (see~\eqref{eq:26}), is asymptotic
to $N(E_0+I_\Lambda)|\Lambda|$. Theorem~\ref{thr:7} states that they
are distributed uniformly. Thus, the reference mean spacing between
localization centers is of size
$(|\Lambda|/N(E_0+I_\Lambda)|\Lambda|)^{1/d}
=(N(E_0+I_\Lambda))^{-1/d}$. This motivates the following
definition. Define the empirical distribution of center spacing to be
the random number
\begin{equation}
  \label{DCS}
  DCS(s;I_\Lambda,\Lambda,\omega)=\frac{\#\left\{j;\
      \begin{aligned}
        &E_j(\omega,\Lambda)\in I_\Lambda,\\
        \sqrt[d]{N(E_0+I_\Lambda)}\displaystyle&\min_{i\neq
          j}|x_i(\omega)-x_j(\omega)| \geq s
      \end{aligned}
    \right\} }
  {N(E_0+I_\Lambda,\Lambda,\omega)}.
\end{equation}
We prove an analogue of Theorem~\ref{thr:1}, namely
\begin{Th}
  \label{thr:DCS}
  Assume (IAD), (W), (M) and (Loc) hold.  Pick $E_0\in I$ such that,
  for some $\tilde\rho\in[0,\rho/(1+d(\rho+1)))$ small enough
  (depending on $\rho$ and $d$), in $U$ some neighborhood of $E_0$,
  one has~\eqref{eq:86}. Assume that, for some $\nu\in(0,1)$, one has
  \begin{equation}
    \label{eq:73}
    N(E_0+I_\Lambda)\,|\Lambda|\vers_{|\Lambda|\to+\infty}+\infty\quad
    \text{and}\quad N(E_0+I_\Lambda)\,\log^{d/\nu}|\Lambda|
    \vers_{|\Lambda|\to+\infty}0.
  \end{equation}
  Then, as $|\Lambda|\to+\infty$, in probability,
  $DCS(s;I_\Lambda,\Lambda,\omega)$ converges uniformly to the
  distribution $x\mapsto e^{-s^d}$, that is, for any $\varepsilon>0$,
  \begin{equation}
    \label{eq:82}
    \P\left(\left\{\omega;\
        \sup_{s\geq0}\left|DCS(s;E_0+I_\Lambda,\Lambda,\omega)-e^{-s^d} 
        \right|\geq\varepsilon\right\}\right)\vers_{\Lambda\nearrow\R^d} 0.
  \end{equation}
\end{Th}
\noindent Of course, as Theorem~\ref{thr:DCS} is the counterpart of
Theorem~\ref{thr:1}, Theorem~\ref{thr:9} also has its counterpart for
localization centers.
\subsection{Another point of view}
\label{sec:another-point-view}
In the present section, we want to adopt a different point of view on
the spectral statistics. Instead of discussing the statistics of the
eigenvalues of the random system restricted to some finite box in the
large box limit, we will describe the spectral statistics of the
infinite system in the localized phase. Therefore, we first need to
explicit what we mean with the localized phase for the random
Hamiltonian on the whole space i.e. state the appropriate replacement
for assumption (Loc) in this setting.\\
Let $I\subset\R$ be an interval. We assume
\begin{description}
\item[(Loc')] there exists $\xi\in(0,1]$, $q>0$ and $\gamma>0$ such that,
  with probability $1$, if $E\in I\cap\sigma(H_\omega)$ and $\varphi$
  is a normalized eigenfunction associated to $E$ then, for some
  $x(E,\omega)\in\R^d$ or $\Z^d$, a maximum of $x\mapsto\|\varphi\|_x$, for
  some $C_\omega>0$, one has, for $x\in\R^d$,
    \begin{equation}
      \label{eq:7}
      \|\varphi\|_x\leq C_\omega (1+|x(E,\omega)|^2)^{q/2}e^{-\gamma
        |x-x(E,\omega)|^\xi};
    \end{equation}
    moreover, one has $\esp(C_\omega)<+\infty$.
\end{description}
As above $x(E)$ is called {\it a center of localization} for energy
$E$ or for the associated  eigenfunction $\varphi$.\\
It is well established that (Loc') holds in any interval contained in
the region of complete localization. The first proof is due
to~\cite{MR97m:47002} for the discrete Anderson model where they show
that (Loc') is a consequence of the fractional moment method
\cite{MR1301371,MR1244867}; there $\xi=1$. The proof extends to
continuous Hamiltonians thanks to
\cite{MR2207021}. In~\cite{MR1627657}, the multiscale analysis is
shown to imply \eqref{eq:7} with $\xi=1$ but with no control on
$\E(C_\omega)$. That the multiscale analysis yields (Loc') for any
$\xi<1$ follows from \cite[Corollary~3]{MR2203782} and
\cite[Eq.~(4.17)]{MR2203782} to see that $\E(C_\omega)<\infty$.
\par Pick $I$ an interval where the Hamiltonian $H_\omega$ is
localized i.e. satisfies (Loc'). Assume that, $\omega$-almost surely,
$\sigma(H_\omega)\cap I=I$. Hence, any subinterval of $I$ contains
infinitely many eigenvalues and to define statistics, we need to
enumerate these eigenvalues in a way or another. To do this, we use
the localization centers. First, we prove
\begin{Pro}
  \label{pro:1}
  Assume (Loc') for some $\xi\in(0,1]$ and fix $q>2d$. Then, there
  exists $\gamma>0$ such that, $\omega$-almost surely, there exists
  $C_\omega>1$, $\E(C_\omega)<\infty$, such that
  \begin{enumerate}
  \item if $x(E,\omega)$ and $x'(E,\omega)$ are two centers of
    localization for $E\in I$ then, for some $C_d>0$ (depending only
    on $d$),
    \begin{equation}
      \label{eq:8}
      |x(E,\omega)-x'(E,\omega)|\leq \gamma^{-1/\xi}\log^{1/\xi}
      \left(C_d C_\omega \langle x(E,\omega)\rangle^{q}\langle 
        x'(E,\omega)\rangle^{q}\frac{1}{\gamma^{d/\xi}}\right). 
    \end{equation}
  \item for $L\geq1$, pick $(I_L)_L$ a sequence of intervals,
    $I_L\subset I$, such that, for some $\varepsilon>0$, one has
    $L^{d-\varepsilon}N(I_L)\to+\infty$ and
    $N(I_L)|I_L|^{-1-\rho}\to+\infty$ where $\rho$ is given by
    $\mathrm{(M)}$ and $N(I_L)$ by \eqref{defIDS}; if $N(I_L,L)$
    denotes the number of eigenvalues of $H_\omega$ having a center of
    localization in $\Lambda_L$, then
    \begin{equation}
      \label{eq:9}
      N(I_L,L)=N(I_L)\,|\Lambda_L|\,(1+o(1)).
    \end{equation}
  \end{enumerate}
\end{Pro}
\noindent For $L\geq1$, pick $I_L\subset I$ such that $L^{d-\eps}
N(I_L)\to +\infty$ for some $\eps>0$. In view of
Proposition~\ref{pro:1}, there are only finitely many eigenvalues of
$H_\omega$ in $I_L$ having a localization center in $\Lambda_L$. Thus,
we can consider their levelspacings: let us enumerate these
eigenvalues as $E_1(\omega,L)\leq E_2(\omega,L)\leq \cdots\leq
E_N(\omega,L)$ where we repeat them according to multiplicity. Define
the empirical distributions $DLS$ and $DLS'$ as in~\eqref{eq:4}
and~\eqref{eq:25} for the eigenvalues of $H_\omega$ in $I_L$ having a
localization center in $\Lambda_L$. We prove
\begin{Th}
  \label{thr:10}
  Assume (IAD), (W), (M), (Loc') hold. One has
  \begin{itemize}
  \item if $E_0\in I_L$ s.t.~\eqref{eq:60} is satisfied for some
    $\tilde\rho\in[0,\rho/(1+d\rho))$ and $|I_L|\to0$ and
    satisfies~\eqref{eq:43}, then, $\omega$-almost surely, for
    $x\geq0$,
    \begin{equation}
      \label{eq:30}
      \lim_{L\to+\infty}DLS(x;I_L,\omega,L)=e^{-x};
    \end{equation}
  \item if, for all $L$ large, $|I_L|=J$ such that $N(J)>0$ and $\nu$
    is continuous on $J$ then, $\omega$-almost surely, one has
  \begin{equation}
    \label{eq:10}
    \lim_{L\to+\infty}DLS'(x;I_L,\omega,L)=g_J(x),
  \end{equation}
  where $g_J$ is defined in~\eqref{eq:31}.
\end{itemize}
\end{Th}
\noindent We see that the level spacings distribution of the
eigenvalues of $H_\omega$ having center of localization in $\Lambda_L$
have the same limits as those of the eigenvalues of
$H_\omega(\Lambda)$. This is a consequence of the localization
assumption (Loc').
\subsection{The local level statistics}
\label{sec:loc-level-statistics}
We now exploit of approximation of eigenvalues by iid ones to revisit
and extend previous results on the convergence to the Poisson law of
rescaled eigenvalues and centers of localization.\\
For $L\in\N$, recall that $\Lambda=\Lambda_L$ and that
$H_{\omega}(\Lambda)$ is the operator $H_\omega$ restricted to
$\Lambda$ with periodic boundary conditions. The notation
$|\Lambda|\to+\infty$ is a shorthand for considering
$\Lambda=\Lambda_L$ in the limit $L\to+\infty$.
\par Denote the eigenvalues of $H_\omega(\Lambda)$ ordered
increasingly and repeated according to multiplicity by
$E_1(\omega,\Lambda)\leq E_2(\omega,\Lambda)\leq \cdots\leq
E_N(\omega,\Lambda)\leq\cdots $.
\par Let $E_0$ be an energy in $I$.  The {\it unfolded local level
  statistics} near $E_0$ is the point process defined by
\begin{equation}
  \label{eq:13}
  \Xi(\xi; E_0,\omega,\Lambda) = 
  \sum_{j\geq1} \delta_{\xi_j(E_0,\omega,\Lambda)}(\xi),
\end{equation}
where
\begin{equation}
  \label{eq:11}
  \xi_j(E_0,\omega,\Lambda)=|\Lambda|(N(E_j(\omega,\Lambda))-N(E_0)).
\end{equation}
The numbers $(|\Lambda| N(E_j(\omega,\Lambda)))_j$ are called the
unfolded eigenvalues of $H_\omega(\Lambda)$ (see
e.g.~\cite{MR2352280,Mi:08} for more details).\\
The unfolded local level statistics are described by
\begin{Th}
  \label{thr:3}
  Assume (IAD), (W), (M) and (Loc) hold.  Pick $\tilde\rho$
  satisfying~\eqref{condrho} where $\rho$ is defined by (M).\\
  Pick $E_0$ be an energy in $I$ such that the integrated density of
  states satisfies
  \begin{equation}
    \label{eq:60}
    \forall a>b,\ \exists C(a,b)>0,\ \exists \varepsilon_0>0,
    \forall\varepsilon\in(0,\varepsilon_0),\quad
    |N(E_0+a\varepsilon)-N(E_0+b\varepsilon)|\geq C(a,b)
    \varepsilon^{1+\tilde\rho}.
  \end{equation}
  When $|\Lambda|\to+\infty$, the point process
  $\Xi(E_0,\omega,\Lambda)$ converges weakly to a Poisson process on
  $\R$ with intensity the Lebesgue measure.
\end{Th}
\noindent If one assume that $N$ is differentiable at $E_0$ and that
its derivative $\nu(E_0)$ is positive i.e.
\begin{equation}
  \label{eq:85}
  0<\nu(E_0):=\lim_{E\to E_0}\frac{N(E)-N(E_0)}{E-E_0},
\end{equation}
it is easy to check that~\eqref{eq:60} is satisfied with
$\tilde{\rho}=0$ and that, for $E_j-E_0$ small,
\begin{equation*}
 \xi_j(E_0,\omega,\Lambda)=
 |\Lambda|\nu(E_0)(E_j(\omega,\Lambda)-E_0)(1+o(1)).
\end{equation*}
Thus, one recovers the convergence to a Poisson process when the point
process~\eqref{eq:13} is replaced by the one defined by the points
$(|\Lambda|\nu(E_0)(E_j(\omega,\Lambda)-E_0))_j$.\\
Theorem~\ref{thr:3} under the additional assumption~\eqref{eq:85} was
first obtained in~\cite{MR84e:34081} for a special one dimensional random
Schr{\"o}dinger model on the real line. For the discrete Anderson model,
Theorem~\ref{thr:10} was proved in~\cite{MR97d:82046} under the
assumption~\eqref{eq:85}.\\
Our method of proof is different from that
of~\cite{MR97d:82046} and, in spirit, closer to that
of~\cite{MR84e:34081} and to the physical heuristics. Clearly,
for~\eqref{eq:60} to be satisfied, we do not need $N$ to be
differentiable at $E_0$ nor its derivative to be positive. E.g. if $N$
satisfies $N(E)=N(E_0)+c(E-E_0)^{1+\tilde\rho}(1+o(1))$ near $E_0$
then~\eqref{eq:60} is satisfied. Condition~\eqref{eq:60} asks that at
a given scale, $N$ behaves roughly uniformly near $E_0$. Note however
that, if $N$ is not differentiable at $E_0$, then the local statistics
of the eigenvalues themselves will not be Poissonian
anymore. \vskip.1cm
\noindent Our method yields a uniform version of Theorem~\ref{thr:3}
to which we now turn.
\subsubsection{Uniform Poisson convergence over small intervals}
\label{sec:strong-poiss-conv}
Fix $\alpha\in(\alpha_{d,\rho,\tilde\rho},1)$ (recall that
$\alpha_{d,\rho,\tilde\rho}$ is defined in~\eqref{eq:59}). The uniform
version of the Poisson process is a version that holds uniformly over
an interval of energy, say, $I$ centered at $E_0$ such that
$N(I)\asymp|\Lambda|^{-\alpha}$. Such an interval is much larger than
an interval satisfying $N(I)\asymp|\Lambda|^{-1}$.  This is the main
improvement of Theorem~\ref{thr:2} below over Theorem~\ref{thr:3} or
the statements found in~\cite{MR2299191,MR97d:82046,MR84e:34081}. It
is natural to wonder what is the largest size of interval in which a
result like Theorem~\ref{thr:2}. We do not know the answer to that
question.\\
Let $I_\Lambda(E_0,\alpha)$ be the interval such that
$N(I_\Lambda(E_0,\alpha))$ be centered at $N(E_0)$ of length
$2|\Lambda|^{-\alpha}$.  Denote by $N_\Lambda(\omega,E_0):=\tr
\car_{I_\Lambda(E_0,\alpha)}(H_\omega(\Lambda))$ the number of
eigenvalues of $H_\omega(\Lambda)$ in $I_\Lambda(E_0,\alpha)$.  For
$1\leq j\leq N_\Lambda(\omega,E_0)$, define the unfolded local
eigenvalues $\xi_j(\omega,\Lambda)$ by~\eqref{eq:11}.  Hence, for all
$1\leq j\leq N_\Lambda(\omega,E_0)$, one has
$\xi_j(\omega,\Lambda)\in |\Lambda|^{1-\alpha}\cdot[-1,1]$.\\
We then prove
\begin{Th}
  \label{thr:2}
  Assume (IAD), (W), (M) and (Loc) hold. Let $E_0$ be an energy in $I$
  such that, for some $\tilde\rho$ such that~\eqref{condrho} holds
  true and such that
  \begin{equation}
    \label{eq:116}
    \tilde{\rho}\geq\rho\,\frac{1-d\rho}{1+d\rho}. 
  \end{equation}
  the integrated density of states satisfies
  \begin{equation}
    \label{eq:107}
    \begin{split}
    \forall \delta\in(0,1),\ \exists C(\delta)>0,\ \exists \varepsilon_0>0,\ 
    &\forall\varepsilon\in(0,\varepsilon_0),\ \forall a\in[-1,1],\\
    &|N(E_0+(a+\delta)\varepsilon)-N(E_0+a\varepsilon)|\geq C(\delta)
    \,\varepsilon^{1+\tilde\rho}.      
    \end{split}
  \end{equation}
  Pick $\alpha\in(\alpha_{d,\rho,\tilde\rho},1)$. Then, there exists
  $\delta>0$ such that, for any sequence of intervals
  $I_1=I^\Lambda_1,\dots,I_p=I^\Lambda_p$ in
  $|\Lambda|^{1-\alpha}\cdot[-1,1]$ (here, $p$ may depend on $\Lambda$
  and be arbitrarily large) satisfying
  \begin{equation}
    \label{eq:32}
    \inf_{j\not= k}\mathrm{dist}(I_j, I_k)\geq \e^{-|\Lambda|^\delta},
  \end{equation}
  we have, for any sequences of integers
  $k_1=k^\Lambda_1,\cdots,k_p=k^\Lambda_p\in\N^p$,
  \begin{equation}
    \label{eq:110}
    \lim_{|\Lambda|\to+\infty}
    \left|\P\left(\left\{\omega;\
        \begin{aligned}
          &\#\{j;\ \xi_j(\omega,\Lambda)\in
          I_1\}=k_1\\&\vdots\hskip3cm\vdots\\
          &\#\{j;\ \xi_j(\omega,\Lambda)\in I_p\}=k_p
        \end{aligned}
      \right\}\right)-\frac{|I_1|^{k_1}}{k_1!}\e^{-|I_1|}\cdots
    \frac{|I_p|^{k_p}}{k_p!} \e^{-|I_p|}\right|=0.
  \end{equation}
  In particular, $ \Xi(\xi; E_0,\omega,\Lambda) $ defined in
  \eqref{eq:13} converges weakly to a Poisson point process with
  intensity Lebesgue.
\end{Th}
\noindent Note that, in Theorem~\ref{thr:2}, we do not require the
limits
\begin{equation*}
  \begin{split}
    &\lim_{|\Lambda|\to+\infty} \frac{|I_1|^{k_1}}{k_1!} \e^{-|I_1|}=
    \lim_{|\Lambda|\to+\infty}
    \frac{|I^\Lambda_1|^{k^\Lambda_1}}{k^\Lambda_1!}
    \e^{-|I^\Lambda_1|},\    \dots,\\
    &\lim_{|\Lambda|\to+\infty} \frac{|I_p|^{k_p}}{k_p!} \e^{-|I_p|}=
    \lim_{|\Lambda|\to+\infty}
    \frac{|I^\Lambda_p|^{k^\Lambda_p}}{k^\Lambda_p!}
    \e^{-|I^\Lambda_p|}
  \end{split}
\end{equation*}
to exist. \\
Condition~\eqref{eq:116} imposes no restriction upon
condition~\eqref{condrho} if we know that the Minami estimate (M)
holds for all $\rho$ in $(0,1)$. This is the case for all the models
we know of for which the Minami estimate is proved
(see~\cite{MR97d:82046,MR2290333,MR2360226,MR2505733,Kl:11a} and
refereces therein).
\subsubsection{Asymptotic independence of the local processes}
\label{sec:asympt-indep-local}
Once Theorem~\ref{thr:3} is known, it is natural to wonder how the
point processes obtained at distinct energies relate to each other. To
understand this, we assume
\begin{description}
\item[(GM)] for $J\subset K\subset I$, one has
  \begin{equation}
    \label{eq:28}
    \E\left[\text{tr}(\car_J(H_\omega(\Lambda)))
      \cdot\text{tr}(\car_K(H_\omega(\Lambda))-1)\right]\leq
    C |J|\,|K|\,|\Lambda|^2.
  \end{equation}
\item[(D)] for $\beta\in(0,1)$ and $\{E_0,E'_0\}\subset I$
  s.t. $E_0\not=E'_0$, when $L\to+\infty$ and $\ell\asymp L^\beta$,
  one has
  \begin{equation}
    \label{eq:27}
    \P\left(\left\{
        \begin{aligned}
          \sigma(H_\omega(\Lambda_\ell))\cap
          (E_0+L^{-d}[-1,1])\not=\emptyset,\\
          \sigma(H_\omega(\Lambda_\ell))\cap
          (E'_0+L^{-d}[-1,1])\not=\emptyset
        \end{aligned}\right\}\right)= o\left((\ell/L)^d\right).
  \end{equation}
\end{description}
In their nature, assumptions (GM) and (D) are similar: they
state that the probability to have two eigenvalues constrained to some
intervals is much smaller than that of having a single eigenvalue in
an interval. Note that $(\ell/L)^d$ is the order of magnitude of the
right hand side in Wegner's estimate (W) for $H_\omega(\Lambda_\ell)$
and the interval $E_0+L^{-d}[-1,1]$.\\
Assumption (GM) was proved to hold for the discrete Anderson model
in~\cite{MR2505733}. In~\cite{Kl:10}, it is proved that
assumption (D) holds for the discrete Anderson model in dimension 1 at
any two distinct energies, and, in any dimension, for energies
sufficiently far apart from each other.\\
Under these assumptions, we have
\begin{Th}
  \label{thr:4}
  Assume (IAD), (W), (GM), (Loc), and (D) hold.  Pick $E_0\in I$ and
  $E'_0\in I$ such that $E_0\not=E'_0$ and~\eqref{eq:60} is satisfied
  at $E_0$ and $E'_0$.\\
  When $|\Lambda|\to+\infty$, the point processes
  $\Xi(E_0,\omega,\Lambda)$ and $\Xi(E'_0,\omega,\Lambda)$, defined
  in~\eqref{eq:13}, converge weakly respectively to two independent
  Poisson processes on $\R$ with intensity the Lebesgue measure. That
  is, for $U_+\subset\R$ and $U_-\subset\R$ compact intervals and
  $\{k_+,k_-\}\in\N\times\N$, one has
  \begin{equation*}
    \mathbb{P}\left(\left\{\omega;\
        \begin{aligned}
          &\#\{j;\xi_j(E_0,\omega,\Lambda)\in U_+\}=k_+\\
          &\#\{j;\xi_j(E'_0,\omega,\Lambda)\in U_-\}=k_-
        \end{aligned}
      \right\}
    \right)\vers_{\Lambda\to\Z^d} \left(\frac{|U_+|^{k_+}}{k_+!}
      \e^{-|U_+|} \right) 
    \left(   \frac{|U_-|^{k_-}}{k_-!} \e^{-|U_-|} \right).
  \end{equation*}
\end{Th}
\noindent Theorem~\ref{thr:4} naturally leads to wonder how far the
energies $E_0$ and $E'_0$ need to be from each other with respect to
the scaling used to renormalize the eigenvalues
for such a result to still hold.\\
We prove 
\begin{Th}
  \label{thr:8}
  Assume (IAD), (W), (GM), (Loc). Pick $E_0\in I$ such
  that~\eqref{eq:60} is satisfied Assume moreover that the
  density of states $\nu$ is continuous at $E_0$.\\
  Consider two sequences of energies, say $(E_\Lambda)_\Lambda$ and
  $(E'_\Lambda)_\Lambda$ such that
  \begin{enumerate}
  \item  $\displaystyle E_\Lambda\vers_{\Lambda\to\Z^d}E_0$ and
    $\displaystyle E'_\Lambda\vers_{\Lambda\to\Z^d}E_0$,
  \item $\displaystyle |\Lambda|\cdot|N(E_\Lambda)-N(E'_\Lambda)|
    \vers_{\Lambda\to\Z^d}+\infty$.
  \end{enumerate}
  Then, the point processes $\Xi(E_\Lambda,\omega,\Lambda)$ and
  $\Xi(E'_\Lambda,\omega,\Lambda)$, defined in~\eqref{eq:13}, converge
  weakly respectively to two independent Poisson processes on $\R$
  with intensity the Lebesgue measure.
\end{Th}
\noindent Theorem~\ref{thr:8} shows that, in the localized regime,
eigenvalues that are separated by a distance that is asymptotically
infinite when compared to the mean spacing between the eigenlevels,
behave as independent random variables. There are no interactions
except at very short distances.

Assumption (2) can clearly not be omitted in Theorem~\ref{thr:8}; it
suffices to consider e.g. $E_\Lambda, E'_\Lambda$
s.t. $N(E_\Lambda)=N(E'_\Lambda)+ a|\Lambda|^{-1}$ to see that the two
limit random processes are obtained as a shift from one another.
\subsection{The joint (energy - localization center) statistics}
\label{sec:local-cent-stat}
Recall that $E_1(\omega,\Lambda)\leq E_2(\omega,\Lambda)\leq
\cdots\leq E_N(\omega,\Lambda)$ denote the eigenvalues of
$H_{\omega}(\Lambda)$ ordered increasingly and repeated according to
multiplicity. Recall Lemma~\ref{le:3} and Lemma~\ref{le:1}: it states that, to an
eigenvector associated to $E_j(\omega,\Lambda)$, we can associate a
center of localization that we denote by $x_j(\omega,\Lambda)$.
\subsubsection{Uniform Poisson convergence for the joint
  (energy,center)-distribution}
\label{sec:strong-poiss-conv-1}
We now place ourselves in the same setting as in
section~\ref{sec:strong-poiss-conv}. We prove
\begin{Th}
  \label{thr:11}
  Assume (IAD), (W), (M) and (Loc) hold. Let $E_0$ be an energy in $I$
  such that~\eqref{eq:107} holds for some
  $\tilde\rho\in[0,\rho/(1+d\rho))$. Pick
  $\alpha\in(\alpha_{d,\rho,\tilde\rho},1)$.  Then, there exists
  $\delta>0$ such that,
  \begin{itemize}
  \item for any sequences of intervals
    $I_1=I^\Lambda_1,\dots,I_p=I^\Lambda_p$ in
    $|\Lambda|^{1-\alpha}\cdot[-1,1]$ satisfying~\eqref{eq:32},
  \item for any sequences of cubes
    $C_1=C^\Lambda_1,\dots,C_p=C^\Lambda_p$ in $[-1/2,1/2]^d$,
  \end{itemize}
  one has, for any sequences of integers
  $k_1=k^\Lambda_1,\cdots,k_p=k^\Lambda_p\in\N^p$,
  \begin{equation}
    \label{eq:109}
    \lim_{|\Lambda|\to+\infty} \left|\P\left(\left\{\omega;\
          \begin{aligned}
            &\#\left\{j;
              \begin{aligned}
                \xi_j(\omega,\Lambda)&\in I_1 \\
                x_j(\omega,\Lambda)/L&\in C_1 \end{aligned}
            \right\}=k_1\\&\vdots\hskip3cm\vdots\\
            &\#\left\{j;
              \begin{aligned}
                \xi_j(\omega,\Lambda)&\in I_p \\
                x_j(\omega,\Lambda)/L&\in C_p
              \end{aligned}
            \right\}=k_p
          \end{aligned}
        \right\}\right)-
      \prod_{n=1}^pe^{-|I_n||C_n|}\frac{(|I_n||C_n|)^{k_n}}{k_n!}\right|=0,
  \end{equation}
  where the $\xi_j(\omega,\Lambda_L)$'s are defined in \eqref{eq:11}.
  \\
  In particular the point process defined as
  \begin{equation}
    \label{jointprocess}
    \Xi^2_\Lambda(\xi,x;E_0,\Lambda) = \sum_{j=1}^N
    \delta_{\xi_j(\omega,\Lambda_L)}(\xi)
    \otimes\delta_{x_j(\omega,\Lambda)/L}(x)
  \end{equation}
  converges weakly to a Poisson point process on $\R\times \R^{d}$ with
  intensity the Lebesgue measure.
\end{Th}
\noindent The joint (energy,center)-distribution given by $
\Xi^2_\Lambda(\xi,x;E_0,\Lambda)$ in \eqref{jointprocess} have been
studied in~\cite{MR2299191}, where they prove it converges weakly to a
Poisson process. \\
We point out that in Theorem~\ref{thr:11} intervals $I_j$'s and cubes
$C_j$'s may depend on $\Lambda$. But the limit only depends on the
product $|I_j||C_j|$. We shall exploit this fact in the next result.
\subsubsection{Covariant scaling joint (energy,center)-distribution}
\label{sec:covar-scal-joint}
Fix $\xi\in(0,1)$ and an increasing sequence of scales
$\ell=(\ell_{\Lambda})_{\Lambda}$ such that
\begin{equation}
  \label{eq:14}
  \frac{\ell_{\Lambda}}{\log^{1/\xi}|\Lambda|}\vers_{|\Lambda|\to+\infty}
  +\infty \quad\text{ and }\quad \ell_{\Lambda}\leq |\Lambda|^{1/d}.
\end{equation}
Pick $E_0\in I$ so that $\nu(E_0)>0$. Consider the point
process
\begin{equation}
  \label{eq:15}
  \Xi^2_\Lambda(\xi,x;E_0,\ell) = \sum_{j=1}^N
  \delta_{\ell_\Lambda^d[N(E_j(\omega,\Lambda))-N(E_0)]}(\xi)
  \otimes\delta_{x_j(\omega)/\ell_\Lambda}(x). 
\end{equation}
The process is valued in $\R\times\R^d$; actually, if
$c\,\ell_{\Lambda}\geq|\Lambda|^{1/d}$, it is valued in
$\R\times(-c,c)^d$. Assume it exists and define the limit
\begin{equation}
  \label{eq:16}
  c_\ell=\lim_{|\Lambda|\to+\infty}|\Lambda|^{1/d}\ell_\Lambda^{-1}
  \in [1,+\infty].
\end{equation}
Note that if $\ell_\Lambda=L$, we recover~\eqref{jointprocess}. We
prove
\begin{Th}
  \label{thr:5}
  Assume (IAD), (W), (M) and (Loc) hold.  Let $E_0$ be an energy in
  $I$ such that~\eqref{eq:60}  holds for some
  $\tilde\rho\in[0,\rho/(1+d\rho))$. The point process
  $\Xi^2_\Lambda(\xi,x;E_0,\ell)$ converges weakly to a Poisson
  process on $\R\times(-c_\ell,c_\ell)^d$ with intensity measure the
  Lebesgue measure.
\end{Th}
\noindent 
As a result, we see that, once the energies and the localization
centers are scaled covariantly, the convergence to a Poisson process
is true at any scale that is essentially larger than the localization
width. This covariant scaling is very natural; it is the one
prescribed by the Heisenberg uncertainty principle: the more precision
we require in the energy variable, the less we can afford in the space
variable. In this respect, the energies behave like a homogeneous
symbol of degree $d$. This is quite different from what one obtains in
the case of the Laplace operator.
\subsubsection{Non-covariant scaling joint (energy,center)-distribution}
\label{sec:non-covar-scal}
One can also study what happens when the energies and localization
centers are not scaled covariantly. Consider two increasing sequences
of scales, say $\ell=(\ell_{\Lambda})_{\Lambda}$ and
$\tilde\ell=(\tilde\ell_{\Lambda})_{\Lambda}$. Pick $E_0$ an energy in
$I$ such that~\eqref{eq:60} holds for some
$\tilde\rho\in[0,\rho/(1+d\rho))$. Consider the point process
\begin{equation}
  \label{eq:29}
  \Xi^2_\Lambda(\xi,x;E_0,\ell,\tilde\ell) = \sum_{j=1}^N
  \delta_{\ell_\Lambda^d[N(E_j(\omega,\Lambda))-N(E_0)]}(\xi)
  \otimes\delta_{x_j(\omega)/\tilde\ell_\Lambda}(x). 
\end{equation}
Then, we prove
\begin{Th}
  \label{thr:7}
  Assume (IAD), (W), (M) and (Loc) hold.   Let $E_0$ be an energy in
  $I$ such that~\eqref{eq:60}  holds for some
  $\tilde\rho\in[0,\rho/(1+d\rho))$.  Assume the sequences of increasing scales
  $\ell=(\ell_{\Lambda})_{\Lambda}$ and
  $\tilde\ell=(\tilde\ell_{\Lambda})_{\Lambda}$ satisfy~\eqref{eq:14}.
  Assume that
  \begin{equation}
    \label{eq:70}
    \text{if }\ell'=o(L)\text{ then
    }\frac{\ell_{\Lambda_{L+\ell'}}}{\ell_{\Lambda_L}}
    \vers_{|\Lambda|\to+\infty}1\text{ and }
    \frac{\tilde\ell_{\Lambda_{L+\ell'}}}{\tilde\ell_{\Lambda_L}}
    \vers_{|\Lambda|\to+\infty}1.
  \end{equation}
  Let $J$ and $C$ be bounded measurable sets respectively in $\R$ and
  $(-c_{\tilde\ell},c_{\tilde\ell})^d\subset\R^d$. One has
  \begin{enumerate}
  \item if, for some $\delta>0$, one has
    $\displaystyle\frac{\tilde\ell_\Lambda}{\ell_\Lambda}
    \leq|\Lambda|^{-\delta}$, then $\omega$ almost surely, for
    $\Lambda$ sufficiently large, 
    \begin{equation*}
      \int_{J\times C}\Xi^2_\Lambda(\xi,x;E_0,\ell,\tilde\ell)d\xi
      dx=0.
    \end{equation*}
  \item if, for some $\delta>0$, one has
  $\displaystyle\frac{\tilde\ell_\Lambda}{\ell_\Lambda}
  \geq|\Lambda|^{\delta}$, then $\omega$ almost surely,
    \begin{equation*}
      \left(\frac{\ell_\Lambda}{\tilde\ell_\Lambda}\right)^{-d}
      \int_{J\times C}\Xi^2_\Lambda(\xi,x;E_0,\ell,\tilde\ell)d\xi
      dx\vers_{|\Lambda|\to+\infty}|J|\cdot|C|.
    \end{equation*}
  \end{enumerate} 
\end{Th}
\noindent Theorem~\ref{thr:7} proves that the local energy levels and
the localization centers become uniformly distributed in large energy
windows if one conditions the localization centers to a much larger
window. On the other hand, for a typical sample, if one looks for
eigenvalues in an energy interval much smaller than the correctly
scaled one with localization center in a cube, then asymptotically
there are none.\\
If one replaces the polynomial growth or decay conditions of the ratio
of scales $\ell_\Lambda/\tilde\ell_\Lambda$ by the condition that they
tend to $0$ or $\infty$, or if one omits condition~\eqref{eq:70}, the
results stays valid except for the fact that the convergence is not
almost sure anymore but simply holds in probability; actually, one has
convergence in some $L^p$ norm (see Remark~\ref{rem:3} in
section~\ref{sec:proof-theorem7}).
\subsection{Outline of the article}
\label{sec:outline-article}
To complete this section, let us now briefly describe the
architecture of the remaining parts of the paper that consist in the
proofs of all the results stated in sections~\ref{sec:introduction}
and~\ref{sec:results}.\\
We start in section~\ref{sec:asymptotics} with the computation of two
important quantities related to our approximation scheme. Consider a
cube $\Lambda$ and an energy interval $I$ such that
$|I|\cdot|\Lambda|$ is small. In section~\ref{sec:asymptotics}, we
compute
\begin{itemize}
\item the probability that $H_\omega(\Lambda)$ has exactly any eigenvalue
  in $I$,
\item the distribution of this eigenvalue conditioned on the fact
  that it is unique.
\end{itemize}
This distribution is used in the sequel to approximate the eigenvalue
and localization center processes. \\
Section~\ref{sec:cutting-pieces} is devoted to the proof of the
approximation theorems, Theorems~\ref{thr:vsmall1}
and~\ref{thr:vbig1}, in section~\ref{sec:cutting-pieces}
Section~\ref{sec:stat-conv} is devoted to the proof of the results on
the spectral statistics. In section~\ref{sec:proof-prop-theor}, we
derive the results for the full Hamiltonian i.e. we develop the other
point of view presented in
section~\ref{sec:another-point-view}. Finally, the
appendix~\ref{sec:appendix} is devoted to various technical results
used in the course of the proofs, included a description of equivalent
finite volume localization properties.



\section{The local distribution of eigenvalues}
\label{sec:asymptotics}
In this section, we compute the distribution of unfolded eigenvalues.
\subsection{The distribution of unfolded eigenvalues}
\label{sec:distr-renorm-eigenv}
Pick $1\ll\ell'\ll\ell$. Consider a cube $\Lambda$ of side-length
$\ell$ i.e. $\Lambda=\Lambda_{\ell}$ and an interval
$I_\Lambda=[a_\Lambda,b_\Lambda]\subset I$ (i.e. $I_\Lambda$ is
contained in the localization region). Consider the following random
variables:
\begin{itemize}
\item $X=X(\Lambda,I_\Lambda)=X(\Lambda,I_\Lambda,\ell')$ is the
  Bernoulli random variable
\begin{equation*}
  X=\car_{H_\omega(\Lambda)\text{ has exactly one
      eigenvalue in }I_\Lambda\text{ with localization center in }
    \Lambda_{\ell-\ell'}};
\end{equation*}
\item $\tilde E=\tilde E(\Lambda,I_\Lambda)$ is the eigenvalue of
  $H_\omega(\Lambda)$ in $I_\Lambda$ conditioned on $X=1$;
\item $\tilde\xi=\tilde\xi(\Lambda,I_\Lambda)=(\tilde
  E(\Lambda,I_\Lambda)-a_\Lambda)/|I_\Lambda|$.
\end{itemize}
Clearly $\tilde\xi$ is valued in $[0,1]$; let $\tilde\Xi$ be its distribution
function.\\
In the present section, we will describe the distribution of these
random variables as $|\Lambda|\to+\infty$ and $|I_\Lambda|\to0$. We
prove
\begin{Le}
  \label{lemasympt}
  Assume (W), (M) and (Loc) hold. For any $\nu\in(0,1)$, one has
  \begin{equation}
    \label{eq:51}
    \left|\pro(X=1)-N(I_\Lambda)|\Lambda|\right|\lesssim
    (|\Lambda||I_\Lambda|)^{1+\rho}+N(I_\Lambda)|\Lambda|\ell'\ell^{-1}
      +|\Lambda|e^{-(\ell')^\nu}
  \end{equation}
  where $N(E)$ denotes the integrated density of states of
  $H_\omega$.\\
  One has, for all $x,y\in[0,1]$,
  \begin{equation}
    \label{eq:50}
    \left|(\tilde\Xi(x)-\tilde\Xi(y))\,\P(X=1)\right|\lesssim
    |x-y||I_\Lambda||\Lambda|.
  \end{equation}
  Moreover, setting $N(x,y,\Lambda):=[N(a_\Lambda+x|I_\Lambda|)-
  N(a_\Lambda+y|I_\Lambda|)]|\Lambda|$, one has
  \begin{equation}
    \label{eq:52}
    \left|(\tilde\Xi(x)-\tilde\Xi(y))\,\P(X=1)-N(x,y,\Lambda)\right|
    \lesssim (|\Lambda||I_\Lambda|)^{1+\rho}
    +|N(x,y,\Lambda)|\,\ell'\ell^{-1}+ |\Lambda|e^{-(\ell')^\nu}. 
  \end{equation}
\end{Le}
\noindent This lemma differs from the usual computation of the DOS in
the sense that the size of the interval decays as the thermodynamic
limit is taken. A joint limit in the volume and the size of the
interval has to be taken here. The price we pay for this joint limit
is that we shall restrict ourselves to the localization regime, while
the IDS exists in a broader region.\\
First let us note that, when we will use Lemma~\ref{lemasympt} in
conjunction with Theorems~\ref{thr:vsmall1} or~\ref{thr:vbig1}, the
role of $\Lambda$ will be played by the cube $\Lambda_\ell$.\\
Of course, estimates~\eqref{eq:51} and~\eqref{eq:52} are of interest
mainly if their right hand side which is to be understood as an error
is smaller than the main term. In~\eqref{eq:51}, the main restriction
comes from the requirement that $N(I_\Lambda)|\Lambda|\gg
(|\Lambda||I_\Lambda|)^{1+\rho}$ which is essentially a requirement
that $N(I_\Lambda)$ should not be too small with respect to
$|I_\Lambda|$ (similar to that found in Theorems~\ref{thr:vsmall1}
and~\ref{thr:vbig1}). Lemma~\ref{lemasympt} will be used in
conjunction with Theorems~\ref{thr:vsmall1} and~\ref{thr:vbig1}. The
cube $\Lambda$ in Lemma~\ref{lemasympt} will be the cube
$\Lambda_\ell$ in Theorems~\ref{thr:vsmall1} and~\ref{thr:vbig1}.
Therefore, the requirements induced by the other two terms are less
restrictive.\\
In~\eqref{eq:52}, the main restriction comes from the requirement that
$N(x,y,\Lambda)\gg (|\Lambda||I_\Lambda|)^{1+\rho}$. This is
essentially a requirement on the size of $|x-y|$. It should not be too
small. On the other hand, we expect the spacing between the
eigenvalues of $H_\omega(\Lambda_L)$ to be of size $|\Lambda_L|^{-1}$
(we keep the notations of Theorem~\ref{thr:vbig1} and recall that the
cube $\Lambda$ in Lemma~\ref{lemasympt} will be the cube
$\Lambda_\ell$ in Theorem~\ref{thr:vbig1}). So to distinguish between
the eigenvalues, one needs to be able to know $\tilde\Xi$ up to
resolution $|x-y||I_\Lambda|\asymp|\Lambda_L|^{-1}$. This will force us to
use Lemma~\ref{lemasympt} on intervals $I_\Lambda$ such that
$N(I_\Lambda)\asymp|\Lambda|^{-\alpha}$ for some $\alpha\in(0,1)$ close
to $1$ (see the discussion following Theorem~\ref{thr:vbig1} and
section~\ref{sec:some-prel-cons}).\\
To prove Lemma~\ref{lemasympt}, we will use
\begin{Le}
  \label{le:9}
  Assume (W) and (Loc) hold in $I$ a compact interval. For
  $\nu\in(0,1)$ and $1\ll\ell'\ll\ell$, let $N(J,\ell,\ell')$ be the
  number of eigenvalues of $H_\omega(\Lambda_\ell)$ in $J$ with
  localization center in $\Lambda_{\ell-\ell'}$. Then, there exists
  $C>0$ such that, for $J\subset I$ an interval such that $|J|\geq
  e^{-(\ell')^\nu}$, one has
  \begin{equation}
    \label{eq:96}
    \left|\esp(N(J,\ell,\ell'))-N(J)|\Lambda_\ell|\right|\lesssim
    N(J)|\Lambda_\ell|\ell'\ell^{-1} +\ell^d e^{-(\ell')^\nu}
  \end{equation}
\end{Le}
\noindent We turn to the proofs of Lemma~\ref{le:9} and~\ref{lemasympt}
\begin{proof}[Proof of Lemma~\ref{le:9}]
  Recall Lemma~\ref{le:3} and let $\mathcal{V}_{\Lambda_\ell}$ be the set
  of configurations given in Part~(II) of Lemma~\ref{le:3} for some
  $\nu\in(0,1)$ given. Outside $\mathcal{V}_{\Lambda_\ell}$, we bound the
  number of eigenvalues of $H_\omega(\Lambda)$ by
  $C|\Lambda_\ell|$. Thus, one has
  \begin{equation}
    \label{eq:97}
    \esp(\car_{\omega\not\in\mathcal{V}_{\Lambda_\ell}}N(J,\ell,\ell'))\lesssim \ell^d
    e^{-\ell^\nu}.
  \end{equation}
  Assume now that $\omega\in\mathcal{V}_{\Lambda_\ell}$. It follows from
  Lemma~\ref{lemcenter} that for such $\omega$'s, one has
  \begin{equation}
    \label{eq:99}
    \tr(\car_{\Lambda_{\ell-2\ell'}}\car_{J_-}(H_\omega))
    +O\left(\ell^d e^{-(\ell')^\nu}\right) 
    \leq N(J,\ell,\ell'))\leq
    \tr(\car_{\Lambda_\ell}\car_{J_+}(H_\omega))    +O\left(\ell^d
      e^{-(\ell')^\nu}\right) 
  \end{equation}
  where, for some $C>0$, $J_+=J+Ce^{-(\ell')^\nu}[-1,1]$ and
  $J_-=J\setminus[(\R\setminus J)+C e^{-(\ell')^\nu}[-1,1]]$. Note
  that $|J|-2Ce^{-(\ell')^\nu}\leq |J_-|\leq |J_+|\leq
  |J|+2Ce^{-(\ell')^\nu}$.\\
  Taking the expectation of the right hand side of~\eqref{eq:99},
  using the covariance for the operator $H_\omega$ and the Wegner
  estimate (W), we compute
  \begin{equation}
    \label{eq:100}
    \esp(\tr(\car_{\Lambda_{\ell}}\car_{J_+}(H_\omega)))
    =N(J_+)|\Lambda_{\ell}|=N(J)|\Lambda_{\ell}|
    +O\left(\ell^d e^{-(\ell')^\nu}\right).
  \end{equation}
  The left hand side is estimated in the same way. Plugging this back
  into the expectation of~\eqref{eq:99} and using~\eqref{eq:97}, $|\Lambda_{\ell-2\ell'}|= |\Lambda_{\ell}|(1+C\ell'\ell^{-1})$, and
  the assumption that $|J|\geq e^{-(\ell')^\nu}$ easily
  yields~\eqref{eq:96}. This completes the proof of Lemma~\ref{le:9}.
\end{proof}
\begin{proof}[Proof of Lemma~\ref{lemasympt}]
  Using the notations of Lemma~\ref{le:9}, note that
  $\P(X=1)=\P\set{N(I_\Lambda,\ell,\ell')=1}$. First, we relate
  $\P\set{N(I_\Lambda,\ell,\ell')=1}$ to $\E[N(J,\ell,\ell')]$. To do
  so, we follow the ideas used~\cite{MR2413200,MR2505733} to estimate
  the probability for $H_\omega(\Lambda)$ to have an eigenvalue in
  $J$. We notice that, as $N(J,\ell,\ell')$ is integer valued
  \begin{equation*}
    \begin{split}
      \E[N(J,\ell,\ell')]-\P(X=1)&=
      \E[N(J,\ell,\ell')]-\P\set{N(I_\Lambda,\ell,\ell')=1}\\&=
      \sum_{k=2}^\infty k\,\P\{N(I_\Lambda,\ell,\ell')=k\}\\ &\leq
      \sum_{k=2}^\infty k(k-1)\P\{ \tr
      \car_{I_\Lambda}(H_\omega(\Lambda)) =k\} \\ & = \E \set{\tr
        \car_{I_\Lambda}(H_\omega(\Lambda))\,\pa{ \tr
          \car_{I_\Lambda}(H_\omega(\Lambda))-1}}.
    \end{split}
  \end{equation*}
  Thus, by our assumption (M), we know
  \begin{equation}
    \label{eq:114}
    0\leq \E[N(J,\ell,\ell')]-\P(X=1)
    \leq C |\Lambda|^{1+\rho}|I_\Lambda|^{1+\rho} 
  \end{equation}
  The evaluation of $\E[N(J,\ell,\ell')]$ is then given by
  Lemma~\ref{le:9}. This yields~\eqref{eq:51}.\\
  The estimate~\eqref{eq:50} is an immediate consequence of the Wegner
  estimate (W) and the normalization of $\tilde\Xi$. \\
  Set
  $I_{x,y,\Lambda}=[a_\Lambda+x|I_\Lambda|,a_\Lambda+y|I_\Lambda|]$. To
  prove~\eqref{eq:52}, we write
  \begin{equation}
    \label{eq:55}
    \begin{split}
      \left|\P\{N(I_{x,y,\Lambda},\ell,\ell')=1\}
        -(\tilde\Xi(x)-\tilde\Xi(y))\P(X=1)\right|&\leq
      \P\left\{
        \begin{aligned}
          H_\omega(\Lambda)\text{ has at least}\\\text{ 2 eigenvalues in
          }I_\Lambda
        \end{aligned}
      \right\}\\
      &\lesssim (|\Lambda_\ell||I_\Lambda|)^{1+\rho}
    \end{split}
  \end{equation}
  using (M).\\
  Replacing $I_\Lambda$ with the interval $I_{x,y,\Lambda}$ in the
  estimation of $\P(N(I_\Lambda,\ell,\ell')=1)$ yields the estimation
  of the probability $\P(N(I_{x,y,\Lambda},\ell,\ell')=1)$ and, thus,
  completes the proof of~\eqref{eq:52}.\\
  The proof of Lemma~\ref{lemasympt} is complete.
\end{proof}
\begin{Rem}
  \label{rem:10}
  The gist of Lemma~\ref{lemasympt} is that the local distribution of
  $\tilde E$ is that of the density of states i.e., in~\eqref{eq:52},
  the remainder terms should be negligible with respect to
  $N(x,y,\Lambda)$. Clearly, this will only be the case if
  $N(x,y,\Lambda)\gg (|\Lambda||I_\Lambda|)^{1+\rho}$. This imposes a
  condition on the size of $|y-x|$ i.e. $y-x$ can't be too small. This
  restriction will be made clear in the following result.\\
  If one uses the improved Minami estimates of~\cite{MR2505733}
  in~\eqref{eq:55}, one can improve 
  \begin{multline*}
    \left|(\tilde\Xi(x)-\tilde\Xi(y))\,\P(X=1)-N(x,y,\Lambda)\right|\\
    \leq C\left(|x-y|^{1+\rho}|\Lambda|^{1+\rho}|I_\Lambda|^{1+\rho}+
      |\Lambda|^2|I_\Lambda||I_{x,y,\Lambda}|+N(x,y,\Lambda)\ell'\ell^{-1}
      +|\Lambda|e^{-(\ell')^\nu}\right)
  \end{multline*}
  and thus take advantage of the possible smallness of
  $I_{x,y,\Lambda}$ compared to $I_\Lambda$. So this will lift the
  above restriction, at least if for $J\subset I_\Lambda$, one has
  $N(J)\gtrsim |J|$. This can be done in some
  cases~\cite{Ge-Kl:10a}.
\end{Rem}
\noindent We now describe the distribution of the unfolded
eigenvalues. Therefore we slightly change our notations to localize the
quantities near some energy $E_0$. Let $1\ll\ell'\ll\ell$. Pick
$E_0\in I$ such that~\eqref{eq:60} be satisfied and
$I_\Lambda=[a_\Lambda,b_\Lambda]$. Recall that
\begin{itemize}
\item $X=X(\Lambda,E_0+I_\Lambda)=X(\Lambda,E_0+I_\Lambda,\ell')$ is
  the Bernoulli random variable
  \begin{equation*}
    X=\car_{H_\omega(\Lambda)\text{ has exactly one
        eigenvalue in }E_0+I_\Lambda\text{ with localization center in }
      \Lambda_{\ell-\ell'}}
  \end{equation*}
\item $\tilde E=\tilde E(\Lambda,E_0+I_\Lambda)$ is this eigenvalue
  conditioned on $X=1$.
\end{itemize}
Define
\begin{equation}
  \label{eq:88}
  \xi=\frac{N(\tilde
    E)-N(E_0+a_\Lambda)}{N(E_0+b_\Lambda)-N(E_0+a_\Lambda)}
  =\frac{N(\tilde E)-N(E_0+a_\Lambda)}{N(E_0+I_\Lambda)}.
\end{equation}
The random variable $\xi$ is valued in $[0,1]$. Let $\Xi$ be its
distribution function. We prove
\begin{Le}
  \label{le:7}
  Assume (W), (M) and (Loc) hold. Pick $E_0$ such that~\eqref{eq:86}
  holds for $\rho'\in(0,\rho)$. Fix $\nu\in(0,1)$. Assume, moreover, that
  \begin{equation}
    \label{eq:91}
    e^{-(\ell')^\nu}\leq N(E_0+I_\Lambda)
    =o\left(|\Lambda|^{-\rho(1+\rho')/(\rho-\rho')}\right)
    \text{ when }|\Lambda|\to+\infty.
  \end{equation}
  Then, for $1\ll\ell'\ll\ell$ and $(x,y)\in[0,1]^2$ such that
  $\displaystyle|x-y|\gg
  N\left(I_\Lambda\right)^{(\rho-\rho')/(1+\rho')}|\Lambda|^\rho$, one
  has
  \begin{equation}
    \label{eq:89}
    \Xi(x)-\Xi(y)=(x-y)\left(1+O\left(\ell'\ell^{-1}+ e^{-(\ell')^{\nu}}+
        |x-y|^{-1}N(E_0+I_\Lambda)^{\frac{\rho-\rho'}{1+\rho'}}\,
        |\Lambda|^\rho\right)\right).
  \end{equation}
\end{Le}
\noindent Recalling the discussion following Lemma~\ref{lemasympt}, to
be able to perform our analysis of the level spacings, we will
need~\eqref{eq:89} to give a good approximation of $\Xi(x)-\Xi(y)$ for
$|x-y|\ll(N(E_0+I_\Lambda)|\Lambda_L|)^{-1}$ (recall that $\Lambda$ in
Lemma~\ref{le:7} is $\Lambda_\ell$ in
Theorem~\ref{thr:vbig1}). Indeed, by Lemma~\ref{lemasympt} and
Theorem~\ref{thr:vbig1}, the number of eigenvalues of
$H_\omega(\Lambda_L)$ in $I_\Lambda$ is asymptotic to
$N(E_0+I_\Lambda)|\Lambda_L|$ (see Theorem~\ref{thr:16}).
\begin{proof} Recall that the IDS $N$ is monotone by definition and
  Lipschitz continuous thanks to (W). Assumption~\eqref{eq:86} and the
  Wegner estimate (W) guarantee that, for $\Lambda$ sufficiently
  large, for $[a,b]\subset E_0+I_\Lambda$, one has
  \begin{equation}
    \label{eq:92}
    \frac1C(b-a)\leq |N^{-1}([a,b])|\leq(b-a)^{1/(1+\rho')}.
  \end{equation}
  Here, $|N^{-1}([a,b])|$ denotes the Lebesgue measure of the interval
  $N^{-1}([a,b])$.\\
  By the definitions of $\tilde\xi$ and $\xi$ (see the beginning of
  section~\ref{sec:distr-renorm-eigenv} and~\eqref{eq:88}), for
  $x\in[0,1]$, one has
  \begin{equation}
    \Xi(x)=\tilde\Xi
    \left[N^{-1}(N(E_0+a_\Lambda)+x N(E_0+|I_\Lambda|))\right].
  \end{equation}
  By Lemma~\ref{lemasympt} applied e.g. for $\nu$ replaced with
  $(1+\nu)/2$, for $(x,y)$ as in~\eqref{eq:89}, one has
  \begin{equation}
    \Xi(x)-\Xi(y)=(x-y)\frac{1+A+B+C}{1+A'+B'+C'}
  \end{equation}
  where, using~\eqref{eq:92}, the assumption on $|x-y|$
  in~\eqref{eq:89} and the left hand side of~\eqref{eq:91}, we compute
  \begin{align*}
      |A|&\lesssim
      \frac{(|I_\Lambda||\Lambda|)^{1+\rho}}{|x-y|N(I_\Lambda)|\Lambda|}
      \lesssim\frac{N(E_0+I_\Lambda)^{(\rho-\rho')/(1+\rho')}|\Lambda|^\rho}
      {|x-y|}
    ,\quad\quad\quad\quad|B|\lesssim \ell'\ell^{-1},
    \\
     |C|&\lesssim
      \frac{e^{-(\ell')^{(1+\nu)/2}}}{N(E_0+I_\Lambda)|x-y|}\lesssim
      e^{-(\ell')^{\nu}}.
  \end{align*}
  The quantities $|A'|$, $|B'|$ and $|C'|$ are respectively bounded by
  the same bounds as $|A|$ , $|B|$ and $|C|$ for $(x,y)=(0,1)$. This
  completes the proof of Lemma~\ref{le:7}.
\end{proof}
\subsection{The proof of Theorem~\ref{thr:16}}
\label{sec:local-count-funct}
Theorem~\ref{thr:16} will be a consequence of Theorem~\ref{thr:vbig1}
and Lemma~\ref{lemasympt}. We use the notation of
Theorem~\ref{thr:vbig1}. Recall that the number of eigenvalues of
$H_\omega(\Lambda) $ in $I_\Lambda$ is denoted by
$N(I_\Lambda,\Lambda,\omega)$. The control of
$N(I_\Lambda,\Lambda,\omega)$ will be useful in order obtain the level
spacings statistics.\\
Recall the assumptions of Theorem~\ref{thr:16}:
\begin{itemize}
\item $N(I_\Lambda)\,(\log|\Lambda|)^{1/\delta}\to0$ as
  $|\Lambda|\to+\infty$
\item $N(I_\Lambda)\,|\Lambda|^{1-\nu}\to+\infty$ as
  $|\Lambda|\to+\infty$
\item $N(I_\Lambda)\,|I_\Lambda|^{-1-\tilde\rho}\to+\infty$ as
  $|\Lambda|\to+\infty$,
\end{itemize}
For $\delta>0$ sufficiently small, this guarantees that one can pick
$\alpha>1$ large and $\ell=\ell_\Lambda\asymp(\log|\Lambda|)^\alpha$
and $\ell'=\ell'_\Lambda\asymp(\log|\Lambda|)^{1/\xi}$ so that they
fulfill all the assumptions of Theorem~\ref{thr:vbig1}, in
particular,~\eqref{eq:93},~\eqref{eq:102} and~\eqref{condell} for some
$\xi\in(0,1)$ (see also the discussion following
Theorem~\ref{thr:vbig1}). The estimate~\eqref{eq:94} gives the
probability of $\mathcal{Z}_\Lambda$, the set of configurations where
one has a good description of most the eigenvalues. Moreover, for
$\delta>0$ sufficiently small,, the number of eigenvalues of
$H_\omega(\Lambda)$ in $I_\Lambda$ that are not described by
Theorem~\ref{thr:vbig1} is bounded by
\begin{equation}
  \label{eq:53}
  C N(I_\Lambda)|\Lambda|
  \left(N(I_\Lambda)^{\frac{\rho-\tilde{\rho}}{1+\tilde{\rho}} }
    \ell_\Lambda^{d(1+\rho)}+(\ell'_\Lambda)^{d+1}\ell_\Lambda^{-1}\right)\leq
  CN(I_\Lambda)|\Lambda|(\log|\Lambda|)^{-\delta}.
\end{equation}
Consider the boxes $(\Lambda_\ell(\gamma_j))_{1\leq j\leq \tilde N}$
given by Theorem~\ref{thr:vbig1}. Their number, say, $\tilde
N$ satisfies $\tilde{N}=|\Lambda|/|\Lambda_\ell|\,(1+o(1))$.\\
For $1\leq j\leq \tilde N$, let $X_j=X(\Lambda_\ell(\gamma_j),
I_\Lambda)$ i.e. $X_j$ is the Bernoulli random variable equal to 1 if
$H_\omega(\Lambda_\ell(\gamma_j))$ has exactly one eigenvalue in
$I_\Lambda$ with localization center at a distance at least $\ell'$
from $\partial\Lambda$ (see Theorem~\ref{thr:vbig1}) and zero
otherwise. It follows from Lemma~\ref{lemasympt} and the choice for
$(\ell,\ell')$ in Theorem~\ref{thr:vbig1} made above that 
\begin{equation}
  \label{eq:101}
  \P(X_j=1)= N(I_\Lambda)\,|\Lambda_\ell|\,
  \left[1+o\left((\log|\Lambda|)^{-\delta})\right)\right].  
\end{equation}
We have
\begin{equation}
  \label{eq:54}
  \left|N(I_\Lambda,\Lambda,\omega) - N(I_\Lambda)|\Lambda|
  \right| \leq \left| N(I_\Lambda,\Lambda,\omega) -
    \sum_{j=1}^{\tilde{N}} X_j \right| + \left|
    \sum_{j=1}^{\tilde{N}} X_j - N(I_\Lambda)|\Lambda| \right|
\end{equation}
By~\eqref{eq:53}, for $\omega\in\mathcal{Z}_\Lambda$, we have
\begin{equation*}
  \left| N(I_\Lambda,\Lambda,\omega) - \sum_{j} X_j \right|
  \lesssim N(I_\Lambda)|\Lambda|(\log|\Lambda|)^{-\delta}.
\end{equation*}
The second term in the right hand side of~\eqref{eq:54} is then
bounded by a standard large deviation estimate for i.i.d. Bernoulli
variables valued in $\{0,1\}$ with expectation $p=p(\tilde N)$
s.t. $p\in(0,1/2]$ and $p\,\tilde N\sim N(I_\Lambda)|\Lambda|
\to+\infty$ (see e.g.~\cite{MR98m:60001}); for $\delta'\in(0,1/2)$, it
yields, for $|\Lambda|$ sufficiently large,
\begin{align}
  \P\left( \left|\sum_{j=1}^{\tilde{N}} X_j - p\tilde{N}\right| \ge
    (p\tilde{N})^{\delta'} \right) \leq e^{-(p\tilde{N})^{2
      \delta'-1}/4}.
\end{align}
Theorem~\ref{thr:16} follows taking $\delta'$ close to $1/2$,
using~\eqref{eq:101} and noting that $p\tilde{N}\gg\log|\Lambda|$ by
our assumptions on $N(I_\Lambda)$.\qed
\begin{Rem}
  \label{rem:5}
  We can get a more precise version of Theorem~\ref{thr:16} by
  optimizing in the intermediate scale $\ell$, the number of
  eigenvalues we miss in the picture of Theorem~\ref{thr:vbig1},
  namely by choosing $\ell$ so that $Kn=K'n'$. Estimates can even be
  improved by resorting to higher order Minami estimates in order to
  bound the missing eigenvalues (replacing the crude deterministic
  bound given by the Weyl formula).
\end{Rem}
\begin{Rem}
  \label{rem:4}
  If $N(I)\ll |I|^{1+\rho}$, Theorem~\ref{thr:16} still holds if one
  can improve on the Minami estimate, replacing one power of the
  interval length $|I|$ by $N(I)$. This can be done in some
  cases~\cite{Ge-Kl:10a}.
\end{Rem}



\section{The proofs of Theorems~\ref{thr:vbig1} and~\ref{thr:vsmall1}}
\label{sec:cutting-pieces}
\noindent Recall $\Lambda_L$ is a cube of side length $L$.  Let
$I_\Lambda$ be an interval inside $I$ the region of localization.\\
We first prove the useful
\begin{Le}
  \label{lemcenter}
  Assume (IAD), (W), (M) and (Loc). Consider scales $\ell',\ell$ so
  that $(\log |\Lambda_L|)^{1/\xi} \ll \ell'\ll \ell \ll L$, and, for
  some given $\gamma\in\Lambda_L$, consider a box
  $\Lambda_\ell(\gamma)$ such that
  $\Lambda_{\ell-\ell'}(\gamma)\subset\Lambda_L$. Let
  $\mathcal{W}_{\Lambda_L}$ be either the set
  $\mathcal{U}_{\Lambda_L}$ or $\mathcal{V}_{\Lambda_L}$ defined in
  Lemma~\ref{le:3}. For $L$ large enough, we have:
  \begin{enumerate}
  \item for any $\omega\in\mathcal{W}_{\Lambda_L}$, if $E(\omega)$ is
    an eigenvalue of $H_\omega(\Lambda_L)$ with a localized
    eigenfunction in the sense of \eqref{eq:19} with a center of
    localization in $\Lambda_{\ell-\ell'}(\gamma)$, then
    $H_\omega(\Lambda_L\cap\Lambda_{\ell}(\gamma))$ has an eigenvalue
    in a neighborhood of $E(\omega)$ of size
    $\mathcal{O}(\e^{-(\ell')^\xi/2})$; moreover, if
    $\omega\in\mathcal{W}_{\Lambda_\ell(\gamma)}$, the corresponding
    eigenfunction is localized in the sense of \eqref{eq:19}.
  \item Assume now additionally that
    $\Lambda_\ell(\gamma)\subset\Lambda_L$. Then, conversely, for any
    $\omega\in\mathcal{W}_{\Lambda_\ell(\gamma)}$, if $E(\omega)$ is
    eigenvalue of $H_\omega(\Lambda_L)$ in
    $H_\omega(\Lambda_{\ell}(\gamma))$ with an eigenfunction
    exponentially localized, in the sense of \eqref{eq:19}with a
    center of localization in $\Lambda_{\ell-\ell'}(\gamma)$, then
    $H_\omega(\Lambda_L)$ has an eigenvalue in a neighborhood of
    $E(\omega)$ of size $\mathcal{O}(\e^{-(\ell')^\xi/2})$; moreover,
    if $\omega\in\mathcal{W}_{\Lambda_L}$, the corresponding
    eigenfunction is localized in the sense of \eqref{eq:19}.
  \end{enumerate}
  As a consequence of (1), (W) and (M), given an interval $I_\Lambda$,
  \begin{itemize}
  \item the probability to have at least one center of localization in
    $\Lambda_{\ell-\ell'}(\gamma)$ corresponding to an eigenvalue of
    $H_\omega(\Lambda_L)$ in $I_\Lambda$ is bounded by $
    C\left(\P(\mathcal{W}_{\Lambda_L}^c) + |I_\Lambda| \ell^d+\ell^d
      e^{-(\ell')^\xi/2}\right)$;
  \item the probability to have at least two centers of localization
    in $ \Lambda_{\ell-\ell'}(\gamma)$ corresponding to two
    eigenvalues of $H_\omega(\Lambda_L)$ in $I_\Lambda$ is bounded by
    $\displaystyle C\left(\P(\mathcal{W}_{\Lambda_L}^c)
      +(|I_\Lambda|\ell^d)^{(1+\rho)}+L^{d(1+\rho)}
      e^{-\rho(\ell')^\xi/2}\right)$.
  \end{itemize}
\end{Le}
\noindent Similar results can be found in~\cite{Kl:10}.
\begin{Rem}
  \label{rem:12}
  In the first part of Lemma~\ref{lemcenter}, we do not require the
  small cube $\Lambda_\ell(\gamma)$ to lie entirely inside the big
  cube $\Lambda_L$. This will be used in our analysis to treat the
  localization centers near the boundary of $\Lambda_L$.\\
  If one has $\Lambda_\ell(\gamma)\subset\Lambda_L$, then, using
  Lemma~\ref{lemasympt}, the bound on the probability to have at least
  one center of localization in $\Lambda_{\ell-\ell'}(\gamma)$
  corresponding to an eigenvalue of $H_\omega(\Lambda_L)$ in
  $I_\Lambda$ can be improved into
  $C\left(\P(\mathcal{W}_{\Lambda_L}^c) + N(I_\Lambda) \ell^d+\ell^d
    e^{-(\ell')^\xi/2}\right)$.\\
  Finally, we note that, in the last part of Lemma~\ref{lemcenter}, it
  is of importance that the probability that appears, namely
  $\P(\mathcal{W}_{\Lambda_L}^c)$, is the one related to the box
  $\Lambda_L$, and not to a small box of size $\ell$.
\end{Rem}
\begin{proof}
  We start with the point (1). Let $\varphi=\varphi_{\omega,\Lambda_L}$
  be the eigenfunction associated to the center $x(\omega)$, and
  $E(\omega)\in I_\Lambda$ the corresponding eigenvalue. Let
  $\Psi_{\ell}$ be a smooth characteristic function covering the cube
  $\Lambda_{\ell}(\gamma)$, i.e. $\Psi_\ell$ such that $\supp
  \Psi_{\ell}\subset \Lambda_{\ell}(\gamma)$, $\Psi_{\ell}=1$ on
  $\Lambda_{\ell-\frac12\ell'}(\gamma)$,
  $\supp\nabla\Psi_{\ell}\subset \Lambda_{\ell}(\gamma)\setminus
  \Lambda_{\ell-\frac12 \ell'}(\gamma)$.\\
  Since $\omega\in\mathcal{U}_{\Lambda_L}$, we have $\|\Psi_\ell
  \varphi\|\ge \frac12$ for $\ell'$ large enough. Set
  $\eta_\ell:=\Psi_\ell\varphi/\|\Psi_{\ell}\varphi\|$. Then,
  $\eta_\ell$ is an approximate eigenvector of the Hamiltonian
  $H_{\omega}(\Lambda_{\ell}(\gamma))$, in the sense that $
  \|\eta_\ell\|=1$ and
  \begin{equation*}
    \| (H_{\omega}(\Lambda_{\ell}(\gamma))-E ) \eta_\ell \|
    \le 2  \| [H_{\omega}(\Lambda_{\ell}(\gamma)) ,\Psi_\ell]
    \varphi \|
    \lesssim \sup_{\supp\nabla\Psi_{\ell}} |\varphi| \lesssim
    \e^{-(\ell')^\xi}.
  \end{equation*}
  It follows that $H_{\omega}(\Lambda_{\ell}(\gamma))$ has an
  eigenvalue in the interval $[E-c\e^{-(\ell')^\xi},
  E+c\e^{-(\ell')^\xi}]$.\\
  Let us prove point (2). Recall $\Psi_\ell$ above. Let
  $\varphi=\varphi_{\omega,\Lambda_\ell(\gamma)}\in
  \mathrm{L}^2(\Lambda_\ell(\gamma))$ be the eigenfunction associated
  to $x(\omega)$. Set $\eta_\ell=\Psi_\ell \varphi$ on
  $\Lambda_\ell(\gamma)$ and $\eta_\ell=0$ on
  $\Lambda_L\setminus\Lambda_\ell(\gamma)$. Since
  $\dist(x(\omega)), \partial\Lambda_\ell(\gamma))\ge \ell'$, it is
  immediate to see that $\eta_\ell$ is an approximate eigenfunction of
  $H_{\omega}(\Lambda_L)$ in the sense that
  \begin{equation*}
    \|(H_{\omega}(\Lambda_L) -E)\eta_\ell \|
    \lesssim \sup_{\supp\nabla\Psi_{\ell}} |\varphi|
    \lesssim \e^{-(\ell')^\xi/2}.
  \end{equation*}    
  The first consequence is immediate. For the second, if the two
  eigenvalues of $H_\omega(\Lambda)$ are at a a distance at least
  $e^{-(\ell')^\xi/2}$ from each other, by point (2), they give rise to
  two distinct eigenvalues of $H_\omega(\Lambda_\ell(\gamma))$; thus,
  we can apply (M). If they are closer, we bound the probability using
  (M) for $H_\omega(\Lambda_L)$.\\
  This completes the proof of Lemma~\ref{lemcenter}.
\end{proof}
We now turn to the proofs of Theorems~\ref{thr:vsmall1}
and~\ref{thr:vbig1}.
\begin{proof}[Proof of Theorem~\ref{thr:vsmall1}]
  For $\beta'>0$ sufficiently small (this will be defined precisely
  below) and $\beta>\beta'$ to be chosen later, set $\ell'=L^{\beta'}$
  and $\ell$ so that $(\ell+\ell')k + \ell'= L$, where
  $k=[L^{1-\beta}]$. Note that $\ell=\mathcal{O}(L^\beta)$ in the
  large volume limit. With such definitions we can pick equally
  distributed boxes of size $\ell$ in $\Lambda_L$ (with distance
  $\ell'$ between two neighbors) satisfying the conditions of the
  theorem.\\
  Note that for $L$ large enough $\ell'>R$, so that events based on
  distinct boxes $\Lambda_\ell(\gamma_j)$'s are independent.\\
  In this proof, we shall use the localization property described by
  Lemma~\ref{le:3} Part~(I).
  \\
  Up to a probability less than $C|\Lambda_L|^{(1-\beta) - \beta p}
  \lesssim |\Lambda_L|^{-1}$ provided $p\ge 2\beta^{-1}-1$, we can
  assume that all the boxes $\Lambda_\ell(\gamma)$ satisfy Part~(I) of
  Lemma~\ref{le:3}, since $p$ in Lemma~\ref{le:3} can be chosen
  arbitrary large. Since $\alpha-\alpha_{d,\rho,\tilde\rho} <1$ we can
  neglect this probability.
  \\
  Recall $N(I_\Lambda) \asymp |\Lambda_L|^{-\alpha}$. Since
  $N(I_\Lambda)\lesssim |I_\Lambda|$, we have $
  \P(\mathcal{U}_{\Lambda_L}^c) \lesssim |I_\Lambda| \ell^d $ and
  $\P(\mathcal{U}_{\Lambda_L}^c) \lesssim |I_\Lambda|^{1+\rho}
  \ell^{2d}$, provided $p>0$  in Lemma~\ref{le:3} is large enough.\\
  Let $\mathcal{S}_{\ell,L}$ be the set of boxes
  $\Lambda_{\ell-\ell'}(\gamma_j)\subset\Lambda_L$ containing at least
  two centers of localization of $H_\omega(\Lambda_L)$. It follows from
  Lemma~\ref{lemcenter} and \eqref{eq:60} that,
  \begin{equation}
    \label{proba1}
    \P(\# \mathcal{S}_{\ell,L} \ge 1) \lesssim
    |\Lambda_L|^{1-\beta}  (|\Lambda_L|^{\beta}|I_\Lambda|)^{1+\rho} \le
    |\Lambda_L|^{1+\beta\rho} N(I_\Lambda)^{\frac{1+\rho}{1+\tilde{\rho}}}
    \lesssim    |\Lambda_L|^{1+\beta\rho -
      \alpha\frac{1+\rho}{1+\tilde{\rho}}}.
  \end{equation}
  To insure that all the centers of $H_\omega(\Lambda_L)$ fall inside
  one of the $\Lambda_\ell(\gamma_j)$'s and actually sufficiently well
  inside (by a distance $\ell'$), we define $\Upsilon\subset\Lambda_L$
  as the set $\Lambda_L\setminus\cup_j\Lambda_\ell(\gamma_j)$ enlarged
  by a length $\ell'$. In other terms $\Upsilon=\Lambda_L\setminus
  \cup_j\Lambda_{\ell-\ell'}(\gamma_j)$.  One has $|\Upsilon| \lesssim
  |\Lambda_L| \ell' / \ell$. We consider a partition
  $\displaystyle\Upsilon = \bigcup_{m=1}^{2^{d-1}} \Upsilon_m$, with
  $\Upsilon_m\cap\Upsilon_{m'}=\emptyset$ if $m\neq m'$, each
  $\Upsilon_m$ being a union of boxes of side length $\ell'$ which are
  two by two distant by at least $\ell'$.  For each given $m$, the
  distance between two boxes of $\Upsilon_m$ is larger than $\ell'$,
  so that we can enlarge each box in $\Upsilon_m$ by, say, $\frac1
  {10} \ell'$, except for sides of boxes that coincide with the
  boundary of $\Lambda_L$. It follows from Lemma~\ref{lemcenter} and
  the Wegner estimate that,
  \begin{align}
    \label{proba2}
    & \P(H_\omega(\Lambda_L) \mbox{ has a center of localization in }
    \Upsilon) \\
    & \lesssim \sum_m \P(H_\omega(\Lambda_L) \mbox{ has a center of
      localization in }  \Upsilon_m) \\
    & \lesssim \sum_m |\Upsilon_m| |I_\Lambda| \lesssim |\Upsilon|
    N(I_\Lambda)^{\frac 1{1+\tilde{\rho}}} \lesssim |\Lambda_L|^{1 -
      \frac1d(\beta-\beta') - \frac \alpha{1+\tilde{\rho}}}
  \end{align}
  We, thus, require that
  \begin{equation*}
    \alpha > (1+\tilde{\rho}) \max \left(\frac{1+\beta\rho}{1+\rho}, 1 -
      \frac1d(\beta-\beta') \right). 
  \end{equation*}
  Optimization yields
  \begin{equation}
    \label{alphaopt} 
    \beta = \frac{d\rho + \beta'(1+\rho)}{(d+1)\rho+1} \, ,
    \quad \alpha >  (1+\tilde{\rho}) \left( 1 - \frac{d\rho +
        \beta'(1+\rho)}{(d+1)\rho+1}\right).  
  \end{equation}
  We, thus, require $\alpha_{d,\rho,\tilde\rho}<\alpha$ where
  \begin{equation*}
    \alpha_{d,\rho,\tilde\rho} =  (1+\tilde{\rho}) \left(1 -
      \frac{\rho}{\rho (d+ 1)+1} \right)= (1+\tilde{\rho})
    \frac{\rho d+1}{\rho (d +1)+1} < 1,
  \end{equation*}
  which is our assumption~\eqref{eq:59}.\\
  It follows from \eqref{proba1} and \eqref{proba2} that with
  probability larger than $1- c|\Lambda_L|^{-(\alpha-
    \alpha_{d,\rho,\tilde\rho} )}$, item~(2) of
  Theorem~\ref{thr:vsmall1} holds, as well as the ``only
  if" part of item~(3).\\
  Next, item~(2) of Lemma~\ref{lemcenter} implies that item~(1) as
  well as the ``if" part of item~(3) of Theorem~\ref{thr:vsmall1}
  hold, with probability at least
  $1-c|\Lambda_L|^{-\alpha+\beta-\frac1d(\beta-\beta')}\ge
  1-c|\Lambda_L|^{-(\alpha-\alpha_{d,\rho,\tilde\rho})}$.\\
  This completes the proof of Theorem~\ref{thr:vsmall1}.
\end{proof}
\begin{proof}[Proof of Theorem \ref{thr:vbig1}]
  For a scale $\ell_\Lambda$ given, we set
  $q=[(L-\ell')/(\ell_\Lambda+\ell')]$. Then we may adjust the scale
  $\ell_\Lambda$ by enlarging it to a new scale $\ell$ so that
  $(\ell+\ell')q + \ell'= L$ and $0\le \ell-\ell_\Lambda \lesssim
  \ell_\Lambda^2 / |\Lambda_L| = o(\ell_\Lambda)$. As a consequence, we
  can consider a collection of boxes $\Lambda_\ell(\gamma_j)$ equally
  distant to their closest neighbors by a length $\ell'$ and
  satisfying the description of the theorem. In particular, events
  based on different boxes are independent.\\
  In this proof, we shall use the localization property described by
  Lemma~\ref{le:3} Part~(II).
  \\
  For $L$ sufficiently large, up to a probability $\lesssim
  |\Lambda_L|^{(1-\beta)}\exp^{-\ell^\nu} \le |\Lambda_L|^{-q}$, with
  $q>0$ arbitrary large, we can assume that all the boxes
  $\Lambda_\ell(\gamma)$ satisfy Part~(II) of Lemma~\ref{le:3}.\\
  It follows from \eqref{eq:93} and Lemma~\ref{le:3} that for $\ell$
  large enough, we have $ \P(\mathcal{V}_{\Lambda_L}^c) \lesssim
  |I_\Lambda| \ell^d $ and $\P(\mathcal{V}_{\Lambda_L}^c) \lesssim
  |I_\Lambda|^{1+\rho} \ell^{2d}$.  \\
  Let $\mathcal{S}_{\ell,L}$ be the set of disjoint boxes
  $\Lambda_{\ell}(\gamma_j)\subset\Lambda_L$ containing at least $2$
  centers of localization of $H_\omega(\Lambda_L)$. We set $n=\ell^d$.
  It follows from Lemma~\ref{lemcenter} (taking into account
  $\ell'\ll\ell$) that, using independence and Stirling's formula,
  \begin{equation*}
    \begin{split}
      \P(\sharp (\mathcal{S}_{\ell,L} \geq k) & \lesssim
      \binom{|\Lambda_L|/n}{k}
      (|I_\Lambda| n)^{(1+\rho)k} \\
      & \lesssim (\e |\Lambda_L|/(nk))^k (|I_\Lambda| n)^{(1+\rho)k} =
      \left(\frac{\e |\Lambda_L|}{ k}
        N(I_\Lambda)^{\frac{1+\rho}{1+\rho'} } n^{\rho} \right)^k
      \lesssim 2^{-k},
    \end{split}
  \end{equation*}
  if we choose
  \begin{equation}
    k \ge K:=\left[
      2\e N(I_\Lambda) |\Lambda_L|  (
      N(I_\Lambda)^{\frac{\rho-\rho'}{1+\rho'}} n^{\rho}) \right]+1.
    \label{defK}
  \end{equation}
  Note that
  \begin{align}
    K \asymp \frac{|\Lambda_L|}n \left( N(I_\Lambda)^{\frac{1}{1+\rho'}}
      n\right)^{1+\rho} = o\left(\frac{|\Lambda_L|}n \right)
  \end{align}
  by assumption.
  As a consequence, we get that
  \begin{equation*}
    \P(\# (\mathcal{S}_{\ell,L} ) \ge K)
    \lesssim 2^{-K}. 
  \end{equation*}
  So that, with probability larger than $1 - 2^{-K}$, we can assume
  that the boxes $\Lambda_\ell(\gamma_j)$, except at most $K$ of them,
  contain at most one center of localization.\\
  We now control the number of centers of localization that may be
  contained in these $K$ exceptional boxes. In a box of size $\ell$,
  the deterministic a priori bound on the number of eigenvalues
  guarantees that this number is bounded by $\lesssim \ell^d :=n$ (see
  e.g.~\cite{MR58:12429c}). Using this crude estimate the number of
  eigenvalues we miss with these $K$ boxes is bounded by
  \begin{equation*}
    K n \lesssim N(I_\Lambda) |\Lambda_L|  (
    N(I_\Lambda)^{\frac{\rho-\rho'}{1+\rho'}} n^{1+\rho}) = 
    o(N(I_\Lambda) |\Lambda_L| ),
  \end{equation*}
  provided
  \begin{equation}
    \label{condnl}
    N(I_\Lambda)^{\frac{\rho-\rho'}{1+\rho'}} n^{1+\rho} = o(1).
  \end{equation}
  We now turn to the complement of the $\Lambda_\ell(\gamma_j)$'s,
  that is
  \begin{equation*}
    \Upsilon=\Lambda_L\setminus
    \bigcup_j\Lambda_{\ell-\ell'}(\gamma_j)   ,
  \end{equation*}
  and we consider a partition of $\Upsilon$ in terms of $2^{d-1}$ sets
  of boxes of side length $\ell'$. More precisely,
  $\displaystyle\Upsilon = \bigcup_{m=1}^{2^{d-1}} \Upsilon_m$, with
  $\Upsilon_m\cap\Upsilon_{m'}=\emptyset$ if $m\neq m'$, each
  $\Upsilon_m$ is a union of boxes of side length $\ell'$ which are
  two by two distant by at least $\ell'$.\\
  Clearly, one has $|\Upsilon| \lesssim |\Lambda_L| \ell' / \ell$.\\
  Let $\mathcal{S}'_{\ell,L}$ be the set of boxes
  $\Lambda_{\ell'}(\gamma_j)\in \Upsilon$ containing at least one
  centers of localization of $H_\omega(\Lambda_L)$. From considerations
  similar to those above, for each given $m$, taking into account that
  boxes in $\Upsilon_m$ may be enlarged by, say, $\frac1 {10} \ell'$,
  for the distance between two of them is larger than $\ell'$, it
  follows from Lemma~\ref{lemcenter} and Wegner estimate that, using
  independence and Stirling formula,
  \begin{align*}
    \P(\# (\mathcal{S}'_{\ell,L} \cap \Upsilon_m ) \ge K') &\lesssim
    \sum_{k\ge K'} \binom{|\Upsilon_m|/ (\ell')^d}{k}
    (C |I_\Lambda|  (\ell')^d)^{k} \\
    &\lesssim \left(\frac {C |\Lambda_L| N(I_\Lambda)^{\frac 1{1+\rho'}}
        \ell'}{k\ell} \right)^k \lesssim 2^{-K'},
  \end{align*}
  where $C$ is a constant that varies but only depends on the constant
  appearing in Wegner and $d$, and provided one sets, for $C'>C$,
  \begin{equation}
    \label{defK'}
    K':=\left[ \frac {C'  |\Lambda_L| N(I_\Lambda)^{\frac 1{1+\rho'}}
        \ell'}{\ell}  \right] +1.
  \end{equation}
  Note that
  \begin{align}
    K' \asymp \frac{|\Lambda_L|}n \left(
      N(I_\Lambda)^{\frac{1}{1+\rho'}} n \frac{\ell'}{\ell}\right) =
    o\left(\frac{|\Lambda_L|}n \right)
  \end{align}
  by assumption. As a consequence, one computes
  \begin{equation*}
    \begin{split}
      \P(\# \mathcal{S}'_{\ell,L} \ge 2^{d-1} K') & \le \P(\exists m,
      \, \#(\mathcal{S}'_{\ell,L} \cap \Upsilon_m) \ge K' ) \\& \le
      \sum_{m=1}^{2^{d-1}} \P(\# (\mathcal{S}'_{\ell,L} \cap
      \Upsilon'_m ) \ge K') \lesssim 2^{-K'}.
    \end{split}
  \end{equation*}
  Hence, up to $2^{d-1} K'$ boxes, we can assume that boxes of size
  $\ell'$ in $\Upsilon'$ contain at most one center of localization.\\
  The maximal number of eigenvalues that can be contained in these
  $2^{d-1} K'$ bad boxes is $\lesssim K'n'=o(N(I_\Lambda)|\Lambda_L|)$
  provided
  \begin{equation}
    \label{condnn'}
    \ell \gg (\ell')^{d+1} .
  \end{equation}
  Combining~\eqref{condnl} and~\eqref{condnn'}, we see that the
  intermediate scale $\ell$ has to satisfy~\eqref{condell}.  To
  summarize, we proved that the picture described by
  Theorem~\ref{thr:vbig1} holds with a probability larger than
  \begin{multline}
    \label{probaev}
    1 - \e^{-c K} - \e^{ -cK'}  \ge 1 - \exp\left(-c N(I_\Lambda)
      |\Lambda_L| (N(I_\Lambda)^{\frac{\rho-\rho'}{1+\rho'} }
      \ell^{d\rho}) \right)\\ -\exp\left( -cN(I_\Lambda)|\Lambda_L|
      (N(I_\Lambda)^{-\frac {\rho'}{1+\rho'}} \ell'\ell^{-1})\right)
  \end{multline}
  Moreover, the number of eigenvalues of $H_\omega(\Lambda_L)$ that are
  not described by this picture is bounded by
  \begin{align}
    \label{badev}
    C(Kn + K'n') & \lesssim N(I_\Lambda) |\Lambda_L|\left( (
      N(I_\Lambda)^{\frac{\rho-\rho'}{1+\rho'}} \ell^{d(1+\rho)} +
      N(I_\Lambda)^{-\frac {\rho'}{1+\rho'}}
      \frac{(\ell')^{d+1}}{\ell}\right) \\
    & = o( N(I_\Lambda) |\Lambda_L|),
  \end{align}
  provided \eqref{condell} holds.\\
  This completes the proof of Theorem~\ref{thr:vbig1}.
\end{proof}



\section{The spectral statistics}
\label{sec:stat-conv}
In this section, we prove most of the results on the local spectral
statistics described in section~\ref{sec:results}.\\
The whole of our analysis relies on Theorems~\ref{thr:vsmall1}
and~\ref{thr:vbig1}.
\subsection{Convergence of the local levels statistics}
\label{sec:conv-levels-stat}
We first prove the uniform Poisson convergence, Theorem~\ref{thr:2},
of which Theorem~\ref{thr:3} is an immediate consequence if one takes
into account Remark~\ref{rem:13} to relax assumption~\eqref{eq:107}
into assumption~\eqref{eq:60}.
\subsubsection{Proof of Theorem~\ref{thr:2}}
\label{sec:proof-theorem2}
We keep the notations of section~\ref{sec:strong-poiss-conv}. Recall
that $\Lambda=\Lambda_L$. Under the assumptions of
Theorem~\ref{thr:2}, we can apply Theorem~\ref{thr:vsmall1} to the
interval $I_\Lambda=N^{-1}(N(E_0)+|\Lambda|^{-\alpha}[-1,1])$. Let
$\mathcal{Z}_\Lambda$ be the set of configurations where the
conclusions of Theorem~\ref{thr:vsmall1} for this interval hold and
let $\tilde N$ be the number of the cubes cons\-tructed in
Theorem~\ref{thr:vsmall1} and $(\Lambda_\ell(\gamma_j))_{1\leq j\leq
  \tilde N}$ be those cubes. Then, $\tilde
N=|\Lambda|^{1-\beta}(1-O(|\Lambda|^{-(1-\beta)/d}))$.\\
Recall that $\xi_j(E_0,\omega,\Lambda)$ is defined by~\eqref{eq:11}
and consider the event
\begin{equation*}
  \Omega^\Lambda_{I_1,k_1;I_2,k_2;\cdots;I_p,k_p}:=\bigcap_{l=1}^p
  \left\{\omega;\
    \#\{j;\ \xi_j(\omega,\Lambda)\in I_l\}=k_l \right\}.
\end{equation*}
Pick $\delta d<p:=\beta'\xi$ where $\beta'$ is given by
Theorem~\ref{thr:vsmall1} and $\delta$ by~\eqref{eq:32}.\\
As $\pro(\mathcal{Z}_\Lambda)\to1$, to prove Theorem~\ref{thr:2}, it
suffices to prove that
\begin{equation}
  \label{eq:56}
  \P\left(\Omega^\Lambda_{I_1,k_1;I_2,k_2;\cdots;I_p,k_p}\cap\mathcal{Z}_\Lambda
  \right)-e^{-|I_1|}\frac{|I_1|^{k_1}}{k_1!}\cdots
  e^{-|I_p|}\frac{|I_p|^{k_p}}{k_p!}\vers_{|\Lambda|\to+\infty}0.
\end{equation}
Recall that the function $N$ is non decreasing and continuous; thus,
for $I$ an interval, the set $N^{-1}(I)$ is an interval. For a cube
$\Lambda$ and an interval $I$, define the Bernoulli random variable
$X_{\Lambda_\ell,I}$ by
\begin{equation}
  \label{eq:111}
  X_{\Lambda_\ell,I}=\car_{H_\omega(\Lambda_\ell)\text{ has an
      e.v. in }N^{-1}[N(E_0)+|\Lambda|^{-1}I]\text{ with localization
      center in }\Lambda_{\ell-\ell'}}.
\end{equation}
Here, the length scales $\ell$ and $\ell'$ are taken as in
Theorem~\ref{thr:vsmall1} that is $\ell\asymp L^\beta$ and
$\ell'\asymp L^{\beta'}$. Notice that, using the notations of
section~\ref{sec:asymptotics}, one has $X_{\Lambda_\ell,I}=
X(\Lambda_\ell,N^{-1}[N(E_0)+|\Lambda|^{-1}I],\ell')$.\\
We are first going to eliminate a degenerate case i.e the case when
the length of one interval $|I_j|$ goes to $0$. Assume
$|I_1|\leq\varepsilon$ for some $\varepsilon$ small fixed. Let us
first assume that $k_1\not=0$. Then, by the description given by
Theorem~\ref{thr:vsmall1}, for some $\eta>0$, for $L$ sufficiently
large, one has
\begin{equation*}
  \begin{split}
    \P\left(\Omega^\Lambda_{I_1,k_1;I_2,k_2;\cdots;I_p,k_p}
      \cap\mathcal{Z}_\Lambda \right)&\leq
    \P\left(\Omega^\Lambda_{I_1,k_1}\cap\mathcal{Z}_\Lambda
    \right)\lesssim \tilde N\,
    \P\left(X_{\Lambda_\ell(\gamma),I_1+e^{-L^\eta}[-1,1]}=1\right)
    \\&\lesssim \tilde N\cdot
    N\left(N^{-1}\left[N(E_0)+|\Lambda|^{-1}\left(I_1+e^{-L^\eta}[-1,1]
        \right)\right]\right)|\Lambda_\ell|\\&\hskip1cm+\tilde
    N\left(|\Lambda_\ell|\left|N^{-1}\left[N(E_0)+|\Lambda|^{-1}
          \left(I_1+e^{-L^\eta}[-1,1]
          \right)\right]\right|\right)^{1+\rho}\\
    &\lesssim \varepsilon+|\Lambda||\Lambda_\ell|^\rho
    |\Lambda|^{-(1+\rho)/(1+\tilde\rho)}\varepsilon^{(1+\rho)/(1+\tilde\rho)}
    \\&\lesssim \varepsilon.
  \end{split}
\end{equation*}
To obtain the second inequality, we have used the upper bound provided
by~\eqref{eq:51} of Lemma~\ref{lemasympt}. The third inequality is a
consequence of assumption~\eqref{eq:60}.\\
On the other hand, if $k_1\not=0$, clearly, one has
$e^{-|I_1|}\frac{|I_1|^{k_1}}{k_1!}\lesssim \varepsilon$; moreover,
one always has $e^{-|I|}\frac{|I|^k}{k!}\leq 1$. Thus, $\displaystyle
  e^{-|I_1|}\frac{|I_1|^{k_1}}{k_1!}\cdots
  e^{-|I_p|}\frac{|I_p|^{k_p}}{k_p!}\lesssim\varepsilon$.\\
If $k_1=0$, then, by the computation made above, we have that
\begin{equation*}
  0\leq\P(\Omega^\Lambda_{I_2,k_2;\cdots;I_p,k_p}\cap\mathcal{Z}_\Lambda)-
  \P(\Omega^\Lambda_{I_1,0;I_2,k_2;\cdots;I_p,k_p}\cap\mathcal{Z}_\Lambda)
  \lesssim\varepsilon
\end{equation*}
and $\displaystyle 0\leq e^{-|I_2|}\frac{|I_1|^{k_2}}{k_2!}\cdots
e^{-|I_p|}\frac{|I_p|^{k_p}}{k_p!}-
e^{-|I_1|}e^{-|I_2|}\frac{|I_1|^{k_2}}{k_2!}\cdots
e^{-|I_p|}\frac{|I_p|^{k_p}}{k_p!}\lesssim\varepsilon$. Thus we are
back to estimating $\displaystyle
\P\left(\Omega^\Lambda_{I_2,k_2;\cdots;I_p,k_p}\right)
-e^{-|I_2|}\frac{|I_2|^{k_2}}{k_2!}\cdots
e^{-|I_p|}\frac{|I_p|^{k_p}}{k_p!}$.\\
So, from now on, we assume that all the intervals $(I_j)_j$ have
length large than $\varepsilon$. Define
$I_j^+=I_j\cup[-e^{-L^\eta},e^{-L^\eta}]$ and
$I_j^-=I_j\cap(^cI_j+[-e^{-L^\eta},e^{-L^\eta}])$. Clearly, $I_j^-\subset
I_j\subset I_j^+$. Moreover, by~\eqref{eq:32}, for $L$ sufficiently
large, we have that $I_j^+\cap I_k^+=I_j^-\cap I_k^-=\emptyset$ for
$j<k$ and $|I_j^+|=|I_j|(1+O(e^{-L^\eta}))$ and
$|I_j^-|=|I_j|(1+O(e^{-L^\eta}))$. For $L$ sufficiently large, one has
$|I_j^+|\geq|I_j^-|\geq\varepsilon/2$ for $j\in\{1,\dots,p\}$.\\
By Theorem~\ref{thr:vsmall1}, in particular~\eqref{error}, and the
Wegner estimate (W), we know that
\begin{equation}
  \label{eq:121}
  \begin{aligned}
  \bigcap_{l=1}^p \left\{\omega;\ \#\{j;\
    X_{\Lambda_\ell(\gamma_j),I^-_l}=1\}=k_l
  \right\}\bigcap\mathcal{Z}_\Lambda\subset
  \Omega^\Lambda_{I_1,k_1;I_2,k_2;\cdots;I_p,k_p}\cap\mathcal{Z}_\Lambda,\\
  \Omega^\Lambda_{I_1,k_1;I_2,k_2;\cdots;I_p,k_p}\cap\mathcal{Z}_\Lambda
  \subset\bigcap_{l=1}^p \left\{\omega;\ \#\{j;\
    X_{\Lambda_\ell(\gamma_j),I^+_l}=1\}=k_l
  \right\}\bigcap\mathcal{Z}_\Lambda.    
  \end{aligned}
\end{equation}
Let us first use this to upper bound
$\P(\Omega^\Lambda_{I_1,k_1;I_2,k_2;\cdots;I_p,k_p}
\cap\mathcal{Z}_\Lambda)$ and to show that we may assume $p$ to be
finite (depending on $\varepsilon$). Therefore, we compute
\begin{equation*}
  \begin{split}
    \P(\Omega^\Lambda_{I_1,k_1;I_2,k_2;\cdots;I_p,k_p}
    \cap\mathcal{Z}_\Lambda)&\leq \P \left(\bigcap_{l=1}^p
      \left\{\omega;\ \#\{j;\ X_{\Lambda_\ell(\gamma_j),I^+_l}=1\}=k_l
      \right\}\cap\mathcal{Z}_\Lambda\right) \\&\leq\sum_{\substack{
        K_l\in\{1,\cdots,\tilde N\}\\ \#K_l=k_l,\ 1\leq l\leq
        p\\K_l\cap K_{l'}=\emptyset\text{ if }l\not= l'}}\P
    \left(\bigcap_{l=1}^p \left\{\omega;
        \begin{aligned}
          \forall j\in K_l,\ X_{\Lambda_\ell(\gamma_j),I^+_l}&=1,\\
          \forall j\not\in K_l,\ X_{\Lambda_\ell(\gamma_j),I^+_l}&=0
        \end{aligned}
      \right\}\cap\mathcal{Z}_\Lambda \right)
  \end{split}
\end{equation*}
as, by the definition of $\mathcal{Z}_\Lambda$ (see
Theorem~\ref{thr:vsmall1}), one has
$\P(\{X_{\Lambda_\ell(\gamma_j),I^+_l}=1,
X_{\Lambda_\ell(\gamma_j),I^+_{l'}}=1\}\cap\mathcal{Z}_\Lambda)=0$ if
$l\not=l'$.  As $K_l\cap K_{l'}=\emptyset$ if $l\not= l'$ the random
vectors $((X_{\Lambda_\ell(\gamma_j),I^+_l})_{j\in K_l})_{1\leq l\leq
  p}$ are independent. Thus, one estimates
\begin{equation}
  \label{eq:112}
  \begin{split}
    \P(\Omega^\Lambda_{I_1,k_1;I_2,k_2;\cdots;I_p,k_p}
    \cap\mathcal{Z}_\Lambda)&\leq \sum_{\substack{
        K_l\in\{1,\cdots,\tilde N\}\\ \#K_l=k_l,\ 1\leq l\leq
        p\\K_l\cap K_{l'}=\emptyset\text{ if }l\not= l'}}
    \prod_{l=1}^p\P \left(\left\{\omega;\ \forall j\in K_l,\
        X_{\Lambda_\ell(\gamma_j),I^+_l}=1 \right\}\right)
    \\&\hskip2.2cm\times\P \left(\left\{\omega;\ \forall
        j\not\in\bigcup_{l=1}^p K_l,\ \sum_{l=1}^p
        X_{\Lambda_\ell(\gamma_j),I^+_l}=0 \right\}\right) \\& \leq
    \sum_{\substack{ K_l\in\{1,\cdots,\tilde N\}\\ \#K_l=k_l,\ 1\leq
        l\leq p}} \prod_{l=1}^p\P \left(\left\{\omega;\ \forall j\in
        K_l,\ X_{\Lambda_\ell(\gamma_j),I^+_l}=1 \right\}\right)
    \\&\hskip2.2cm\times\P \left(\left\{\omega;\ \forall
        j\not\in\bigcup_{l=1}^p K_l,\ \sum_{l=1}^p
        X_{\Lambda_\ell(\gamma_j),I^+_l}=0 \right\}\right) \\& \leq
    \binom{\tilde N}{k_1+\cdots+k_p} \prod_{l=1}^p
    \left(p^+_l\right)^{k_l} \\&\hskip2.5cm\times\P
    \left(\left\{\omega;\ \sum_{l=1}^p
        X_{\Lambda_\ell(\gamma_1),I^+_l}=0 \right\}\right)^{\tilde
      N-(k_1+\cdots+k_p)}
  \end{split}
\end{equation}
where $p^\pm_l=\P(X_{\Lambda_\ell(\gamma_1),I^\pm_l}=1)$ for $1\lambda
\leq p$. Here, we have used the fact that, as the cubes
$(\Lambda_\ell(\gamma_j))_j$ are at distance at least $R$ from each
other (see (IAD)), the events $\left(\left\{\omega;\ \#\{j;\
    X_{\Lambda_\ell(\gamma_j),I^-_l}=1\}=k_l \right\}\right)_j$ are
two by two independent for $1\leq l\leq p$.\\
To estimate the last term in~\eqref{eq:112}, we will use
\begin{Le}
  \label{le:10}
  Set
  \begin{equation}
    \label{eq:117}
    \delta:=\alpha\frac{1+\rho}{1+\tilde\rho}-1-\beta\rho.
  \end{equation}
  With the choice of $(I^+_l)_{1\leq l\leq p}$ made above, under the
  assumptions of Theorem~\ref{thr:2}, with our choice of $\ell$ and
  $\ell'$, possibly reducing $\beta$ somewhat, for $L$ sufficiently
  large, one has
  \begin{equation*}
    \P\left(\sum_{l=1}^pX_{\Lambda_\ell(\gamma_1),I^+_l}=0\right)
    = 1-\frac{1-O\left(|\Lambda|^{-\delta}\right)}{\tilde N}
    \sum_{l=1}^p|I_l|. 
  \end{equation*}
\end{Le}
\begin{proof}
  The proof is analogous to that of Lemma~\ref{lemasympt}; it also
  relies on our choice for the intervals $(I_l)_{1\leq l\leq p}$. The
  derivation of~\eqref{eq:114} yields
  \begin{equation*}
    0\leq
    \E\left[N\left(\bigcup_{j=1}^pN^{-1}(N(E_0)+|\Lambda|^{-1}I_l^+)
        ,\ell,\ell'\right)\right]-
    \P\left(\sum_{l=1}^pX_{\Lambda_\ell(\gamma_j),I^+_l}\geq1\right)
    \leq C |\Lambda_\ell|^{1+\rho}|I_\Lambda|^{1+\rho}.
  \end{equation*}
  Thus using Lemma~\ref{le:9} and assumption~\eqref{eq:60}, we obtain
  \begin{equation}
    \label{eq:115}
    \begin{split}
      \P\left(\sum_{l=1}^pX_{\Lambda_\ell(\gamma_j),I^+_l}\geq
        1\right)
      &=N\left(\bigcup_{j=1}^pN^{-1}(N(E_0)+|\Lambda|^{-1}I_l^+)\right)+
      |\Lambda_\ell|^{1+\rho}|N(I_\Lambda)|^{(1+\rho)/(1+\tilde\rho)}
      \\&=|\Lambda|^{\beta-1}\bigcup_{j=1}^p|I_l^+|+
      O(|\Lambda|^{\beta(1+\rho)-\alpha(1+\rho)/(1+\tilde\rho)}).
    \end{split}
  \end{equation}
  To conclude the statement of Lemma~\ref{le:10}, it now suffices to
  recall that $|I_l^+|\geq\varepsilon/2$ and that
  $1+\beta\rho<\alpha\frac{1+\rho}{1+\tilde\rho}$
  i.e. $\beta-1>\beta(1+\rho)-\alpha\frac{1+\rho}{1+\tilde\rho}$. Moreover,
  by the definition of the intervals $(I_l^+)_l$, the main term
  in~\eqref{eq:115} is $|\Lambda|^{\beta-1}\bigcup_{j=1}^p|I_l|$.
\end{proof}
\noindent By~\eqref{eq:51} of Lemma~\ref{lemasympt}, the
assumption~\eqref{eq:107} and the definition of $(I_l^\pm)_l$, we have
that
\begin{equation}
  \label{eq:108}
  p^{\pm}_l=|I_l||\Lambda|^{\beta-1}(1+O(|\Lambda|^{-\delta}))
  =|I_l|\tilde N^{-1}(1+O(|\Lambda|^{-\delta})).
\end{equation}
Setting $k=k_1+\cdots+k_p$ and $I^+=|I^+_1|+\cdots+|I^+_1|$, we note
that by Theorem~\ref{thr:vsmall1}, for $\omega\in\mathcal{Z}_\Lambda$,
we have $k\lesssim N(I_\Lambda)|\Lambda|\lesssim|\Lambda|^{1-\alpha}$
and, by the assumptions of Theorem~\ref{thr:2}, we have $I^+\leq
2|\Lambda|^{1-\alpha}$. Thus, for $L$ large, we have
\begin{equation}
  \label{eq:118}
  k+I^+=o(|\Lambda|^{1-\beta})=o(\tilde N).
\end{equation}
From~\eqref{eq:112},~\eqref{eq:108} and Lemma~\ref{le:10}, at the
possible expense of reducing $\beta$ somewhat, we, thus, obtain the
following upper bound
\begin{equation}
  \label{eq:119}
  \begin{split}
    \P(\Omega^\Lambda_{I_1,k_1;I_2,k_2;\cdots;I_p,k_p}
    \cap\mathcal{Z}_\Lambda)&\leq \binom{\tilde N}{k}
    \left(\frac{I^+}{\tilde N} \right)^k \left(1-\frac{1-
        O\left(|\Lambda|^{-\delta}\right)}{\tilde
        N}I^+\right)^{\tilde N-k}\\
    &\leq \frac{(I^+)^k}{k!}e^{-I^+}e^{O(|\Lambda|^{-\delta}I^+)}
    \lesssim \frac{(I^+)^k}{k!}e^{-I^+}
  \end{split}
\end{equation}
where we have used~\eqref{eq:118} as well as the
assumption~\eqref{eq:116} to obtain
\begin{equation}
  \label{eq:120}
  |\Lambda|^{-\delta}I^+\leq|\Lambda|^{-\delta+1-\alpha}=o(1)
\end{equation}
Indeed, by~\eqref{eq:117}, the definition of $\delta$, and as
$\alpha_{d,\rho,\tilde\rho}<\alpha$ (see~\eqref{eq:59}), one computes
\begin{equation*}
  \begin{split}
    -\delta+1-\alpha=
    -\alpha\frac{1+\rho}{1+\tilde\rho}+1+\beta\rho+1-\alpha&<
    -\frac{(2+\rho+\tilde\rho)(1+d\rho)}{1+(d+1)\rho}+2+\beta\rho \\&=
    -\frac{\tilde\rho(1+d\rho)-\rho(1-d\rho)}{1+(d+1)\rho}+\beta\rho
  \end{split}
\end{equation*}
So, under assumption~\eqref{eq:116}, at the expense of possibly
reducing $\beta$, one has $-\delta+1-\alpha<0$ which, taking into
account~\eqref{eq:120}, implies~\eqref{eq:119}.\\
The bound~\eqref{eq:119} proves that if $k+I^+\to+\infty$ (as
$L\to+\infty$) then
$\P(\Omega^\Lambda_{I_1,k_1;I_2,k_2;\cdots;I_p,k_p}
\cap\mathcal{Z}_\Lambda)\to0$. Note that, for $1\leq l\leq p$, one has
$|I_l^+|\geq\varepsilon/2$, one has
$\P(\Omega^\Lambda_{I_1,k_1;I_2,k_2;\cdots;I_p,k_p}
\cap\mathcal{Z}_\Lambda)\to0$ if $p\to+\infty$ as $L\to+\infty$. On
the other hand, $0$ is clearly also the limit of the product
\begin{equation*}
  \frac{|I_1|^{k_1}}{k_1!}\e^{-|I_1|}\cdots
  \frac{|I_p|^{k_p}}{k_p!} \e^{-|I_p|}
\end{equation*}
in any of these cases (as $|I^+_l|\asymp|I_l|\geq\varepsilon/2$). So
we have proved~\eqref{eq:110} if $k+I^++p\to+\infty$ when $L\to+\infty$.\\
From now on, we assume that $k$, $I^+$ and $p$ are bounded. Let us
prove~\eqref{eq:110} in this case. By~\eqref{eq:121}, using the same
computation as above, we have
\begin{equation*}
  \begin{aligned}
    \sum_{\substack{
        K_l\in\{1,\cdots,\tilde N\}\\ \#K_l=k_l,\ 1\leq l\leq
        p\\K_l\cap K_{l'}=\emptyset\text{ if }l\not= l'}}\P
    \left(\bigcap_{l=1}^p \left\{\omega;
        \begin{aligned}
          \forall j\in K_l,\ X_{\Lambda_\ell(\gamma_j),I^-_l}&=1,\\
          \forall j\not\in \bigcup_{l'\not=l}K_{l'},\
          X_{\Lambda_\ell(\gamma_j),I^-_l}&=0
        \end{aligned}
      \right\}\cap\mathcal{Z}_\Lambda \right) \leq
    \P(\Omega^\Lambda_{I_1,k_1;I_2,k_2;\cdots;I_p,k_p}
    \cap\mathcal{Z}_\Lambda),\\
    \P(\Omega^\Lambda_{I_1,k_1;I_2,k_2;\cdots;I_p,k_p}
    \cap\mathcal{Z}_\Lambda)\leq\sum_{\substack{
        K_l\in\{1,\cdots,\tilde N\}\\ \#K_l=k_l,\ 1\leq l\leq
        p\\K_l\cap K_{l'}=\emptyset\text{ if }l\not= l'}}\P
    \left(\bigcap_{l=1}^p \left\{\omega;
        \begin{aligned}
          \forall j\in K_l,\ X_{\Lambda_\ell(\gamma_j),I^+_l}&=1,\\
          \forall j\not\in\bigcup_{l'\not=l}K_{l'},\
          X_{\Lambda_\ell(\gamma_j),I^+_l}&=0
        \end{aligned}
      \right\}\cap\mathcal{Z}_\Lambda \right).
  \end{aligned}
\end{equation*}
Hence, as $1-\P(\mathcal{Z}_\Lambda)=o(1)$ and $p$ is finite, one has
\begin{equation*}
  \begin{aligned}
    \sum_{\substack{
        K_l\in\{1,\cdots,\tilde N\}\\ \#K_l=k_l,\ 1\leq l\leq
        p\\K_l\cap K_{l'}=\emptyset\text{ if }l\not= l'}}\P
    \left(\bigcap_{l=1}^p \left\{\omega;
        \begin{aligned}
          \forall j\in K_l,\ X_{\Lambda_\ell(\gamma_j),I^-_l}&=1,\\
          \forall j\not\in K_l,\ X_{\Lambda_\ell(\gamma_j),I^-_l}&=0
        \end{aligned}
      \right\} \right)+o(1) \leq
    \P(\Omega^\Lambda_{I_1,k_1;I_2,k_2;\cdots;I_p,k_p}
    \cap\mathcal{Z}_\Lambda),\\
    \P(\Omega^\Lambda_{I_1,k_1;I_2,k_2;\cdots;I_p,k_p}
    \cap\mathcal{Z}_\Lambda)\leq\sum_{\substack{
        K_l\in\{1,\cdots,\tilde N\}\\ \#K_l=k_l,\ 1\leq l\leq
        p\\K_l\cap K_{l'}=\emptyset\text{ if }l\not= l'}}\P
    \left(\bigcap_{l=1}^p \left\{\omega;
        \begin{aligned}
          \forall j\in K_l,\ X_{\Lambda_\ell(\gamma_j),I^+_l}&=1,\\
          \forall j\not\in K_l,\ X_{\Lambda_\ell(\gamma_j),I^+_l}&=0
        \end{aligned}
      \right\} \right)+o(1).
  \end{aligned}
\end{equation*}
For $(K_l)_{1\leq l\leq p}\in\{1,\cdots,\tilde N\}^p$ such that
$(\#K_l)_{1\leq l\leq p}=(k_l)_{1\leq l\leq p}$ and $K_l\cap
K_{l'}=\emptyset$ if $l\not= l'$, one computes
\begin{equation*}
  \begin{split}
  \P\left(\bigcap_{l=1}^p \left\{\omega;
      \begin{aligned}
        \forall j\in K_l,\ X_{\Lambda_\ell(\gamma_j),I^+_l}&=1,\\
        \forall j\not\in\bigcup_{l'\not=l}K_{l'},\
        X_{\Lambda_\ell(\gamma_j),I^+_l}&=0
      \end{aligned}
    \right\} \right)&\leq \prod_{j\not\in\cup_{l=1}^p K_l}
  \P\left(\sum_{l=1}^pX_{\Lambda_\ell(\gamma_j),I^+_l}=0\right)
  \prod_{l=1}^p\prod_{j\in K_l} \P\left(
    X_{\Lambda_\ell(\gamma_j),I^+_l}=1\right) \\&\leq
  \P\left(\sum_{l=1}^pX_{\Lambda_\ell(\gamma_1),I^+_l}=0\right)^{\tilde
    N-(k_1+\cdots+k_p)} \prod_{l=1}^p
  \P(X_{\Lambda_\ell(\gamma_1),I^+_l}=1)^{k_l} \\& \leq
  \left(1-\sum_{l=1}^pp^+_l\right)^{\tilde N-(k_1+\cdots+k_p)}
  \prod_{l=1}^p \left((p^+_l\right)^{k_l}(1+o(1)) \\& \leq
  \left(\prod_{l=1}^pe^{-|I^+_l|}|I^+_l|^{k_l} \right) (\tilde
  N)^{k_1+\cdots+k_p}(1+o(1))
  \end{split}
\end{equation*}
and
\begin{equation*}
  \begin{split}
  \P\left(\bigcap_{l=1}^p \left\{\omega;
      \begin{aligned}
        \forall j\in K_l,\ X_{\Lambda_\ell(\gamma_j),I^-_l}&=1,\\
        \forall j\not\in K_l,\ X_{\Lambda_\ell(\gamma_j),I^-_l}&=0
      \end{aligned}
    \right\} \right)&\geq
  \P\left(\sum_{l=1}^pX_{\Lambda_\ell(\gamma_1),I^-_l}=0\right)^{\tilde
    N-(k_1+\cdots+k_p)} \\&\hskip1cm\times\prod_{j\in\bigcup_l K_l}\left(
    \sum_{l=1}^p \car_{j\in K_l}\P\left(
        \begin{aligned}
          X_{\Lambda_\ell(\gamma_j),I^-_l}=1,\\
          \sum_{l'\not=l}X_{\Lambda_\ell(\gamma_j),I^-_{l'}}=0
        \end{aligned}
      \right) \right)\\&\geq
    \P\left(\sum_{l=1}^pX_{\Lambda_\ell(\gamma_1),I^-_l}=0\right)^{\tilde
      N-(k_1+\cdots+k_p)} \\&\times\prod_{j\in\bigcup_l K_l}\left(
      \sum_{l=1}^p \car_{j\in K_l}\left[\P\left(
        X_{\Lambda_\ell(\gamma_j),I^-_l}=1 \right)-\P\left(
        \begin{aligned}
          X_{\Lambda_\ell(\gamma_j),I^-_l}=1,\\
          X_{\Lambda_\ell(\gamma_j),I^-_{l'}}=1
        \end{aligned}
      \right) \right]\right) \\&\geq
  \left(\prod_{l=1}^pe^{-|I^+_l|}|I^-_l|^{k_l} \right) (\tilde
  N)^{k_1+\cdots+k_p}(1+o(1))
  \end{split}
\end{equation*}
as $k_1+\cdot+k_p$ is bounded and
\begin{equation*}
  \P\left(
        X_{\Lambda_\ell(\gamma_j),I^-_l}=1 \right)-\P\left(
        \begin{aligned}
          X_{\Lambda_\ell(\gamma_j),I^-_l}=1,\\
          X_{\Lambda_\ell(\gamma_j),I^-_{l'}}=1
        \end{aligned}
      \right)=p_l^-\left(1+O\left(|\Lambda|^{-\delta}\right)\right).
\end{equation*}
On the other hand, as $k_1+\cdots+k_p$ is bounded, when $N\to+\infty$,
one has
\begin{equation*}
      \sum_{\substack{
        K_l\in\{1,\cdots,\tilde N\}\\ \#K_l=k_l,\ 1\leq l\leq
        p\\K_l\cap K_{l'}=\emptyset\text{ if }l\not= l'}}1=
    \prod_{l=1}^p\binom{\tilde N}{k_l}(1+o(1))
\end{equation*}
Thus, for $L$ sufficiently large, we obtain
\begin{equation*}
  \begin{split}
    e^{-|I^-_1|}\frac{|I^-_1|^{k_1}}{k_1!}\cdots
    e^{-|I^-_p|}\frac{|I^-_p|^{k_p}}{k_p!}(1+o(1))&\leq
    \P(\Omega^\Lambda_{I_1,k_1;I_2,k_2;\cdots;I_p,k_p}
    \cap\mathcal{Z}_\Lambda)\\&\leq
    e^{-|I^+_1|}\frac{|I^+_1|^{k_1}}{k_1!}\cdots
    e^{-|I^+_p|}\frac{|I^+_p|^{k_p}}{k_p!}(1+o(1))
  \end{split}
\end{equation*}
Now, recalling that $|I^\pm_l|=|I_l|+O\left(e^{-L^\eta}\right)$, we
get
\begin{equation*}
  \lim_{|\Lambda|\to+\infty}
  \left|\P\left(\Omega^\Lambda_{I_1,k_1;I_2,k_2;\cdots;I_p,k_p}\right)
    -e^{-|I_1|}\frac{|I_1|^{k_1}}{k_1!}\cdots
    e^{-|I_p|}\frac{|I_p|^{k_p}}{k_p!}\right|=0
\end{equation*}
and the proof of Theorem~\ref{thr:2} is complete.\qed
\begin{Rem}
  \label{rem:13}
  In the present proof, assumption~\eqref{eq:60} does not suffice to
  guarantee~\eqref{eq:108}: indeed, as we did not fix the intervals
  $(I_j)_j$, we want a result uniform over all the intervals of not
  too small length in some neighborhood of $E_0$; thus, we use
  assumption~\eqref{eq:107}, a uniform version of
  assumption~\eqref{eq:60}.\\
  To prove Theorem~\ref{thr:3} however it suffices to consider fixed
  intervals $(I_j)_{1\leq j\leq p}$; thus, assumption~\eqref{eq:60}
  suffices to guarantee that~\eqref{eq:108} holds.
\end{Rem}
\subsubsection{The asymptotic independence}
\label{sec:asympt-indep}
We now turn to the proof of Theorem~\ref{thr:4}. Decompose $\Lambda$
into the boxes constructed in Theorem~\ref{thr:vsmall1}
i.e. $\displaystyle \Lambda=\cup_{\gamma\in\Gamma_\Lambda}
\Lambda_\ell(\gamma)$ ; the scale $\ell$ is determined by
Theorem~\ref{thr:vsmall1} and
$N:=\#\Gamma_\Lambda\sim|\Lambda|\ell^{-d}$. The set of boxes thus
obtained are the same for both energies $E$ and $E'$. For $I\subset
\R$ a compact set, define the random variables
\begin{equation}
  \label{eq:17}
  X_\gamma(E,I)=
  \begin{cases}
    1\text{ if } H_\omega(\Lambda_\ell(\gamma))\text{
      has an e.v. in }N^{-1}[N(E_0)+|\Lambda|^{-1}I],\\
    0\text{ if not}.
  \end{cases}
\end{equation}
Then, to prove Theorem~\ref{thr:4}, it suffices to prove that, for
$(k,k')\in\N^2$ and $I\subset\R$ and $I'\subset\R$ two compact sets,
one has
\begin{multline}
  \label{eq:18}
  \lim_{|\Lambda|\to+\infty}\esp
  \left[\left(\sum_{\gamma\in\Gamma_\Lambda}X_\gamma(E,I)\right)^k
    \left(\sum_{\gamma\in\Gamma_\Lambda}X_\gamma(E',I')\right)^{k'}
  \right]\\=\lim_{|\Lambda|\to+\infty}\esp
  \left[\left(\sum_{\gamma\in\Gamma_\Lambda}X_\gamma(E,I)\right)^k\right]
  \esp\left[\left(\sum_{\gamma\in\Gamma_\Lambda}X_\gamma(E',I')\right)^{k'}
  \right].
\end{multline}
Therefore, using the independence of the cubes, we expand the sums
\begin{equation*}
  \begin{split}
    S_N(k,k')=&\esp\left[\left(\sum_{\gamma\in\Gamma_\Lambda}
        X_\gamma(E,I)\right)^k
      \left(\sum_{\gamma\in\Gamma_\Lambda}X_\gamma(E',I')\right)^{k'}
    \right]\\&=
    \sum_{\gamma_1,\cdots,\gamma_k}\sum_{\gamma'_1,\cdots,\gamma'_{k'}}
    \esp\left[X_{\gamma_1}(E,I)\cdots X_{\gamma_k}(E,I)\cdot
      X_{\gamma'_1}(E',I')\cdots X_{\gamma'_{k'}}(E',I') \right]\\&=
    G_N(k,k')+R_N(k,k')
  \end{split}
\end{equation*}
where
\begin{equation*}
  G_N(k,k')=\sum_{\{\gamma_1,\cdots,\gamma_k\}\cap
    \{\gamma'_1,\cdots,\gamma'_{k'}\}=\emptyset}
  \esp\left[X_{\gamma_1}(E,I)\cdots X_{\gamma_k}(E,I)\cdot
    X_{\gamma'_1}(E',I')\cdots X_{\gamma'_{k'}}(E',I') \right]
\end{equation*}
and
\begin{equation*}
  \begin{split}
    R_N(k,k')&=\sum_{\gamma} \sum_{\substack{\gamma_1,\cdots,\gamma_l
        \\\gamma'_1,\cdots,\gamma'_{l'}}}
    \esp\left[X_{\gamma}(E,I)X_{\gamma}(E',I')\cdot\prod_{j=1}^{k-1}
      X_{\gamma_j}(E,I)\cdot\prod_{j'=1}^{k'-1}
      X_{\gamma'_{j'}}(E',I') \right]\\&= \sum_{\gamma}
    \sum_{l=1}^{k-1}\sum_{l'=1}^{k'-1}\binom{l}{k-1}
    \binom{l'}{k'-1}\\&\hskip2cm
    \sum_{\substack{\gamma\not\in\{\gamma_1,\cdots,\gamma_l\}
        \\\gamma\not\in\{\gamma'_1,\cdots,\gamma'_{l'}\}}}
    \esp\left[X_{\gamma}(E,I)X_{\gamma}(E',I')\cdot\prod_{j=1}^l
      X_{\gamma_j}(E,I)\cdot\prod_{j'=1}^{l'} X_{\gamma'_{j'}}(E',I')
    \right]\\&=
    \sum_{\gamma}\esp\left[X_{\gamma}(E,I)X_{\gamma}(E',I')\right]
    \sum_{l=1}^{k-1}\sum_{l'=1}^{k'-1}\binom{l}{k-1}\binom{l'}{k'-1}
    \\&\hskip2cm\sum_{\substack{\gamma\not\in\{\gamma_1,\cdots,\gamma_l\}
        \\\gamma\not\in\{\gamma'_1,\cdots,\gamma'_{l'}\}}}
    \esp\left[\prod_{j=1}^l X_{\gamma_j}(E,I)\cdot\prod_{j'=1}^{l'}
      X_{\gamma'_{j'}}(E',I') \right].
  \end{split}
\end{equation*}
Hence,
\begin{equation}
  \label{eq:65}
  R_N(k,k')=\sum_{\gamma}\esp\left[X_{\gamma}(E,I)X_{\gamma}(E',I')\right]
    \sum_{l=1}^{k-1}\sum_{l'=1}^{k'-1}\binom{l}{k-1}
    \binom{l'}{k'-1}S_{N-1}(l,l').
\end{equation}
On the other hand, one computes
\begin{equation*}
  \begin{split}
    S_N(k)S_N(k')&=
    \esp\left[\left(\sum_{\gamma\in\Gamma_\Lambda}X_\gamma(E,I)\right)^k\right]
    \esp\left[\left(\sum_{\gamma\in\Gamma_\Lambda}X_\gamma(E',I')\right)^{k'}\right]
    \\
    &=\sum_{\gamma_1,\cdots,\gamma_k}\esp\left[X_{\gamma_1}(E,I)\cdots
      X_{\gamma_k}(E,I) \right] \sum_{\gamma'_1,\cdots,\gamma'_{k'}}
    \esp\left[ X_{\gamma'_1}(E',I')\cdots X_{\gamma'_{k'}}(E',I')
    \right]\\
    &=G_N(k,k') +Q_N(k,k')
  \end{split}
\end{equation*}
where
\begin{equation}
  \label{eq:66}
  \begin{split}
  Q_N(k,k')&=\sum_{\gamma} \sum_{\gamma_1,\cdots,\gamma_{k-1}}
    \esp\left[X_{\gamma}(E,I)X_{\gamma_1}(E,I)\cdots
      X_{\gamma_{k-1}}(E,I) \right]\\&\hskip1cm\cdot
    \sum_{\gamma'_1,\cdots,\gamma'_{k'-1}} \esp\left[X_{\gamma}(E',I')
      X_{\gamma'_1}(E',I')\cdots X_{\gamma'_{k'-1}}(E',I') \right]\\
    &=\sum_{\gamma}\esp\left[X_{\gamma}(E,I)\right]
    \esp\left[X_{\gamma}(E',I')\right]
    \left[\sum_{l=1}^{k-1}\binom{l}{k-1}S_{N-1}(l)\right]^2.
  \end{split}
\end{equation}
Hence,
\begin{equation}
  \label{eq:64}
  \begin{split}
    S_N(k,k')-S_N(k)S_N(k')=R_N(k,k')-Q_N(k,k').
  \end{split}
\end{equation}
By Cauchy-Schwartz, one has
\begin{equation*}
  S_N(k,k')\leq\sqrt{S_N(2k)S_N(2k')}.  
\end{equation*}
On the other hand, as $\pro(X_\gamma=1)\leq C/N$ by
Lemma~\ref{lemasympt}, one computes
\begin{equation*}
    S_N(k)\leq \sum_{j=0}^N\binom{j}{N}\left(\frac CN\right)^j\,j^k
    \leq C_k\sum_{j=0}^N\binom{j}{N}\,\left(\frac{Ck}N\right)^j
    \leq C_k\,e^{N\log(1+Ck/N)}\leq C_k\,e^{Ck}<+\infty
\end{equation*}
where $\displaystyle C_k=e^{k(\log k-1-\log\log k)}$ for $k\geq2$.\\
So, for any $k$ and $k'$, one has
\begin{equation*}
  \sup_{N\geq1}|S_N(k)|+|S_N(k,k')|<+\infty,
\end{equation*}
Thus, using~\eqref{eq:65},~\eqref{eq:66} and~\eqref{eq:64}, one obtains
\begin{multline*}
  |S_N(k,k')-S_N(k)S_N(k')|\\\leq C_{k,k'}\max\left(\sum_{\gamma}
    \esp\left[X_{\gamma}(E,I)\right]
    \esp\left[X_{\gamma}(E',I')\right],\sum_{\gamma}
    \esp\left[X_{\gamma}(E,I)X_{\gamma}(E',I')\right] \right)
\end{multline*}
Hence,~\eqref{eq:18} and thus Theorem~\ref{thr:4} and~\ref{thr:8}
follow from the following two properties
\begin{equation*}
  \sum_{\gamma}\esp\left[X_{\gamma}(E,I)\right]
  \esp\left[X_{\gamma}(E',I')\right]\vers_{|\Lambda|\to+\infty}0\quad
  \text{and}\quad 
  \sum_{\gamma}\esp\left[X_{\gamma}(E,I)X_{\gamma}(E',I')\right]
  \vers_{|\Lambda|\to+\infty}0.
\end{equation*}
As the operators $(H_\omega(\Lambda_\ell(\gamma)))_{\gamma}$ are
i.i.d, this will be proved if we prove that
\begin{equation}
  \label{eq:20}
  \left(\frac{L}{\ell}\right)^d\esp\left[X_0(E,I)\right]
  \esp\left[X_0(E',I')\right]\vers_{|\Lambda|\to+\infty}0\quad
  \text{and}\quad 
  \left(\frac{L}{\ell}\right)^d\esp\left[X_0(E,I)X_0(E',I')\right]
  \vers_{|\Lambda|\to+\infty}0.
\end{equation}
The first limit in~\eqref{eq:20} is an immediate consequence of the
Wegner estimate as $\ell/L\to0$. The second limit in~\eqref{eq:20} clearly is a
consequence of (D) if $E\not= E'$ are fixed energies, and of (GM) if
$E=E_\Lambda$ and $E'=E'_\Lambda$ satisfy the assumptions of
Theorem~\ref{thr:8}.\\
This completes the proof of Theorem~\ref{thr:4} and~\ref{thr:8}.\qed
\subsection{Study of the (levels, centers) statistics}
\label{sec:conv-poiss-levels}
In this section, we will prove Theorems~\ref{thr:11},~\ref{thr:5}
and~\ref{thr:7}. To control the eigenvalues, as in the previous
section, we use Theorem~\ref{thr:vsmall1}. So we keep the same
notations here.
\subsubsection{The proof of Theorem~\ref{thr:11}}
\label{sec:proof-theorem11}
The proof of Theorem~\ref{thr:11} is very similar to that of
Theorem~\ref{thr:2}. One deals with the case when the length of some
interval $(I_k)_k$ tend to $0$ or $\infty$ as in the proof of
Theorem~\ref{thr:2}; we will not repeat this here. We only indicate
the differences and keep the same notations. We need to estimate the
probability of the event
\begin{equation*}
  \Omega^\Lambda_{(I_l,C_l,k_l)_{1\leq l\leq p}}:=\bigcap_{l=1}^p
  \left\{\omega;\
    \#\left\{j;\
      \begin{aligned}
        \xi_j(\omega,\Lambda)&\in I_l\\ x_j(\omega,\Lambda)/L&\in C_l
      \end{aligned}
    \right\}=k_l \right\}.
\end{equation*}
Let us first deal with the case when the volume of one of the cubes
$(C_k)_k$ tends to $0$; assume now that $|C_1|\leq\varepsilon$,
$|I_1|\geq\varepsilon$ and that $k_1\geq1$. Then, we have that
\begin{equation*}
  \begin{split}
    \P\left(\Omega^\Lambda_{(I_l,C_l,k_l)_{1\leq l\leq
          p}}\cap\mathcal{Z}_\Lambda \right)&\leq
    \P\left(\Omega^\Lambda_{I_1,C_1,k_1}\cap\mathcal{Z}_\Lambda
    \right)\\&\lesssim \#\{\gamma\in L C_1;\
    \Lambda_\ell(\gamma)\text{ in decomposition}\}
    \P\left(X_{\Lambda_\ell(\gamma),I_1+e^{-L^\eta}[-1,1]}=1\right)
    \\&\lesssim L^d|C_1|\tilde N\cdot\ell^d \\&\lesssim
    \varepsilon.
  \end{split}
\end{equation*}
The case when $k_1=0$ is then dealt with as in the proof of
Theorem~\ref{thr:2}.\\
As in the proof of Theorem~\ref{thr:2}, one shows that if
\begin{equation*}
 p+(k_1+\cdots+k_p)+(|I_1|+\cdots+|I_p|)+(|C_1|+\cdots+|C_p|)\to+\infty 
\end{equation*}
as $L\to+\infty$, then, both terms in~\eqref{eq:109} converge to $0$
in the large $L$ limit.\\
The degenerate cases having been removed, the same reasoning as in the
proof of Theorem~\ref{thr:2} yields
\begin{equation*}
  \begin{split}
    \prod_{l=1}^p\P(\{\omega;\ \#\{j;\ &\gamma_j/L\in C^-_l\text{ and }
    X_{\Lambda_\ell(\gamma_j),I^-_l}=1\}=k_l\})
    -(1-\P(\mathcal{Z}_\Lambda))\\
    \hskip3cm&\leq\P(\Omega^\Lambda_{(I_l,C_l,k_l)_{1\leq l\leq p}}
    \cap\mathcal{Z}_\Lambda) \\\hskip3cm&\leq
    \prod_{l=1}^p\P(\{\omega;\ \#\{j;\ \gamma_j/L\in C^+_l\text{ and }
    X_{\Lambda_\ell(\gamma_j),I^+_l}=1\}=k_l
    \})+(1-\P(\mathcal{Z}_\Lambda))
  \end{split}
\end{equation*}
where $C^+_l=C_l+[-\ell/L,\ell/L]^d$ and $\R^d\setminus
C^-_l=\R^d\setminus (C_l+(-\ell/L,\ell/L)^d)$.\\
Hence, the same computations as in the proof of Theorem~\ref{thr:2}
also yields
\begin{equation*}
  \begin{split}
    \prod_{l=1}^p\binom{\tilde
      N^-_l}{k_l}(p^-_l)^{k_l}(1-p^-_l)^{\tilde N^-_l-k_l}&+o(1)\leq
    \P(\Omega^\Lambda_{(I_l,C_l,k_l)_{1\leq l\leq p}}
    \cap\mathcal{Z}_\Lambda) \\
    &\leq \prod_{l=1}^p\binom{\tilde
      N^+_l}{k_l}(p^+_l)^{k_l}(1-p^+_l)^{\tilde N^+_l-k_l}+o(1)
  \end{split}
\end{equation*}
where
\begin{equation*}
  N^+_l=|C_l|\tilde N(1+O(L^{-1+\beta})),\quad N^-_l=|C_l|\tilde
  N(1+O(L^{-1+\beta})).
\end{equation*}
One concludes in the same way as in the proof of
Theorem~\ref{thr:2}. This completes the proof of
Theorem~\ref{thr:11}.\qed
\subsubsection{The proof of Theorem~\ref{thr:5}}
\label{sec:proof-theorem5}
Let $\ell_\Lambda$ be the scale defined in
section~\ref{sec:local-cent-stat} satisfying~\eqref{eq:14}. To prove
Theorem~\ref{thr:5}, it is sufficient to prove that, for $(I_j)_{1\leq
  j\leq l}$ disjoint segments of $\R$ and $(C_j)_{1\leq j\leq l}$
disjoint compact cubes in $[-c_\ell,c_\ell]$ (see~\eqref{eq:16}), one
has
\begin{multline}
  \label{eq:33}
  \pro\left(\forall 1\leq j\leq l,\ \#\left\{n;\
      \begin{aligned}
        N(E_n(\omega,\Lambda))&\in N(E_0)+\ell^{-d}_\Lambda I_j\\
        x_n(\omega)&\in\ell_\Lambda C_j
      \end{aligned}\right\}\geq k_j\right)\\\vers_{|\Lambda|\to+\infty}
  \prod_{j=1}^l\left(\sum_{k\geq
      k_j}e^{-|I_j||C_j|}\frac{(|I_j||C_j|)^{k}}{k!}\right).
\end{multline}
This is a consequence of Proposition~\ref{pro:2} and
Theorem~\ref{thr:11}. Indeed, we first pick $(C_j^\pm)$ cubes
s.t. $\overline{C_j^-}\subset \overset{\circ}{C_j}\subset
\overline{C_j}\subset \overset{\circ}{C_j^+}$ and
$\overline{C_j^+}\cap\overline{C_k^+}=\emptyset$ for $j\not=k$. For
$\ell_\Lambda$ large, the cubes $\displaystyle\left(\ell_\Lambda
  \overline{C_j^+}\right)_{1\leq j\leq l}$ are at distance at least
$R$ from one another (see (IAD)), thus, the $(H_\omega(\ell_\Lambda
C_j))_{1\leq j\leq l}$ are two by two independent so as
$(H_\omega(\ell_\Lambda C^-_j))_{1\leq j\leq l}$ and
$(H_\omega(\ell_\Lambda C^+_j))_{1\leq j\leq l}$.  Using
Proposition~\ref{pro:2}, we have that
\begin{multline*}
  \prod_{j=1}^l\pro\left(\#\left\{n;\ N(E_n(\omega,\ell_\Lambda C_j^-))\in
      N(E_0)+\ell^{-d}_\Lambda I^-_j\right\}\geq k_j\right)
  -(1-\pro(\mathcal{Z}_\Lambda))-\pro_- \\\leq\pro\left(\forall 1\leq j\leq
    l,\ \#\left\{n;\
      \begin{aligned}
        N(E_n(\omega,\Lambda))&\in N(E_0)+\ell^{-d}_\Lambda I_j\\
        x_n(\omega)&\in\ell_\Lambda C_j
      \end{aligned}\right\}\geq k_j\right)\\\leq
  \prod_{j=1}^l\pro\left(\#\left\{n;\ N(E_n(\omega,\ell_\Lambda C_j^+))\in
      N(E_0)+\ell^{-d}_\Lambda I^+_j \right\}\geq k_j\right)
  +(1-\pro(\mathcal{Z}_\Lambda))+\pro_+
\end{multline*}
where 
\begin{equation*}
  \pro_\pm:=\pro\left\{\exists j;
    \begin{aligned}
      H_\omega(\ell_\Lambda C_j)\text{ has an eigenvalue in }
      \ell^{-d}_\Lambda I_j\\\text{ with a localiz. center in }\ell_\Lambda(
      C_j^\pm\Delta C_j)
    \end{aligned}
  \right\}\lesssim \sum_{j=1}^l |C_j^\pm\Delta C_j|,
\end{equation*}
the last bound being a consequence of the Wegner estimate (W).\\
Now, using Theorem~\ref{thr:11}, we compute the asymptotics of
$\pro(\#\{n;\ N(E_n(\omega,\ell_\Lambda C_j^\pm))\in
N(E_0)+\ell^{-d}_\Lambda I^\pm_j\})$; then, we let $C_j^\pm$ tend to
$C_j$ and $I_j^\pm$ tend to $I_j$ to get the desired result. This
completes the proof of Theorem~\ref{thr:5}.\qed
\subsubsection{The proof of Theorem~\ref{thr:7}}
\label{sec:proof-theorem7}
It is sufficient to consider the case of $J$ a (non empty) segment and
$C$ a (non empty) cube. As the operators we considered are defined
with periodic boundary conditions, we can restrict ourselves to cubes
containing $0$. We can also assume $J$ is of the form $[0,a]$ or
$[-a,0]$ for some $a>0$. Hence, the sets $\tilde\ell_\Lambda C$ are
increasing and the sets $E_0+\ell_\Lambda J$ are decreasing.\\
Pick $\ell^{\pm}=(\ell^{\pm}_\Lambda)_\Lambda$
and $\tilde\ell^{\pm}=(\tilde\ell^{\pm}_\Lambda)_\Lambda$ such that
\begin{equation}
 \label{eq:34}
 \pm\frac{\ell^\pm_\Lambda-\ell_\Lambda}{\log^{1/\xi}|\Lambda|}\to+\infty,\quad
 \frac{\ell^\pm_\Lambda}{\ell_\Lambda}\to1,\quad
 \pm\frac{\tilde\ell^\pm_\Lambda-\tilde\ell_\Lambda}{\log^{1/\xi}|\Lambda|}
 \to+\infty,\quad
 \frac{\tilde\ell^\pm_\Lambda}{\tilde\ell_\Lambda}\to1.
\end{equation}
Let $\chi_A$ be the characteristic function of $A$. Compute
\begin{equation*}
  \tr(\chi_{\tilde\ell_\Lambda
    C}\car_{N^{-1}[N(E_0)+(\ell_\Lambda)^{-d}J]}(H_\omega(\Lambda)))
  =\sum_{\substack{E_n(\omega,\Lambda)\in \sigma(H_\omega(\Lambda))\\
      N(E_n(\omega,\Lambda))\in N(E_0)+\ell^{-d}_\Lambda J}}
  \|\chi_{\tilde\ell_\Lambda C}\varphi_n(\omega,\Lambda)\|^2.
\end{equation*}
If $\varphi_n(\omega,\Lambda)$ has its localization center in
$\tilde\ell_\Lambda C$, by~\eqref{eq:34} and~\eqref{eq:19}, one has
\begin{equation*}
  \|\chi_{\tilde\ell^+_\Lambda C}\varphi_n(\omega,\Lambda)\|^2
  =1+O(|\Lambda|^{-\infty}),
\end{equation*}
and if it has its localization center outside $\tilde\ell_\Lambda C$,
then
\begin{equation*}
  \|\chi_{\tilde\ell^-_\Lambda C}\varphi_n(\omega,\Lambda)\|^2
  =O(|\Lambda|^{-\infty}).
\end{equation*}
Hence, as the number of eigenvalues of $H_\omega(\Lambda)$ in
$N^{-1}[N(E_0)+(\ell_\Lambda)^{-d}J]$ is bounded by $C|\Lambda|$ for
some $C>0$, by Lemma~\ref{le:1}, we get that, for any $p>0$, for
$\Lambda$ sufficiently large, with a probability at least
$1-|\Lambda|^{-p}$, one has
\begin{multline}
  \label{eq:35}
  \begin{aligned}
    \int_{J\times C}\Xi^2_\Lambda(\xi,x;E_0,\ell,\tilde\ell)d\xi
    dx-|\Lambda|^{-p}&\leq\tr(\chi_{\tilde\ell^+_\Lambda C}
    \car_{N^{-1}[N(E_0)+(\ell_\Lambda)^{-d}J]}(H_\omega(\Lambda))) \\
    &\leq \int_{J\times C}\Xi^2_\Lambda(\xi,x;E_0,\ell,\tilde\ell)d\xi
    dx+|\Lambda|^{-p}.
  \end{aligned}
\end{multline}
Partitioning $\tilde\ell^\pm_\Lambda C$ into cubes of side length $1$
and using the covariance and~\eqref{eq:34}, for $\chi_0$ taken as in
section~\ref{sec:asymptotics}, we get that
\begin{multline*}
  (\tilde\ell_\Lambda)^d \,|C|\,\E_{\ell_\Lambda}(1+o(1))\leq
  \E(\tr(\chi_{\tilde\ell^-_\Lambda C}
  \car_{N^{-1}[N(E_0)+(\ell_\Lambda)^{-d}J]}(H_\omega(\Lambda))))\\\leq
  \E(\tr(\chi_{\tilde\ell^+_\Lambda
    C}\car_{N^{-1}[N(E_0)+(\ell_\Lambda)^{-d}J]}
  (H_\omega(\Lambda))))\leq(\tilde\ell_\Lambda)^d
  \,|C|\,\E_{\ell_\Lambda}(1+o(1)).
\end{multline*}
where $\E_{\ell_\Lambda}=\E(\tr(\chi_0\car_{N^{-1}[N(E_0)+(\ell_\Lambda)^{-d}J]}
(H_\omega(\Lambda))))$.  The computations done in
section~\ref{sec:asymptotics} show that
\begin{equation*}
  \E_{\ell_\Lambda}=(\ell_\Lambda)^{-d}|J|(1+o(1)).
\end{equation*}
Taking the expectation in~\eqref{eq:35}, we immediately obtain
\begin{gather}
  \label{eq:67}
  \esp\left(\int_{J\times C}\Xi^2_\Lambda(\xi,x;E_0,\ell,\tilde\ell)
    d\xi\right)\leq C
  \left(\frac{\tilde\ell_\Lambda}{\ell_\Lambda}\right)^d\text{ if }
  \frac{\tilde\ell_\Lambda}{\ell_\Lambda}
  \vers_{|\Lambda|\to+\infty}0,\\
  \label{eq:68}
  \left(\frac{\ell_\Lambda}{\tilde\ell_\Lambda}\right)^d
  \esp\left(\int_{J\times C}\Xi^2_\Lambda(\xi,x;E_0,\ell,\tilde\ell)
    d\xi\right)\vers_{|\Lambda|\to+\infty}|J|\cdot|C|
  \text{ if }\frac{\tilde\ell_\Lambda}{\ell_\Lambda}
  \vers_{|\Lambda|\to+\infty}+\infty.
\end{gather}
Assume now that
$\frac{\tilde\ell_\Lambda}{\ell_\Lambda}\geq|\Lambda|^\rho$.  Pick two
scales $(\ell'_\Lambda)_\Lambda$ and $(\ell''_\Lambda)_\Lambda$ such
that, for some $\rho'>0$
\begin{equation}
  \label{eq:37}
  \ell'_\Lambda\geq|\Lambda|^{\rho'},\quad
  \frac{\tilde\ell_\Lambda}{\ell'_\Lambda}\geq|\Lambda|^{\rho'}
  \quad\text{and}\quad
  \frac{\ell'_\Lambda}{\ell_\Lambda}\geq|\Lambda|^{\rho'}.
\end{equation}
Partition the cubes $\tilde\ell^\pm_\Lambda C$ into cubes of
side length asymptotic to $\ell'_\Lambda$: let
$\Gamma^\pm_\Lambda=(\ell'_\Lambda\Z^d)\cap(\tilde\ell^\pm_\Lambda C)$ and
\begin{equation*}
  \tilde\ell^\pm_\Lambda C=\bigcup_{\gamma\in\Gamma_\Lambda}
  C_{\gamma,\ell'_\Lambda}
  \quad\text{where}\quad
  C^\pm_{\gamma,\ell'_\Lambda}=\gamma+\ell'_\Lambda[-1/2,1/2]^d.
\end{equation*}
Then, we have
\begin{equation*}
  \tr(\chi_{\tilde\ell^\pm_\Lambda C}
  \car_{N^{-1}[N(E_0)+(\ell_\Lambda)^{-d}J]}(H_\omega(\Lambda)))
  =\sum_{\gamma\in\Gamma^\pm_\Lambda}
  \tr(\chi_{C^\pm_{\gamma,\ell'_\Lambda}}\car_{E_0+(\ell_\Lambda)^{-d}
      J}(H_\omega(\Lambda)))    
\end{equation*}
Thus
\begin{equation*}
  \tr^2(\chi_{\tilde\ell^\pm_\Lambda C}
  \car_{N^{-1}[N(E_0)+(\ell_\Lambda)^{-d}J]}(H_\omega(\Lambda)))
  \\=\sum_{\gamma\in\Gamma_\Lambda}
  \sum_{\gamma'\in\Gamma_\Lambda}
  T(\gamma,J,\ell'_\Lambda,\ell_\Lambda,\Lambda)
  T(\gamma',J,\ell'_\Lambda,\ell_\Lambda,\Lambda)
\end{equation*}
where
\begin{equation*}
  T(\gamma)=T(\gamma,J,\ell'_\Lambda,\ell_\Lambda,\Lambda)=
  \tr(\chi_{C^\pm_{\gamma,\ell'_\Lambda}}
  \car_{N^{-1}[N(E_0)+(\ell_\Lambda)^{-d}J]}(H_\omega(\Lambda))). 
\end{equation*}
We prove
\begin{Le}
  \label{le:8}
  If $|\gamma-\gamma'|\geq2\ell'_\Lambda$ then
  \begin{equation*}
    \left|\esp\left(T(\gamma)\cdot T(\gamma')\right)-
      \esp\left(T(\gamma)\right)\cdot
      \esp\left(T(\gamma')\right)\right|
    \leq C e^{-(\ell'_\Lambda)^{1/\xi}/C}.
  \end{equation*}
\end{Le}
\noindent Hence, we have
\begin{equation*}
  \esp\left|\tr(\chi_{\tilde\ell^\pm_\Lambda C}
    \car_{N^{-1}[N(E_0)+(\ell_\Lambda)^{-d}J]}(H_\omega(\Lambda)))-
    \esp\left(\tr(\chi_{\tilde\ell^\pm_\Lambda
        C}\car_{N^{-1}[N(E_0)+(\ell_\Lambda)^{-d}J]}(H_\omega(\Lambda)))\right)
  \right|^2\leq C(\ell'_\Lambda)^d
\end{equation*}
By~\eqref{eq:37} and~\eqref{eq:68}, we get, for some $\rho>0$,
\begin{equation}
  \label{eq:71}
  \esp\left|\left(\frac{\ell_\Lambda}{\tilde\ell_\Lambda}\right)^d
    \tr(\chi_{\tilde\ell^\pm_\Lambda C}
    \car_{N^{-1}[N(E_0)+(\ell_\Lambda)^{-d}J]}(H_\omega(\Lambda)))-
    |J|\cdot|C|
  \right|^2\leq C|\Lambda|^{-\rho}
\end{equation}
If $\tilde\ell_\Lambda/\ell'_\Lambda
\leq|\Lambda|^{-\rho}$,~\eqref{eq:67} becomes
\begin{equation}
  \label{eq:69}
\esp\left(\int_{J\times C}\Xi^2_\Lambda(\xi,x;E_0,\ell,\tilde\ell)
    d\xi\right)\leq C|\Lambda|^{-d\rho}.  
\end{equation}
Now choose the scales
$\tilde\ell^+_{\Lambda_{L^p}}=\tilde\ell_{\Lambda_{(L+1)^p}}$ and
$\ell^+_{\Lambda_{L^p}}=\ell_{\Lambda_{(L+1)^p}}$. By~\eqref{eq:70},
$\tilde\ell^+_{\Lambda_{L^p}}/\tilde\ell_{\Lambda_{L^p}}\to1$ and
$\ell^+_{\Lambda_{L^p}}/\ell_{\Lambda_{L^p}}\to1$. Moreover, the
estimates~\eqref{eq:69} and~\eqref{eq:71} hold for the pairs of scales
$(\tilde\ell^+_{\Lambda_{L^p}},\ell_{\Lambda_{L^p}})$, and
$(\tilde\ell_{\Lambda_{L^p}},\ell^+_{\Lambda_{L^p}})$. Moreover, for
$p$ large enough, they are summable. Thus, we have proved that
\begin{itemize}
\item in case (1) of Theorem~\ref{thr:7}, $\omega$ almost surely, for
  $L$ sufficiently large,
  \begin{equation*}
    \int_{J\times C}\Xi^2_{\Lambda_{L^p}}(\xi,x;E_0,\ell,\tilde\ell^+)d\xi
    dx=0.
  \end{equation*}
\item in case (2) of Theorem~\ref{thr:7}, $\omega$ almost surely,
  \begin{gather*}
    \left(\frac{\ell^+_{\Lambda_{L^p}}}{\tilde\ell_{\Lambda_{L^p}}}\right)^{-d}
    \int_{J\times C}\Xi^2_{\Lambda_{L^p}}(\xi,x;E_0,\ell,\tilde\ell^+)d\xi
    dx\vers_{|\Lambda|\to+\infty}|J|\cdot|C|,\\
    \left(\frac{\ell_{\Lambda_{L^p}}}{\tilde\ell^+_{\Lambda_{L^p}}}\right)^{-d}
    \int_{J\times C}\Xi^2_{\Lambda_{L^p}}(\xi,x;E_0,\ell^+,\tilde\ell)d\xi
    dx\vers_{|\Lambda|\to+\infty}|J|\cdot|C|
  \end{gather*}
\end{itemize}
For $L^p\leq k\leq(L+1)^p$, as the sequences $(\ell_\Lambda)_\Lambda$
and $(\tilde\ell_\Lambda)_\Lambda$ are increasing, one has
\begin{equation*}
  \begin{split}
  \int_{J\times
    C}\Xi^2_{\Lambda_{L^p}}(\xi,x;E_0,\ell,\tilde\ell^+)d\xi dx&\geq
  \int_{J\times C}\Xi^2_{\Lambda_k}(\xi,x;E_0,\ell,\tilde\ell)d\xi dx
  \\&\geq \int_{J\times
    C}\Xi^2_{\Lambda_{L^p}}(\xi,x;E_0,\ell^+,\tilde\ell)d\xi dx.    
  \end{split}
\end{equation*}
By~\eqref{eq:70}, for $L^p\leq k\leq(L+1)^p$,
$\tilde\ell^+_{\Lambda_{L^p}}\sim
\tilde\ell^+_{\Lambda_k}\sim\tilde\ell_{\Lambda_{L^p}}$ and
$\ell^+_{\Lambda_{L^p}}\sim
\ell^+_{\Lambda_k}\sim\ell_{\Lambda_{L^p}}$. Hence, we get
\begin{itemize}
\item in case (1) of Theorem~\ref{thr:7}, $\omega$ almost surely, for
  $L$ sufficiently large,
  \begin{equation*}
    \int_{J\times C}\Xi^2_{\Lambda}(\xi,x;E_0,\ell,\tilde\ell^+)d\xi
    dx=0.
  \end{equation*}
  As $\int_{J\times C}\Xi^2_{\Lambda}(\xi,x;E_0,\ell,\tilde\ell^+)d\xi
    dx$ is an integer, this implies that this integer is $0$.
\item in case (2) of Theorem~\ref{thr:7}, $\omega$ almost surely,
  \begin{gather*}
    \left(\frac{\ell_\Lambda}{\tilde\ell_\Lambda}\right)^{-d}
    \int_{J\times C}\Xi^2_\Lambda(\xi,x;E_0,\ell,\tilde\ell)d\xi
    dx\vers_{|\Lambda|\to+\infty}|J|\cdot|C|.
  \end{gather*}
\end{itemize}
This completes the proof of Theorem~\ref{thr:7}.\qed
\begin{Rem}
  \label{rem:3}
  If we don't assume that either $\tilde\ell_\Lambda/\ell'_\Lambda
  \leq|\Lambda|^{-\rho}$ or $\tilde\ell_\Lambda/\ell'_\Lambda
  \leq|\Lambda|^{-\rho}$, but merely that either tends to $0$, or we
  do not assume condition~\eqref{eq:70}
  then,~\eqref{eq:67},~\eqref{eq:35} and~\eqref{eq:71} show
  nevertheless that
  \begin{gather*}
    \esp\left(\int_{J\times C}\Xi^2_\Lambda(\xi,x;E_0,\ell,\tilde\ell)
      d\xi\right)\to0\quad\text{ if }\quad
    \frac{\tilde\ell_\Lambda}{\ell_\Lambda}
    \vers_{|\Lambda|\to+\infty}0,\\
    \esp\left|\left(\frac{\ell_\Lambda}{\tilde\ell_\Lambda}\right)^d
      \int_{J\times C}\Xi^2_\Lambda(\xi,x;E_0,\ell,\tilde\ell) d\xi -
      |J|\cdot|C| \right|^2\to0 \quad\text{ if
    }\quad\frac{\tilde\ell_\Lambda}{\ell_\Lambda}
    \vers_{|\Lambda|\to+\infty}+\infty.
  \end{gather*}
  This implies convergence in probability.
\end{Rem}
\begin{proof}[Proof of Lemma~\ref{le:8}]
  Define $\tilde C^\pm_{\gamma,\ell'_\Lambda}
  =C^\pm_{\gamma,\ell'_\Lambda} +\ell'_\Lambda[-1/2,1/2]^d$; hence, in
  view of the computations done in the proof of Lemma~\ref{le:9} and
  in section~\ref{sec:proof-proposition1}, there exists $C>0$ such
  that, for $\delta>0$, we have
  \begin{multline}
    \label{eq:72}
    \esp\left|\tr(\chi_{C_{\gamma,\ell'_\Lambda}}
      \car_{N^{-1}[N(E_0)+(\ell_\Lambda)^{-d}J]}
      (H_\omega(\Lambda)))-\tr(\chi_{C^\pm_{\gamma,\ell'_\Lambda}}
      \car_{N^{-1}[N(E_0)+(\ell_\Lambda)^{-d}J]}(H_\omega(\tilde
      C^\pm_{\gamma,\ell'_\Lambda})))\right|\\\leq C
    \left(\delta(\ell_\Lambda)^d+\delta^{-C}(\ell_\Lambda)^{dC}
      e^{-(\ell'_\Lambda)^{1/\xi}/C}\right)
  \end{multline}
  Pick $\delta=e^{-(\ell'_\Lambda)^{1/\xi}/C'}$ for some $C'>C$. As
  $|\gamma-\gamma'|\geq 2\ell'_\Lambda$, the operators
  $H_\omega(\tilde C^\pm_{\gamma,\ell'_\Lambda})$ and $H_\omega(\tilde
  C^\pm_{\gamma',\ell'_\Lambda})$ are stochastically independent of
  each other. Now, as by standard arguments of Schr{\"o}dinger operator
  theory (see e.g.~\cite{Re-Si:79}), $T(\gamma)$ and $T(\gamma')$ are
  bounded by $C(\ell'_\Lambda)^d$,~\eqref{eq:72} yields the result of
  Lemma~\ref{le:8} and completes its proof.
\end{proof}
\subsection{Study of the level spacings statistics}
\label{sec:conv-level-spac}
We will first prove Theorem~\ref{thr:1}. To do so, we first use a the
reduction constructed in Theorem~\ref{thr:vbig1} and study the
spacings for the approximated eigenvalues given by
Theorem~\ref{thr:vbig1}. Then, we derive the statistics described in
Theorem~\ref{thr:1} from those computations. Using the estimates
obtained in the proof of Theorem~\ref{thr:1}, we will prove
Theorem~\ref{thr:9}.
\subsubsection{Some preliminary considerations}
\label{sec:some-prel-cons}
We now use Theorem~\ref{thr:vbig1}. The length scale $\ell'_\Lambda$
(the localization radius) is determined by Theorem~\ref{thr:vbig1}
i.e.  $\ell'=\ell'_\Lambda=(R\log|\Lambda|)^{1/\xi}$ where
$\xi\in(0,1)$ can be chosen arbitrary by (Loc). Recall that
$\rho'\in[0,\rho/(1+d(\rho+1)))$ is fixed by
assumption~\eqref{eq:86}. \\
We first assume that the integrated density of
states of the interval $E_0+I_\Lambda$, that is, $N(E_0+I_\Lambda)$
satisfies 
\begin{equation}
  \label{eq:104}
  N(E_0+I_\Lambda)\asymp |\Lambda|^{-\alpha}
\end{equation}
for some $\alpha\in(0,1)$. When this is not the case, then by
assumption~\eqref{eq:43}, we know that $N(E_0+I_\Lambda)\gg
|\Lambda|^{-\alpha}$ for any $\alpha\in(0,1)$; thus, we can partition
$E_0+I_\Lambda$ into intervals satisfying~\eqref{eq:104} for some
$\alpha\in(0,1)$.\\
We pick the scale $\ell=\ell_\Lambda$ and so that,
\begin{equation}
  \label{eq:90}
  \ell_\Lambda\asymp N(E_0+I_\Lambda)^{-\nu}.
\end{equation}
We now show that $\nu$ can be chosen in $(0,1/d)$ so that all the
assumptions of Theorem~\ref{thr:vbig1} (in particular~\eqref{eq:93}
and~\eqref{condell}) and Lemmas~\ref{lemasympt} and~\ref{le:7} (in
particular~\eqref{eq:91}) are satisfied. Note that, when applying
Lemmas~\ref{lemasympt} or~\ref{le:7}, in~\eqref{eq:91}, the cube
$\Lambda$ is the cube $\Lambda_\ell$. \\
These requirements yield the following conditions on the exponents
\begin{equation}
  \label{eq:105}
  \begin{split}
    \frac1{1+\rho'}>d\nu&,\quad \quad
    0<1-\frac1\alpha+ \frac{\rho-\rho'}{1+\rho'}-\nu d\rho,\\
    \frac{\rho-\rho'}{1+\rho'}>\nu d(1+\rho)&, \quad \quad
    0<1-\frac1\alpha+ \nu-\frac{\rho'}{1+\rho'}.
  \end{split}
\end{equation}
As $\displaystyle\frac{1+\rho}{1+\rho'}>1>\frac{\rho-\rho'}{1+\rho'}$,
to obtain~\eqref{eq:105} for some $\nu\in(0,1/d)$ and
$\alpha\in(0,1)$, it suffices that
\begin{equation*}
   \frac{\rho-\rho'}{(1+\rho')(1+\rho)}>\frac{d\rho'}{1+\rho'}.
\end{equation*}
This is satisfied as $\rho'\in[0,\rho/(1+d(\rho+1)))$.\\
Moreover, recalling the discussion following Lemma~\ref{lemasympt}, we
want $I_\Lambda$ and $\ell_\Lambda$ to be such that
$(N(I_\Lambda)|\Lambda_L|)^{-1}\gg
N\left(I_\Lambda\right)^{(\rho-\rho')/(1+\rho')}|\Lambda_{\ell_\Lambda}|^\rho$
that is, using~\eqref{eq:104} and~\eqref{eq:90}, we need that
\begin{equation*}
 \frac1\alpha<1+\frac{\rho-\rho'}{1+\rho'}-d\nu\rho.
\end{equation*}
This condition is fulfilled by~\eqref{eq:105}. From now on, we assume
that $\ell=\ell_\Lambda$ and $N(E_0+I_\Lambda)$ satisfy~\eqref{eq:104}
and~\eqref{eq:90} for such $\alpha$ and $\nu$. Pick $\beta\in(0,1/2)$ be
such that
\begin{equation}
  \label{eq:106}
  \left(\frac1\alpha-1\right)(1+2\beta)=\frac{\rho-\rho'}{1+\rho'}-d\nu\rho.
\end{equation}
Let $\mathcal{Z}_\Lambda$ be the set of realizations for which the
conclusions of Theorem~\ref{thr:vbig1} hold. We know that, for any
$p>0$, if $L$ is sufficiently large, one has $\P(\mathcal{Z}_\Lambda)
\geq1-O(|\Lambda|^{-p})$ (see~\eqref{eq:94}).\\
Let $\tilde N$ be the number of the good cubes (i.e. cubes that
determine eigenvalues of $H_\omega(\Lambda)$) constructed in
Theorem~\ref{thr:vbig1} and $(\Lambda_\ell(\gamma_j))_{1\leq j\leq
  \tilde N}$ be those cubes. Then, $\tilde
N=|\Lambda|\ell_\Lambda^{-d}(1+o(1))$. As before, define the following
random variables:
\begin{itemize}
\item $X_{j}=X_{j}(\ell,E_0+I_\Lambda)$ is the Bernoulli
  random variable
  \begin{equation*}
    X_j=\car_{H_\omega(\Lambda_\ell(\gamma_j))\text{ has exactly one
        eigenvalue in }E_0+I_\Lambda\text{ with localization center in }
      \Lambda_{\ell-\ell'}};
  \end{equation*}
  here, $\ell'\asymp(\log|\Lambda|)^{1/\xi}\ll\ell=\ell_\Lambda$ (see
  the discussion above);
\item $\tilde E_j=\tilde E_j(\ell,E_0+I_\Lambda)$ is this
  eigenvalue conditioned on $X_j=1$.
\end{itemize}
Assume $I_\Lambda=[a_\Lambda,b_\Lambda]$ and define
\begin{equation}
  \label{eq:87}
  \xi_j=\frac{N(\tilde E_j)-N(E_0+a_\Lambda)}{N(E_0+b_\Lambda)-N(E_0+a_\Lambda)}.
\end{equation}
Note that $\xi_j$ is valued in $[0,1]$. Let $\Xi$ denote the common
distribution function of the $(\xi_j)_{1\leq j\leq k}$. It was studied
in Lemma~\ref{le:7}, the assumptions of which are satisfied (note that
the set $\Lambda$ in~\eqref{eq:91} in Lemma~\ref{le:7} is the set
$\Lambda_\ell(\gamma_j)$). The error term in~\eqref{eq:89} is then of
order $(\log|\Lambda|)^{-\beta}$ for some $\beta>0$.\\
We first study the spacings for i.i.d. copies of the random variables
$(\xi_j)_{1\leq j\leq k}$. Let $(\overline{\xi}_j)_{1\leq j\leq k}$
denote the $(\xi_j)_{1\leq j\leq k}$ ordered increasingly and define
\begin{equation}
  \label{eq:44}
  DLS_\xi(x,k;E_0+I_\Lambda,\omega,\Lambda)=\frac1{k-1}\#\{1\leq j\leq k;\
    \overline{\xi}_{j+1}-\overline{\xi}_j>x/k\}.
\end{equation}
Our main technical result is
\begin{Le}
  \label{le:2}
  Pick $E_0\in I$ such that~\eqref{eq:86} be satisfied. Pick
  $(I_\Lambda)_\Lambda$ be intervals and $(\ell_\Lambda)$ length
  scales such that~\eqref{eq:104} and~\eqref{eq:90} be satisfied for
  $(\nu,\alpha)$ satisfying~\eqref{eq:105}. Define $\beta$
  by~\eqref{eq:106} and let $N_\Lambda:=N(E_0+I_\Lambda)|\Lambda|$ and
  $K_\Lambda=N_\Lambda^{\beta}$.\\
  Then, there exists $C>0$ such that, for $|\Lambda|$ sufficiently
  large, for $N_\Lambda K^{-1}_\Lambda\leq k\leq N_\Lambda K_\Lambda$,
  one
  \begin{equation}
    \label{eq:45}
    \sup_{K_\Lambda^{-1}\leq x\leq K_\Lambda}
    \esp\left(\left|DLS_\xi(x,k;E_0+I_\Lambda,\omega,\Lambda)-
        D(k,E_0+I_\Lambda,\Lambda)\right|^2\right)\leq \frac Ck.
  \end{equation}
  where we have defined
  \begin{equation}
    \label{eq:46}
    D(k,E_0+I_\Lambda,\Lambda):=\int_0^1(1-\Xi(y+x/k)+\Xi(y))^{k-1}d\Xi(y).
  \end{equation}
  Let $(J_\Lambda)_\Lambda$ be a sequence of intervals such that
  $\displaystyle\sup_{E\in J_\Lambda}|E|\to0$ as
  $|\Lambda|\to+\infty$. Then
  \begin{equation}
    \label{eq:6}
    \sup_{N_\Lambda K^{-1}_\Lambda\leq k\leq
      N_\Lambda K_\Lambda}
    \sup_{\substack{I_\Lambda\subset J_\Lambda\\I_\Lambda
        \text{ as in Lemma~\ref{le:2}}}}
    \left|D(k,E_0+I_\Lambda,\Lambda)-e^{-x}\right|
    \vers_{|\Lambda\|\to+\infty}0.
  \end{equation}
\end{Le}
\begin{proof}
  To analyze the spacings of the $(\xi_j)_{1\leq j\leq k}$, we use the
  computations of section 7 in~\cite{MR0216622} (see in particular
  (7.3)) that yield
  \begin{gather*}
      \int_0^1(1-\Xi(y+x/k)+\Xi(y))^{k-1}d\Xi(y)
      =\esp(DLS_\xi(x,k;E_0+I_\Lambda,\omega,\Lambda))
      =D(k,E_0+I_\Lambda,\Lambda)
    \\\intertext{and}
    \begin{split}
      &\esp\left([DLS_\xi(x,k;E_0+I_\Lambda,\omega,\Lambda)]^2\right)=O\left(\frac1k\right)+
      \\&+2\int_{\R}\int^{+\infty}_{y+x/k}
      (1-\Xi(y+x/k)+\Xi(y)-\Xi(z+x/k)+\Xi(z))^{k-2}d\Xi(z)d\Xi(y).
    \end{split}
  \end{gather*}
  Fix $\nu\in(\xi,1)$. By~\eqref{eq:104} and~\eqref{eq:90} as
  $\ell'=\ell'_\Lambda\asymp(\log|\Lambda|)^{1/\xi}$, one has
  $K_\Lambda^{-1}\gg \ell^d_\Lambda e^{-(\ell'_\Lambda)^{\nu}}$ for
  $|\Lambda|$ large. Thus, Lemma~\ref{le:7} and~\eqref{eq:106} yield
  that, for $y-x\geq (N_\Lambda K_\Lambda)^{-1}$, one has
  \begin{equation}
    \label{eq:48}
    \Xi(y)-\Xi(x)=(y-x) \left(1+o\left(K_\Lambda^{-1}\right)\right)
  \end{equation}
  Hence, for $K_\Lambda^{-1}\leq x\leq K_\Lambda$ and $N_\Lambda
  K^{-1}_\Lambda\leq k\leq N_\Lambda K_\Lambda$, one has
  \begin{equation*}
    \begin{split}
      &\esp(DLS^2_\xi(x,k;I_\Lambda,\omega,\Lambda))\\&=\int_0^1\int_0^1
      (1-\Xi(y+x/k)+\Xi(y)-\Xi(z+x/k)+\Xi(z))^{k-2}d\Xi(z)d\Xi(y)\\
      &- 2\int_{\R}\int^{y+x/k}_y
      (1-\Xi(y+x/k)+\Xi(y)-\Xi(z+x/k)+\Xi(z))^{k-2}d\Xi(z)d\Xi(y)+O\left(\frac1k\right)
      \\&=\int_0^1\int_0^1
      (1-\Xi(y+x/k)+\Xi(y)-\Xi(z+x/k)+\Xi(z))^{k-2}d\Xi(z)d\Xi(y)
      +O\left(\frac1k\right).
    \end{split}
  \end{equation*}
  Compute
  \begin{equation*}
    \begin{split}
      &(1-\Xi(y+x/k)+\Xi(y)-\Xi(z+x/k)+\Xi(z))
      \\&=(1-\Xi(y+x/k)+\Xi(y))(1-\Xi(z+x/k)+\Xi(z))
      \\&\hskip3cm-(\Xi(z+x/k)-\Xi(z))(\Xi(y+x/k)-\Xi(y))\\&
      =(1-\Xi(y+x/k)+\Xi(y))(1-\Xi(z+x/k)+\Xi(z))+O(k^{-2}).
    \end{split}
  \end{equation*}
  Hence, plugging this into the previous formula, we get
  \begin{equation*}
    \begin{split}
      &\esp(DLS^2_\xi(x,k;I_\Lambda,\omega,\Lambda))= \left[\int_0^1
        (1-\Xi(y+x/k)+\Xi(y))^{k-2}d\Xi(y)
      \right]^2+O\left(\frac1k\right)\\&= \left[\int_0^1
        (1-\Xi(y+x/k)+\Xi(y))^{k-1}
        \left(1+O\left(\frac1k\right)\right)d\Xi(y)
      \right]^2+O\left(\frac1k\right)
      \\&=D(k,E_0+I_\Lambda,\Lambda)^2+O\left(\frac1k\right)
    \end{split}
  \end{equation*}
  This completes the proof of~\eqref{eq:45}.\\
  Fix $E_0$ such that~\eqref{eq:86} be satisfied and let us
  prove~\eqref{eq:6}. For $I_\Lambda\subset J_\Lambda$, by
  Lemma~\ref{le:7}, for $x\in[0,1]$, we have,
  \begin{equation*}
    \sup_{y\in[0,1-x/k]}\left|\Xi(y+x/k)-\Xi(y)-\frac{x}{k}\right|
    \leq \frac{x}{k}\alpha_k\quad\text{and}\quad
    \sup_{y\in[1-x/k,1]}\left|\Xi(y+x/k)-\Xi(y)\right|\leq
    \frac{x(1+\alpha_k)}k
 \end{equation*}
 where $\displaystyle \sup_{N_\Lambda K^{-1}_\Lambda\leq k\leq
   N_\Lambda K_\Lambda} \sup_{\substack{I_\Lambda\subset
     J_\Lambda\\I_\Lambda \text{ as in Lemma~\ref{le:2}}}}
 |\alpha_k|\to0$ as $|\Lambda|\to+\infty$.\\
 Hence, for $k$ large,
 \begin{equation*}
    \left|D(k,E_0+I_\Lambda,\Lambda)-e^{-x}\right|\leq
    \int_0^{1-x/k}\left|e^{C\alpha_k}-1\right|d\Xi(y)+
    \int_{1-x/k}^1e^{Cx(1+\alpha_k)/k}d\Xi(y)\leq 
    C\left(\alpha_k+\frac1k\right).
 \end{equation*}
 This completes the proof of~\eqref{eq:6} and, thus, of
 Lemma~\ref{le:2}.
\end{proof}
\subsubsection{The proof of Theorem~\ref{thr:1}}
\label{sec:proof-theorem1}
By a classical result (see e.g.~\cite{MR1492447}), as the functions we
are considering are monotonous and as $x\mapsto e^{-x}$ is continuous
on $[0,+\infty)$, it suffices to prove the almost sure pointwise
convergence of $DLS(x;E_0+I_\Lambda,\omega,\Lambda)$ to $e^{-x}$ (for
$x>0$).\\
Fix $\xi'<\xi$. Pick $I_\Lambda$ as in Theorem~\ref{thr:1}. We
distinguish two cases:
\begin{enumerate}
\item if $I_\Lambda$ satisfies~\eqref{eq:104} for $\alpha$ chosen as
  in section~\ref{sec:some-prel-cons} (see the discussion in the
  beginning of that section), we can apply Theorem~\ref{thr:vbig1} and
  Lemma~\ref{le:2} to $I_\Lambda$ itself. We then set $K_\Lambda=1$
  and $I_{1,\Lambda}=I_\Lambda$;
\item if not, as already mentioned, we cover $I_\Lambda$ with disjoint
  intervals $(I_{k,\Lambda})_{1\leq k\leq K_\Lambda}$
  satisfying~\eqref{eq:104} of such that, for each $I_{k,\Lambda}$, we
  can apply Theorem~\ref{thr:vbig1} and Lemma~\ref{le:2}; as there are
  at most $o(|\Lambda|^{\alpha})$ such intervals (recall that $N$ is
  the distribution function of an absolutely continuous measure), the
  probability estimate for the set of good configurations
  $\mathcal{Z}_\Lambda$ (i.e. those for which one has the description
  given by Theorem~\ref{thr:vbig1} and Lemma~\ref{le:2} for all the
  intervals $(I_{k,\Lambda})_{1\leq k\leq K_\Lambda}$) still satisfies
  $\P(\mathcal{Z}_\Lambda)\geq1-O(|\Lambda|^{-p})$ for some arbitrary
  $p>0$.
\end{enumerate}
This yields
\begin{equation}
  \label{eq:42}
  \left|DLS(x;E_0+I_\Lambda,\omega,\Lambda)-\sum_{k=1}^{K_\Lambda}
    DLS(x;E_0+I_{k,\Lambda},\omega,\Lambda)
    \frac{N_\Lambda(E_0+I_{k,\Lambda})}{N_\Lambda(I_\Lambda)}\right|
  \leq\frac{K_\Lambda}{N_\Lambda(I_\Lambda)}
\end{equation}
where we recall that $N_\Lambda(I)$ is the (random) number of
eigenvalues of $H_\omega(\Lambda)$ in $I$.\\
We first study $DLS(x;E_0+I_{k,\Lambda},\omega,\Lambda)$. \\
By~\eqref{eq:4},~\eqref{eq:26} and the approximation of the
eigenvalues given by Theorem~\ref{thr:vbig1}, for
$\omega\in\mathcal{Z}_\Lambda$, if $J=\#\{1\leq j\leq \tilde N;
X_j=1\}$ then we have
\begin{equation}
  \label{eq:49}
  D^-(\omega,E_0+I_{k,\Lambda},\Lambda)+\alpha_\Lambda\geq
  DLS(x;E_0+I_{k,\Lambda},\omega,\Lambda)
  \geq D^+(\omega,E_0+I_{k,\Lambda},\Lambda)-\alpha_\Lambda
\end{equation}
where
\begin{itemize}
\item $\displaystyle D^{\pm}(\omega,E_0+I_{k,\Lambda},\Lambda;J)
  :=\frac{\#\{1\leq j\leq J;\
    \overline{\xi}_{j+1}-\overline{\xi}_j\geq
    x/N_{k,\Lambda}\pm|\Lambda|^{-2}\}}{N_{k,\Lambda}}$;
\item $(\xi_j)_j$ are the renormalized eigenvalues defined
  in~\eqref{eq:87} at the beginning of section~\ref{sec:conv-level-spac}
  for the energy interval $E_0+I_{k,\Lambda}$;
\item $N_{k,\Lambda}=N(E_0+I_{k,\Lambda})\,|\Lambda|$ and
  $\alpha_\Lambda\to0$ as $|\Lambda|\to+\infty$.
\end{itemize}
Define the random variables
\begin{equation*}
  D^{\pm}(\omega,E_0+I_{k,\Lambda},\Lambda)=D^{\pm}(\omega,E_0+I_{k,\Lambda},\Lambda;
  \#\{1\leq j\leq \tilde N; X_j=1\})\,\car_{\mathcal{Z}_\Lambda}.
\end{equation*}
We prove
\begin{Le}
  \label{le:6}
   One has
  \begin{equation}
    \label{eq:62}
    \sup_{1\leq k\leq K_\Lambda}\left[
      (N_{\Lambda}^k)^{1/4}\cdot\E\left(|D^{\pm}(\omega,E_0+I_{k,\Lambda},\Lambda)-
        D(N_{k,\Lambda},E_0+I_{k,\Lambda},\Lambda)|^2\right)\right]\leq
    C .
  \end{equation}
\end{Le}
\noindent Before proving Lemma~\ref{le:6}, let us use it to complete
the proof of Theorem~\ref{thr:1}. We first prove the almost sure
convergence for a subsequence.\\
For $\Lambda=\Lambda_{L^{5/\rho}}$, by the first condition
in~\eqref{eq:43}, summing the estimate of Lemma~\ref{le:6} for $1\leq
k\leq K_{\Lambda_{L^{5/\rho}}}$, we get that
\begin{equation*}
  \sum_{k=1}^{K_{\Lambda_{L^{5/\rho}}}}\sum_{L\geq1}
  \E\left(\left|D^{\pm}(\omega,E_0+I_{k,\Lambda_{L^{5/\rho}}},
      \Lambda_{L^{5/\rho}})
      -D(N_{\Lambda_{L^{5/\rho}}},E_0+I_{k,\Lambda_{L^{5/\rho}}},
      \Lambda_{L^{5/\rho}})\right|^2\right)<+\infty.
\end{equation*}
By~\eqref{eq:6} of Lemma~\ref{le:2}, as $E_0\in I$
satisfying~\eqref{eq:86}, we know that
\begin{equation}
  \label{eq:63}
  \sup_{1\leq k\leq K_\Lambda}|D(N_\Lambda,E_0+I_{k,\Lambda},
  \Lambda)- e^{-x}|\vers_{|\Lambda|\to+\infty}0.
\end{equation}
Hence, $\omega$-almost surely, 
\begin{equation*}
  \sup_{1\leq k\leq K_{\Lambda_{L^{5/\rho}}}}
  |D^{\pm}(\omega,E_0+I_{k,\Lambda_{L^{5/\rho}}},\Lambda_{L^{5/\rho}})
  -e^{-x}|\vers_{L\to+\infty}0.
\end{equation*}
Thus, by~\eqref{eq:49}, $\omega$-almost surely
\begin{equation}
  \label{eq:57}
  \sup_{1\leq k\leq K_{\Lambda_{L^{5/\rho}}}}
  |DLS(x;E_0+I_{k,\Lambda},\omega,\Lambda)-e^{-x}|\vers_{L\to+\infty}0.
\end{equation}
Plugging this into~\eqref{eq:42}, we obtain that, $\omega$-almost
surely
\begin{equation*}
  DLS(x;E_0+I_{\Lambda_{L^{5/\rho}}},\omega,\Lambda_{L^{5/\rho}})
  \vers_{L\to+\infty}e^{-x}.
\end{equation*}
To derive the almost sure convergence of
$(DLS(x;E_0+I_\Lambda,\omega,\Lambda))_\Lambda$, we use
\begin{Le}
  \label{le:4}
  Fix $p>0$ and $r>0$. Let $(I_\Lambda)_\Lambda$ be as
  Theorem~\ref{thr:1} or $I_\Lambda=J$, $J$ as in
  Theorem~\ref{thr:9}. Recall that
  $N_{\Lambda_L}=N_{\Lambda_L}(I_L)$. \\
  Then, $\omega$-almost surely, for $L$ sufficiently large, if
  $L'\in[L^p,(L+1)^p]$, except for at most $o(N_{\Lambda_{L^p}})$ of
  them, all the eigenvalues of $H_\omega(\Lambda_{L'})$ in
  $E_0+I_{\Lambda_{L'}}$ are at a distance less than $L^{-r}$ to an
  eigenvalue of $H_\omega(\Lambda_{L^p})$ in $E_0+I_{\Lambda_{L^p}}$.
\end{Le}
\noindent Pick $p=5/\rho$ and $r>d$. For $L^p\leq L'\leq (L+1)^p$, by
assumption~\eqref{eq:43} on $I_\Lambda$, one has
$N(E_0+I_{\Lambda_{L'}})=N(E_0+I_{\Lambda_{L^p}}) (1+o(1))$; so
$N_{\Lambda_{L'}}=N_{\Lambda_{L^p}}(1+o(1))$.  Lemma~\ref{le:4} then
implies that, $\omega$-almost surely, for $L^p\leq L'\leq (L+1)^p$ and
$L$ sufficiently large, one has
\begin{equation*}
  DLS(x;I_{\Lambda_{L'}},\omega,\Lambda_{L'})=
  DLS(x;I_{\Lambda_{L^p}},\omega,\Lambda_{L^p})+o(1).
\end{equation*}
Hence, $\omega$-almost surely, $DLS(x;I_{\Lambda},\omega,\Lambda)$
converges to $e^{-x}$ when $|\Lambda|\to+\infty$. This completes the
proof of Theorem~\ref{thr:1}.\qed
\begin{Rem}
  \label{rem:2}
  If one does not assume one of the conditions in~\eqref{eq:43}, the
  estimate~\eqref{eq:62} and Lemma~\ref{le:4} are not sufficient
  anymore to obtain the almost sure convergence for
  $DLS(x;I_{\Lambda},\omega,\Lambda)$. Nevertheless, plugging the
  estimates~\eqref{eq:63},~\eqref{eq:62} and~\eqref{eq:49}
  into~\eqref{eq:42}, we see that
  \begin{equation*}
    \esp(|DLS(x;I_{\Lambda},\omega,\Lambda)-e^{-x}|)
    \vers_{|\Lambda|\to+\infty}0.
  \end{equation*}
  This guarantees the convergence in probability.
\end{Rem}
\begin{proof}[Proof of Lemma~\ref{le:6}]
  We now fix $1\leq k\leq K_\Lambda$ and, to alleviate the notations,
  we will write $I_\Lambda:=I_{k,\Lambda}$ and
  $N_\Lambda=N_{k,\Lambda}=N(I_\Lambda)|\Lambda|$. The interval
  $I_\Lambda$ satisfies the requirements to apply
  Theorem~\ref{thr:vbig1} and Lemma~\ref{le:2} (see the discussion in
  the beginning of section~\ref{sec:some-prel-cons}).\\
  As $D^\pm$ is bounded by $C|\Lambda|$, as the $(X_j)_j$ are
  i.i.d. and as $\P(\mathcal{Z}_\Lambda)\leq|\Lambda|^{-p}$ for any
  $p>0$, we compute
  \begin{multline*}
    \E\left(|D^{\pm}(\omega,I_\Lambda,\Lambda)-
      D(N_\Lambda,E_0+I_\Lambda,\Lambda)|^2\right)
    \\\leq \frac1{|\Lambda|^p}+ \sum_{k=0}^{\tilde N}\binom{\tilde
      N}{k}\P(X_1=1)^{k}(1-\P(X_1=1))^{\tilde
      N-k}\cdot\\\hskip2cm\cdot
    \E\left(|DLS_\xi(x\pm|\Lambda|^{-2},k;E_0+I_\Lambda,\omega,\Lambda)
      -D(N_\Lambda,E_0+I_\Lambda,\Lambda)|^2\right).
  \end{multline*}
  Recall that $N_\Lambda=N(I_\Lambda)|\Lambda|$. As
  $\P(X_1=1)=N(I_\Lambda)\ell_\Lambda^d(1+o(1))$ and $\tilde
  N=|\Lambda|\ell_\Lambda^{-d}(1+o(1))$, for any $\varepsilon>0$,
  there exists $\delta>0$ such that
  \begin{equation}
    \label{eq:41}
    \sum_{|k-N_\Lambda|\geq N^{1/2+\varepsilon}_\Lambda}
    \binom{\tilde N}{k}\P(X_1=1)^{k}(1-\P(X_1=1))^{\tilde
      N-k}\leq e^{-N^\delta_\Lambda}.
  \end{equation}
  Hence, by the first condition in~\eqref{eq:43}, one gets
  \begin{multline}
    \label{eq:36}
    \E\left(|D^{\pm}(\omega,I_\Lambda,\Lambda)-D(N_\Lambda,E_0+I_\Lambda,\Lambda)|^2\right)
    \\ \leq \frac1{|\Lambda|^p}+e^{-N^{\delta}_\Lambda/2}
    +2\sum_{|k-N_\Lambda|<
      N^{1/2+\varepsilon}_\Lambda}|D(k,E_0,\Lambda)
    -D(N_\Lambda,E_0+I_\Lambda,\Lambda)|^2\\
    +2\sum_{|k-N_\Lambda|< N^{1/2+\varepsilon}_\Lambda}
    \E\left(|DLS_\xi(x\pm|\Lambda|^{-2},k;E_0+I_\Lambda,\omega,\Lambda)
      -D(k,E_0,\Lambda)|^2\right).
  \end{multline}
  For $|k-N_\Lambda|< N^{1/2+\varepsilon}_\Lambda$, let us estimate
  $D(k,E_0,\Lambda) -D(N_\Lambda,E_0+I_\Lambda,\Lambda)$ that is
  \begin{multline}
    \label{eq:47}
    D(k,E_0,\Lambda) -D(N_\Lambda,E_0+I_\Lambda,\Lambda)\\=
    \int_0^1\left((1-\Xi(y+x/k)+\Xi(y))^{k-1}-
      (1-\Xi(y+x/N_\Lambda)+\Xi(y))^{N_\Lambda-1}\right)d\Xi(y) .
  \end{multline}
  Using~\eqref{eq:48}, for $|k-N_\Lambda|<
  N^{1/2+\varepsilon}_\Lambda$, compute
  \begin{equation*}
    \begin{split}
      &\left|(1-\Xi(y+x/k)+\Xi(y))^{k-1}-
        (1-\Xi(y+x/N_\Lambda)+\Xi(y))^{N_\Lambda-1}\right|\\
      &=(1-\Xi(y+x/k)+\Xi(y))^{k-1}\cdot\\
      &\hskip1cm\cdot\left|1- \left(\frac{1-\Xi(y+x/N_\Lambda)+\Xi(y)}
          {1-\Xi(y+x/k)+\Xi(y)}\right)^{N_\Lambda}
        (1-\Xi(y+x/N_\Lambda)+\Xi(y))^{N_\Lambda-k}\right| \\ &\leq C
      \left(1-e^{-C(N_\Lambda-k)/N_\Lambda}e^{-C(N_\Lambda-k)/k}\right)
      \leq C N_\Lambda^{-1/2+\varepsilon}.
    \end{split}
  \end{equation*}
  Hence, plugging this into~\eqref{eq:47}, one obtains
  \begin{equation*}
    |D(k,E_0,\Lambda) -D(N_\Lambda,E_0+I_\Lambda,\Lambda)|\leq C
    N_\Lambda^{-1/2+\varepsilon}.
  \end{equation*}
  Plugging this and the result of Lemma~\ref{le:2} into~\eqref{eq:36},
  we obtain
  \begin{equation*}
    \E\left(|D^{\pm}(\omega,I_\Lambda,\Lambda)-
      D(N_\Lambda,E_0+I_\Lambda,\Lambda)|^2\right)\leq
    C N_\Lambda^{-1/2+3\varepsilon}.
  \end{equation*}
  We pick $\varepsilon=1/12$ to completes the proof of
  Lemma~\ref{le:6}.
\end{proof}
\begin{proof}[Proof of Lemma~\ref{le:4}] Pick $\alpha\in(0,1)$ and
  $q>1$. By Lemma~\ref{lemcenter}, with a probability at least
  $1-L^{-q}$, the eigenvalues of $H_\omega(\Lambda_{L'})$ in
  $E_0+I_{\Lambda_{L'}}$ that are at a distance more than $L^{-r}$ to
  an eigenvalue of $H_\omega(\Lambda_{L^p})$ in
  $E_0+I_{\Lambda_{L^p}}$ fall into two categories:
  \begin{enumerate}
  \item either they belong to $E_0+(I_{\Lambda_{L^p}}\setminus
    I_{\Lambda_{L'}})$ which may be empty,
  \item or they have a localization center that belongs e.g. to the
    annulus $\Lambda_{L'}\setminus \Lambda_{L^p-L^\alpha}$.
  \end{enumerate}
  The numbers of eigenvalues in the first category is bounded by
  $o(N_{\Lambda_{L^p}})$ as a consequence of the second part of
  assumption~\eqref{eq:43}. The numbers of eigenvalues in the second
  category is bounded by $o(N_{\Lambda_{L^p}})$. The
  Borel-Cantelli lemma then implies Lemma~\ref{le:4}.
\end{proof}
\subsubsection{The proof of Theorem~\ref{thr:9}}
\label{sec:proof-theorem6}
The proof is analogous to that of Theorem~\ref{thr:1}. As in that
proof, it suffices to prove the almost sure pointwise convergence of
$DLS'(x;J,\omega,\Lambda)$ to $g_{\nu,J}(x)$ for $x>0$.\\
Fix $\xi'<\xi$. Pick $J=[a,b]$ as in Theorem~\ref{thr:9}. We can then
cover it with disjoint intervals $(I_{j,\Lambda})_{1\leq j\leq
  J_\Lambda}$ of length $|J||\Lambda|^{-\alpha}$ (here $\alpha$ is
chosen as in the beginning of section~\ref{sec:some-prel-cons}) and
containing the point $e_{j,\Lambda}:=a+(j-1/2)|J||\Lambda|^{-\alpha}$
(so that $J_\Lambda\asymp|\Lambda|^{\alpha}$) and such that, for each
$I_{j,\Lambda}$, we can apply Lemma~\ref{le:2}. This yields
\begin{equation}
  \label{eq:61}
  \left|DLS'(x;J,\omega,\Lambda)-
    \sum_{j=1}^{J_\Lambda}DLS\left(\nu(e_{j,\Lambda})\,x/N(J);
      I_{j,\Lambda},\omega,\Lambda\right)
    \frac{N_\Lambda(I_{j,\Lambda})}{N_\Lambda(J)}\right|\leq
  \frac{CJ_\Lambda}{|\Lambda|}
\end{equation}
where $DLS'$ is defined in~\eqref{eq:25}.\\
Using the uniform continuity of $\nu$, the same computations as in the
proof of Theorem~\ref{thr:1} yield the following analogue
of~\eqref{eq:57}: for all $\varepsilon>0$,
\begin{equation*}
  \sup_{\substack{1\leq j\leq J_\Lambda
      \\\nu(e_{j,\Lambda})\geq\varepsilon\\
      \Lambda=\Lambda_{L^{5/\rho}}}}
  \left|DLS\left(\nu(e_{j,\Lambda})\,x/N(J);
      I_{j,\Lambda},\omega,\Lambda\right)-
    e^{-\nu(e_{j,\Lambda})x/N(J)}\right|\vers_{L\to+\infty}0,
  \quad\omega-\text{a.s.}.
\end{equation*}
The large deviation principle for the eigenvalue counting function,
Theorem~\ref{thr:16}, ensures that, $\omega$ almost surely, for
$|\Lambda|$ large, $N_\Lambda(J)\geq |\Lambda||J|/2$ and
\begin{equation*}
  \lim_{L\to+\infty}\sup_{\substack{1\leq j\leq
      J_\Lambda\\\nu(e_{j,\Lambda})\geq\varepsilon}}
  \left|\frac{N_\Lambda(I_{j,\Lambda})}
    {N_\Lambda(J)}\cdot\frac{N(J)}{\nu(e_{j,\Lambda})|I_{j,\Lambda}|}
    -1\right|=0.
\end{equation*}
Moreover, using the uniform continuity of $\nu$ on $J$, we have that
for any $\delta>0$, there exists $\varepsilon>0$ such that, for $L$
sufficiently large, one has
\begin{equation*}
  \sup_{\substack{1\leq j\leq
      J_\Lambda\\\nu(e_{j,\Lambda})\leq\varepsilon}}
  \left|\frac{N_\Lambda(I_{j,\Lambda})}
    {|I_{j,\Lambda}|\,N_\Lambda(J)}\right|\leq \delta.
\end{equation*}
Hence, by~\eqref{eq:61}, $\omega$-almost surely, there exists $C>0$
such that $\varepsilon>0$, one has
\begin{gather*}
  \sum_{\substack{1\leq j\leq J_\Lambda
      \\\nu(e_{j,\Lambda})\leq\varepsilon\\
      \Lambda=\Lambda_{L^{5/\rho}}}}DLS\left(\nu(e_{j,\Lambda})\,x/N(J);
    I_{j,\Lambda},\omega,\Lambda\right)
  \frac{N_\Lambda(I_{j,\Lambda})}{N_\Lambda(J)}\leq \delta
  \sum_{\substack{1\leq j\leq J_\Lambda
      \\\nu(e_{j,\Lambda})\leq\varepsilon\\
      \Lambda=\Lambda_{L^{5/\rho}}}}|I_{j,\Lambda}|\leq C\delta,
  \\ \intertext{thus}
  \limsup_{\substack{L\to+\infty\\\Lambda=\Lambda_{L^{5/\rho}}}}
  \left|DLS'(x;J,\omega,\Lambda)-\sum_{\substack{1\leq j\leq J_\Lambda
        \\\nu(e_{j,\Lambda})\geq\varepsilon\\
        \Lambda=\Lambda_{L^{5/\rho}}}}\frac{\nu(e_{j,\Lambda})
      |I_{j,\Lambda}|}{N(J)}e^{-\nu(e_{j,\Lambda})x}\right|\leq
  C\delta.
\end{gather*}
On the other hand, as $\nu$ is continuous, non negative on $J$ and
$x>0$, one has, for any $\delta>0$, taking $\varepsilon>0$
sufficiently small, one has
\begin{equation*}
  \left|\sum_{\substack{1\leq j\leq J_\Lambda
        \\\nu(e_{j,\Lambda})<\varepsilon}}\frac{\nu(e_{j,\Lambda})
      |I_{j,\Lambda}|}{N(J)}e^{-\nu(e_{j,\Lambda})x/N(J)}\right|
  \leq \delta
\end{equation*}
and 
\begin{equation*}
  \begin{split}
    \lim_{L\to+\infty}\sum_{\substack{1\leq j\leq J_\Lambda\\
        \Lambda=\Lambda_{L^{5/\rho}}}} \frac{\nu(e_{j,\Lambda})
      |I_{j,\Lambda}|}{N(J)}e^{-\nu(e_{j,\Lambda})\,x/N(J)}
    &=\frac{|J|}{N(J)}\int_0^1e^{-\nu(a+|J|y)\,x/N(J)}\nu(a+|J|y)dy\\
    &=\frac1{N(J)}\int_J e^{-\nu(\lambda)\,x/N(J)}\nu(\lambda)d\lambda.
  \end{split}
\end{equation*}
Thus, for $\delta>0$, we get
\begin{equation*}
  \limsup_{\substack{L\to+\infty\\\Lambda=\Lambda_{L^{5/\rho}}}}
  \left|DLS'(x;J,\omega,\Lambda)-\int_Je^{-\nu_J(\lambda)x}
    \nu_J(\lambda)d\lambda\right|\leq\delta.
\end{equation*}
Hence, letting $\delta$ tend to $0$, we get that, $\omega$ almost
surely,
\begin{equation*}
  \lim_{\substack{L\to+\infty\\\Lambda=\Lambda_{L^{5/\rho}}}}
  DLS'(x;J,\omega,\Lambda)=\int_Je^{-\nu_J(\lambda)x}\nu_J(\lambda)d\lambda.
\end{equation*}
To complete the proof of Theorem~\ref{thr:9}, we use Lemma~\ref{le:4}
as in the proof of Theorem~\ref{thr:1}.\qed
\subsection{Study of the localization centers statistics}
\label{sec:conv-local-cent-1}
We prove Theorems~\ref{thr:DCS}.\\ As in the proofs of
Theorems~\ref{thr:1} and~\ref{thr:9}, it suffices to prove the simple
convergence in~\eqref{eq:82} to obtain the uniform convergence as we
are dealing with non increasing functions and the limit is
continuous. Thus, pick $s>0$ and to prove~\eqref{eq:82}, it suffices
to prove
\begin{equation}
  \label{eq:76}
  \esp\left(\left|DCS(s;I_\Lambda,\Lambda,\omega)-e^{-s^d} 
    \right|^2\right)\vers_{\Lambda\nearrow\R^d} 0.
\end{equation}
We apply Theorem~\ref{thr:vbig1} to the cube $\Lambda$ and the
interval $I_\Lambda$ satisfying~\eqref{eq:73}. We keep the same
notations as in Theorem~\ref{thr:vbig1}. Let $\Gamma$ denote the set
of centers $\gamma_i$ constructed in Theorem~\ref{thr:vbig1}. Let
$\Gamma_b$ denote the set of $\gamma\in\Gamma$ that do not not satisfy
(1), (2) and (3) of Theorem~\ref{thr:vsmall1}. Theorem~\ref{thr:vbig1}
states that, for $\omega\in\mathcal{Z}_\Lambda$, one has
\begin{equation*}
  \#\Gamma_b=|\Lambda|\ell^{-d}O\left[(|I_\Lambda|\ell^d)^{1+\rho} 
     + N(E_0+I_\Lambda) \ell^{d-1}\ell'\right]  
\end{equation*}
where $\ell$ and $\ell'$ are defined in Theorem~\ref{thr:vbig1}.\\
By~\eqref{eq:73}, one has
\begin{equation}
  \label{eq:79}
  N(E_0+I_\Lambda)^{-1/d} \gg\ell.
\end{equation}
Thus, if we define
\begin{equation*}
  \Gamma'=\{\gamma\in\Gamma;\
  \text{dist}(\gamma,\Gamma_b)\geq 10 s(N(E_0+I_\Lambda))^{-1/d}\}  
\end{equation*}
then 
\begin{equation}
  \label{eq:38}
  \#\Gamma'=O(\#\Gamma_b\cdot N(E_0+I_\Lambda)^{-1}\ell^{-d})
  =o(|\Lambda|\ell^{-d}).
\end{equation}
For $\gamma\in\Gamma$, define the random variable $X_\gamma$ to be
equal
\begin{itemize}
\item to 1 if $H_\omega(\Lambda_\ell(\gamma))$ has an eigenvalue in
  $I_\Lambda$, and, for all $\gamma'\in\Gamma$ such that
  $|\gamma'-\gamma|\leq (N(E_0+I_\Lambda))^{-1/d}s$,
  $H_\omega(\Lambda_\ell(\gamma))$ has no eigenvalue in $I_\Lambda$,
\item to 0 if not.
\end{itemize}
Then, by~\eqref{eq:79},~\eqref{eq:38} and the estimate given in
Theorem~\ref{thr:vbig1} on the number of eigenvalues of
$H_\omega(\Lambda)$ not associated to an eigenvalue for a cube
$(\Lambda_\ell(\gamma))_{\gamma\in\Gamma}$, we get that, for any
$\varepsilon>0$, for $|\Lambda|$ sufficiently large, one has
\begin{multline}
  \label{eq:80}
  \sup_{\omega\in\mathcal{Z}_\Lambda}
    \left(DCS'(s-\varepsilon;I_\Lambda,\Lambda,\omega)-
      DCS(s;I_\Lambda,\Lambda,\omega)\right)\\+
    \sup_{\omega\in\mathcal{Z}_\Lambda}
    \left(DCS(s;I_\Lambda,\Lambda,\omega)-
      DCS'(s+\varepsilon;I_\Lambda,\Lambda,\omega)\right)
  \leq\varepsilon
\end{multline}
where $DCS(s;I_\Lambda,\Lambda,\omega)$ is defined in~\eqref{DCS} and
\begin{equation*}
  DCS'(s;I_\Lambda,\Lambda,\omega)=\frac1{N(E_0+I_\Lambda,\Lambda,\omega)}
  \sum_{\gamma\in\Gamma}X_\gamma.
\end{equation*}
As $0\leq DCS(s;I_\Lambda,\Lambda,\omega)\leq1$ and $0\leq
DCS'(s;I_\Lambda,\Lambda,\omega)\leq 2$ (for $\Lambda$ large), in view
of~\eqref{eq:80}, to prove~\eqref{eq:76}, it suffices to prove
\begin{equation}
  \label{eq:81}
  \esp\left(\left|DCS'(s;I_\Lambda,\Lambda,\omega)-e^{-s^d} 
    \right|^2\right)\vers_{\Lambda\nearrow\R^d} 0.
\end{equation}
As the Hamiltonians $H_\omega(\Lambda_\ell(\gamma))$ and
$H_\omega(\Lambda_\ell(\gamma))$ are independent when
$\gamma\not=\gamma'$, $(\gamma,\gamma')\in\Gamma^2$, using
Lemma~\ref{lemasympt}, we compute
\begin{align*}
  \E(X_\gamma) &= (1 - N(E_0+I_\Lambda)\ell^d)^{
    (N(E_0+I_\Lambda)\ell^d)^{-1} s^d - 1} N(E_0+I_\Lambda)
  \ell^d + o(|I_\Lambda| \ell^d) \\ & = \e^{-s^d} N(E_0+I_\Lambda)
  \ell^d + o(|I_\Lambda| \ell^d).
\end{align*}
Thus, using Theorem~\ref{thr:16},
one computes
\begin{equation}
  \label{eq:39}
  \begin{split}
    \E(DCS'(s;I_\Lambda,\Lambda,\omega))&=
    \E\left(\frac1{N(E_0+I_\Lambda,\Lambda,\omega)}
      \sum_{\gamma\in\Gamma}X_\gamma\right)
    \\&= \frac1{N(E_0+I_\Lambda)|\Lambda|}\frac{|\Lambda|}{\ell^d}
    \e^{-s^d} N(E_0+I_\Lambda)\ell^d + o(1)= \e^{-s^d}+o(1).
  \end{split}
\end{equation}
On the other hand, by~\eqref{eq:79} and their definition, $X_\gamma$
and $X_{\gamma'}$ are independent when $|\gamma'-\gamma|\leq
2(N(E_0+I_\Lambda))^{-1}s+1$. Hence, using Theorem~\ref{thr:16}, one
computes
\begin{multline*}
  \Delta V:=\E(DCS'(s;I_\Lambda,\Lambda,\omega)^2)
  -\E(DCS'(s;I_\Lambda,\Lambda,\omega))^2
  \\=\frac{1+o(1)}{(N(E_0+I_\Lambda)|\Lambda|)^2}
  \sum_{\gamma\in\Gamma}\sum_{\substack{\gamma'\in\Gamma\\|\gamma'-\gamma|\leq
      2(N(E_0+I_\Lambda))^{-1}s+1}}(\E(X_\gamma
  X_{\gamma'})-\E(X_\gamma)\E( X_{\gamma'})).
\end{multline*}
Thus, as $\#\Gamma\leq C|\Lambda|/\ell^d$, by Lemma~\ref{lemasympt},
one has
\begin{equation}
  \label{eq:40}
  \Delta V\leq \frac2{(N(E_0+I_\Lambda)|\Lambda|)^2}
  \frac{|\Lambda|}{\ell^d}\frac1{\ell^dN(E_0+I_\Lambda)} N(E_0+I_\Lambda)\ell^d=
  \frac2{N^2(E_0+I_\Lambda)|\Lambda|\ell^d}.
\end{equation}
Condition~\eqref{eq:73} then ensures that the choice of the length
$\ell$ in Theorem~\ref{thr:vbig1} can be made such that
$N^2(E_0+I_\Lambda)|\Lambda|\ell^d\to+\infty$ as $|\Lambda|\to+\infty$.\\
Thus,~\eqref{eq:39} and~\eqref{eq:40} imply~\eqref{eq:81}. This
completes the proof of Theorem~\ref{thr:DCS}.\qed



%
\section{Proof of Proposition~\ref{pro:1} and Theorem~\ref{thr:10}}
\label{sec:proof-prop-theor}
We start with the proof of Proposition~\ref{pro:1} then we prove
Theorem~\ref{thr:10}.
\subsection{The proof of Proposition~\ref{pro:1}}
\label{sec:proof-proposition1}
Let us prove point (1). Consider $x(E)$ and $x'(E)$ two centers of
localization for some energy $E\in I$. Let $\varphi$ be a normalized
eigenstate associated to $E$. Then, by the assumption (Loc'), we know
that, for all $x$,
\begin{equation*}
  \|\varphi\|^2_x\leq C_\omega \langle
  x(E)\rangle^{q}e^{-\gamma
    |x-x(E)|^\xi} \langle
  x'(E)\rangle^{q}e^{-\gamma
    |x-x'(E)|^\xi}.
\end{equation*}
Hence, summing over $x$, we get that
\begin{equation*}
  1\leq C_d\, C_\omega\, \langle
  x(E)\rangle^{q}\langle
  x'(E)\rangle^{q}\frac{1}{\gamma^{d/\xi}}e^{-\gamma
    |x(E)-x'(E)|^\xi}.
\end{equation*}
Taking the logarithm, we immediately get~\eqref{eq:8}.\\
Let us now prove point (2). First, note that, by the Wegner estimate
(W), the condition $N(I_L)L^d\to+\infty$ implies that $|I_L|\gtrsim
L^{-d}$ for $L$ large. Hence, if $J$ is an interval centered at $0$
such that $|J|=o(L^{-d})$, then $N(I_L+J)=N(I_L)(1+o(1))$.\\
By definition, one has
\begin{equation*}
  N(I_L,L)=\sum_{\substack{n;\ E_n\in
      I_L\\\gamma_n\in\Lambda_L}}\|\varphi_n\|^2.
\end{equation*}
where, for $n$, $E_n$ is an eigenvalues of $H_\omega$, $\varphi_n$ an
associated eigenfunction and $\gamma_n$ a localization center of
$\varphi_n$ in $\Lambda_L$.\\
Fix $\alpha\in(0,1)$. Using the estimate~\eqref{eq:7} in the same way
as the estimates (Loc) was used in the proof of Lemma~\ref{lemcenter},
we get that
\begin{equation*}
  \left|N(I_L,L)-\sum_{\substack{n;\ E_n\in
        I_L\\\gamma_n\in\Lambda_L}}
    \|\varphi_n\|^2_{\Lambda_{L+L^\alpha}}\right|\leq
    C_\omega L^{q+d} e^{-2\gamma L^{\alpha\xi}}
\end{equation*}
where $\|\cdot\|_{\Lambda_{L+L^\alpha}}$ denotes the $L^2$-norm on
the cube $\Lambda_{L+L^\alpha}$.\\
Hence, we have
\begin{equation*}
  \begin{split}
  N(I_L,L)&\leq \sum_{\substack{n;\ E_n\in
        I_L\\\gamma_n\in\Lambda_L}}
    \|\varphi_n\|^2_{\Lambda_{L+L^\alpha}}+
    C_\omega L^{q+d} e^{-2\gamma L^{\alpha\xi}}\\
    &\leq\tr(\car_{I_L}(H_\omega)\car_{\Lambda_{L+L^\alpha}})+
    C_\omega L^{q+d}e^{-2\gamma L^{\alpha\xi}}.
  \end{split}
\end{equation*}
By standard estimates on Schr{\"o}dinger operators (see
e.g.~\cite{MR2154153}), we know that
\begin{equation}
  \label{eq:74}
  0\leq \tr(\car_{I_L}(H_\omega)\car_{\Lambda_{L+L^\alpha}})\leq CL^d.
\end{equation}
Pick $\chi$ a smooth cut-off function with gradient vanishing outside
$\Lambda_{L+L^\alpha}\setminus \Lambda_{L+L^\alpha/2}$. Then, using
the localization estimate (Loc'), one easily checks that, for $L$
sufficiently large
\begin{itemize}
\item for $n$ such that $E_n\in I_L$ and $\gamma_n\in\Lambda_L$,
  $\chi\varphi_n$ satisfies
  \begin{equation}
    \label{eq:98}
    \|(H_\omega(\Lambda_{L+L^\alpha})-E_n)(\chi\varphi_n)\|\leq
    e^{-L^{\alpha\xi}/C}; 
  \end{equation}
\item the Gram matrix for the family $(\chi\varphi_n)_{E_n\in I_L,\
    \gamma_n\in\Lambda_L}$ satisfies
  \begin{equation}
    \label{eq:75}
    ((\langle\chi\varphi_n,\chi\varphi_m\rangle))_{\substack{E_n\in I_L,\
    \gamma_n\in\Lambda_L\\E_m\in I_L,\
    \gamma_m\in\Lambda_L}}=Id+O\left(e^{-L^{\alpha\xi}/C}\right).
  \end{equation}
\end{itemize}
As~\eqref{eq:74} implies that the number of $n$ such that $E_n\in I_L$
and $\gamma_n\in\Lambda_L$ is bounded by $CL^d$, we get that
\begin{equation*}
  N(I_L,L)\leq \tr(\car_{\tilde I_L}(H_\omega(\Lambda_{L+2L^\alpha})))+
  C_\omega L^{q+d}e^{-2\gamma L^{\alpha\xi}}
\end{equation*}
where $\tilde I_L=I_L+[-e^{L^\alpha\xi/C},e^{L^\alpha\xi/C}]$.\\
Hence, if, as in section~\ref{sec:cutting-pieces}, $N(\omega,
\Lambda_{L+L^\alpha},, I_L)$ denotes the number of eigenvalues of
$H_\omega(\Lambda_{L+L^\alpha})$ in $I_L$, we have
\begin{equation}
  \label{eq:78}
  N(I_L,L)\leq N(\omega,\Lambda_{L+L^\alpha},, I_L)+
  C_\omega L^{q+d}e^{-2\gamma L^{\alpha\xi}}.
\end{equation}
Recall that we assumed $N(I_L)|I_L|^{-1-\rho}\to+\infty$ where $\rho$
is given by (M). The proof of Theorem~\ref{thr:16} then also shows
that, if $N(\omega, \Lambda_{L+L^\alpha}, \Lambda_L, I_L)$ denotes the
number of eigenvalues of $H_\omega(\Lambda_{L+L^\alpha})$ in $I_L$
having a localization center in $\Lambda_L$, then, for any $q$,
\begin{equation*}
  \begin{split}
  \P\left\{  \left| N(\omega, \Lambda_{L+L^\alpha}, \Lambda_L, I_L) -
      N(I_\Lambda)|\Lambda| \right| = o(N(I_\Lambda)|\Lambda|)
  \right\} &\geq 1 - C_q (N(I_\Lambda)|\Lambda|)^{-q}\\
  &\geq 1 - C_q L^{-2}     
  \end{split}
\end{equation*}
as $N(I_L)L^{d-\varepsilon}\to+\infty$ for some $\varepsilon>0$.\\
So, by the Borel-Cantelli Lemma and Theorem~\ref{thr:16}, we know
that, $\omega$-almost surely,
\begin{equation}
  \label{eq:77}
  N(\omega, \Lambda_{L+L^\alpha}, \Lambda_L, I_L)=
  N(\omega, \Lambda_{L+L^\alpha},, I_L)(1+o(1))=N(I_L)L^d(1+o(1)).
\end{equation}
One has
\begin{equation*}
  N(I_L,L)\geq\sum_{\substack{n;\ E_n\in I_L\\\gamma_n\in\Lambda_L}}
  \|\varphi_n\|^2_{\Lambda_{L+L^\alpha}}.
\end{equation*}
Using the same cut-off of the eigenfunctions as above, we then see that
\begin{equation*}
  N(I_L,L)\geq N(\omega, \Lambda_{L+L^\alpha}, \Lambda_L, I_L).
\end{equation*}
Using this estimate,~\eqref{eq:77},~\eqref{eq:78} and
Theorem~\ref{thr:16}, we get~\eqref{eq:9} and complete the proof of
Proposition~\ref{pro:1}.\qed
\subsection{The proof of Theorem~\ref{thr:10}}
\label{sec:proof-theorem}
Let us now consider consider $\omega$ in the full measure set where
the statements of Theorems~\ref{thr:1},~\ref{thr:9} and
Proposition~\ref{pro:1} hold. As in
section~\ref{sec:proof-proposition1}, $\omega$-almost surely, for $L$
large, if $E$ is an eigenvalue of $H_\omega$ with localization center
in $\Lambda_L$, then there exists $E'$ an eigenvalue of
$H_\omega(\Lambda_{L+L^\alpha})$ such that $|E-E'|\leq
e^{-L^\alpha/C}$. To avoid confusion, rename $DLS(x;I_L,\omega,L)$
defined by~\eqref{eq:4} to $DLS^f(x;I_L,\omega,L)$.  Hence, for
$\varepsilon>0$ fixed, we have, for $L$ large enough,
\begin{itemize}
\item when $|I_L|\to0$ satisfying the assumptions of
  Theorem~\ref{thr:10}, we have
  \begin{equation*}
    \begin{split}
    DLS(\nu(E_0)\,x-\varepsilon;I_L,\omega,L)&\geq
    \frac{N(I_L,L)}{N(I_L,\omega,\Lambda_L)}\cdot DLS^f(x;I_L,\omega,L)\\&\geq 
    DLS(\nu(E_0)\,x+\varepsilon;I_L,\omega,L);
    \end{split}
  \end{equation*}
\item when $I_L=J$ a fixed interval satisfying the assumptions of
  Theorem~\ref{thr:10}, we have
  \begin{equation*}
    DLS'(x-\varepsilon;J,\omega,L)\geq
    \frac{N(I_L,L)}{N(J,\omega,\Lambda_L)}
    DLS^f(x;J,\omega,L)\geq DLS'(x+\varepsilon;J,\omega,L).
  \end{equation*}
\end{itemize}
Thus, Theorem~\ref{thr:10} is an immediate corollary of
Theorems~\ref{thr:1} and~\ref{thr:9}, and the second point of
Proposition~\ref{pro:1}. \qed



\section{Appendix}
\label{sec:appendix}
\subsection{The proof of Lemma~\ref{le:5}}
\label{sec:proof-lemma5}
Compute
\begin{equation*}
  \begin{split}
  \int_\R\left|\frac{N(E+x)-N(E+y)}{x-y}-\nu(E)\right|&=
  \int_\R\left|\frac1{x-y}\int_{E+y}^{E+x}(\nu(e)-\nu(E))de\right|
  \\&\leq\int_0^1\left(\int_\R|\nu(E+y+(x-y)e)-\nu(E)| dE\right)de.
  \end{split}
\end{equation*}
As $\nu\in L^1(\R)$, Lebesgue's dominated convergence theorem ensures
that the last integral converges to $0$ when $x$ and $y$ tend to
$0$. Hence, the quotient in the first integral converges to $0$ for
almost every $E$. This completes the proof of Lemma~\ref{le:5}.\qed
\subsection{Local control on the exponential decay of eigenfunctions}
\label{sec:append-local-contr}
In this section, we establish SUDEC and SULE estimates for the
eigenfunctions assciated to an eigenvalue in the localized regime of a
random operator restricted to a finite volume. These are the analogues
of the infinite volume estimates introduced in \cite{MR97m:47002} and
proved in~\cite{MR97m:47002,MR2203782}.
We denote by $\Sigma_{\mathrm{CL}}$ the region of complete
localization for the random operator $H_\omega$. It is defined as the
set of energies $E \in \I$ where we have all the conclusions of the
bootstrap multiscale analysis of~\cite{MR1859605} or the fractional
moment method of~\cite{MR1244867,MR1301371}. These conclusions turn
out to be equivalent properties describing the localized
regime~\cite{MR2078370,MR2203782}.  In Theorem~\ref{thmFVE}, we
provided new equivalent characterizations of the region of complete
localization, involving this time finite volume operators, rather than
the infinite volume one.
\noindent We prove
\begin{Th}
  \label{thmFVE}
  Let $I\subset\Sigma$ be a compact interval and assume that Wegner's
  estimate {\bf (W)} holds in $I$. For $L$ given, consider $\Lambda=
  \Lambda_L(0)$ a cube of side length $L$ centered at $0$, and denote
  by $\varphi_{\omega,\Lambda,j}$, $j=1,\cdots,
  \mathrm{tr}\car_I(H_\omega(\Lambda))$, the normalized eigenvectors
  of $H_\omega(\Lambda)$ with corresponding eigenvalue in $I$. The
  following are equivalent
  \begin{enumerate}
  \item $I\subset \Sigma_{\mathrm{CL}}$

  \item For all $E\in I$, there exists $\theta>3d-1$,
    \begin{equation}
      \label{startMSA}
      \limsup_{L\to\infty} \P\left\{\forall x,y\in\Lambda, \,
        |x-y|\ge \frac L2, \;  \| \chi_x (H_\omega(\Lambda) - E)^{-1}
        \chi_y  \|  \le L^{-\theta} \right\} =1. 
    \end{equation}
  \item For all $\xi\in(0,1)$,
    \begin{equation}
      \label{FVsudecE}
      \sup_{y\in \Lambda} \E\left\{ \sum_{x\in\Lambda} \e^{|x-y|^\xi}
        \sup_j \|  \varphi_{\omega,\Lambda,j} \|_x \|
        \varphi_{\omega,\Lambda,j} \|_y  \right\}<\infty. 
    \end{equation}
  \item There exists $\xi\in(0,1)$,
    \begin{equation}
      \label{FVsudecE2}
      \sup_{y\in \Lambda} \E\left\{ \sum_{x\in\Lambda}  \e^{|x-y|^\xi}
        \sup_j \|  \varphi_{\omega,\Lambda,j} \|_x \|
        \varphi_{\omega,\Lambda,j} \|_y  \right\}<\infty. 
    \end{equation}
  \item For all $\xi\in(0,1)$,
    \begin{equation}
      \label{FVdynloc1}
      \sup_{y\in \Lambda} \E\left\{ \sum_{x\in\Lambda}  \e^{|x-y|^\xi}
        \sup_{\substack{\supp f \subset I \\ |f|\le 1 }} \| \chi_x
        f(H_\omega(\Lambda)) \chi_y \|_2  \right\}<\infty. 
    \end{equation}
  \item There exists $\xi\in(0,1)$,
    \begin{equation}
      \label{FVdynloc2}
      \sup_{y\in \Lambda} \E\left\{ \sum_{x\in\Lambda}  \e^{|x-y|^\xi}
        \sup_{\substack{\supp f \subset I \\ |f|\le 1 }} \| \chi_x
        f(H_\omega(\Lambda)) \chi_y \|_2  \right\}<\infty.
    \end{equation}
  \item There exists $\xi\in(0,1)$,
    \begin{equation}
      \label{FVdynloc3}
      \sup_{y\in \Lambda} \sup_{\substack{\supp f \subset I \\ |f|\le 1 }}
      \E\left\{ \sum_{x\in\Lambda}  \e^{|x-y|^\xi} \| \chi_x
        f(H_\omega(\Lambda)) \chi_y \|_2  \right\}<\infty. 
    \end{equation}
  \item (SUDEC for finite volume with polynomial probability) For all $p>d$, there is $q=q_{p,d}$ so that  for all $\xi\in(0,1)$,
  for any $L$ large enough, the
    following holds with probability at least $1 - L^{-p}$: for any
    eigenvector $\varphi_{\omega,\Lambda,j}$ of $H_{\omega,\Lambda}$,
    with energy in $I$, for any $(x,y)\in\Lambda^2$,  one has
    \begin{equation}
      \label{FVsudecP}
      \|  \varphi_{\omega,\Lambda,j} \|_x \|
      \varphi_{\omega,\Lambda,j} \|_y\le L^{q} \e^{-|x-y|^{\xi}}. 
    \end{equation}
  \item (SULE for finite volume with polynomial probability) For all $p>d$, there is $q=q_{p,d}$ so that  $\xi\in(0,1)$,
   for any $L$ large enough, the
    following holds with probability at least $1 - L^{-p}$: for any
    eigenvector $\varphi_{\omega,\Lambda,j}$ of $H_{\omega,\Lambda}$,
    with energy in $I$, there is a center of localization
    $x_{\omega,\Lambda,j}\in \Lambda$, so that for any $x\in \Lambda$, one has
    \begin{equation}
      \label{FVsule}
      \|  \varphi_{\omega,\Lambda,j} \|_x \le L^q
      \e^{-|x-x_{\omega,\Lambda,j}|^{\xi}}. 
    \end{equation}
  \item (SUDEC for finite volume with subexponential probability) For all $\nu,\xi\in(0,1)$, $\nu<\xi$,
  for any $L$ large enough, the
    following holds with probability at least $1 - \e^{-L^\nu}$: for any
    eigenvector $\varphi_{\omega,\Lambda,j}$ of $H_{\omega,\Lambda}$,
    with energy in $I$, for any $(x,y)\in\Lambda^2$,  one has
    \begin{equation}
      \label{FVsudecP1}
      \|  \varphi_{\omega,\Lambda,j} \|_x \|
      \varphi_{\omega,\Lambda,j} \|_y \le \e^{2L^\nu} \e^{-|x-y|^{\xi}}. 
    \end{equation}
  \item (SULE for finite volume with subexponential probability)For all $\nu,\xi\in(0,1)$, $\nu<\xi$,   for any $L$ large enough, the
    following holds with probability at least $1 - \e^{-L^\nu}$: for any
    eigenvector $\varphi_{\omega,\Lambda,j}$ of $H_{\omega,\Lambda}$,
    with energy in $I$, there is a center of localization
    $x_{\omega,\Lambda,j}\in \Lambda$, so that for any $x\in \Lambda$, one has
    \begin{equation}
      \label{FVsule1}
      \|  \varphi_{\omega,\Lambda,j} \|_x \le  \e^{2L^\nu} 
      \e^{-|x-x_{\omega,\Lambda,j}|^{\xi}}. 
    \end{equation}
  \end{enumerate}
Moreover one can pick $q=p+d$ in (8) and $q=p+\frac 32 d$ in (9).
\end{Th}
\begin{Rem}
  \label{rem:6}
  Theorem~\ref{thmFVE} provides a finite volume analog
  of~\cite[Theorem~1 and Corollary~3]{MR2203782}. If generalized
  eigenfunctions are needed in the infinite volume case, the
  normalized eigenfunctions are good enough in the finite volume one
  for the spectrum intersected with the interval $I$ is discrete.\\
  For discrete models, using e.g. the finite volume fractional moment
  criteria of~\cite{MR2002h:82051}, one can derive bounds of the same
  type as~\eqref{FVsule} where $\xi=1$. This has been done
  in~\cite{Kl:10}.
\end{Rem}
\begin{proof}[Proof of Theorem~\ref{thmFVE}]
  We start by describing precisely $\Sigma_{\mathrm{CL}}$, which is
  the set of energies where the conclusion of the bootstrap multiscale
  analysis of \cite{MR1859605} is valid: $\Sigma_{\mathrm{CL}}$ for
  the random operator $H_\omega$ is defined as the set of $E \in I$
  for which there exists some open interval $I(E)\subset I$, with
  $E\in I(E) $, such that, given any $\zeta\in(0,1)$ and
  $\alpha\in(1,\zeta^{-1})$, there exists a length scale $L_0\in 6
  \mathbb{N}$ and a mass $m>0$, so that if we set $L_{k+1} =
  [L_k^\alpha]_{6\mathbb{N}}$ for $k=0,1,\dots$, we have
  \begin{equation}
    \label{MSAest} 
    \mathbb{P}\,\left\{R\left(m, L_k,
        I(E),x,y\right) \right\}\ge 1 -\mathrm{e}^{-L_{k}^\zeta}
  \end{equation} 
  for all $k=0,1,\ldots$, and $x, y \in \mathbb{Z}^d$ with $|x-y| >
  L_k +R_0$, where
  \begin{equation}
    \label{defsetR}
    R(m,L, I,x,y) = 
    \left\{\omega;\quad
      \begin{aligned}
        \mbox{for every}\; E^\prime \in I \ \mbox{, either} \
        \Lambda_L(x) \\ \mbox{or} \ \Lambda_L(y) \ \mbox{is} \
        \mbox{$(\omega,m, E^\prime)$-regular}
      \end{aligned}
    \right\}.
  \end{equation} 
  Here $[K]_{6\mathbb{N}}= \max \{ L \in 6\mathbb{N}; \; L \le K\}$;
  we work with scales in $6\N$ for convenience; $R_0>0$ is given in
  Assumption (IAD).\\
  Given $E \in \mathbb{R}$, $x \in \mathbb{Z}^d$ and $L \in 6
  \mathbb{N}$, we say that the box $\Lambda_L(x) $ is $(\omega,m,
  E)$-regular for a given $m>0$ if $E \notin
  \sigma(H_{\omega}(\Lambda_L(x) )$ and
  \begin{equation}
    \| \Gamma_{x,L} R_{\omega,x,L}(E)\chi_{x,{{\frac L 3}}} 
    \| \le {\rm e}^{-mL}, \label{regular}
  \end{equation}
  where $R_{\omega,x,L}(E) = (H_{\omega}(\Lambda_L(x) )-E)^{-1}$ and
  $\Gamma_{x,L}$ denotes the charateristic function of the ``belt''
  $\overline{\Lambda}_{L-1}(x) \backslash {\Lambda}_{L-3}(x) $ when
  $\mathcal{H}=\mathrm{L}^2(\mathbb{R}^d,{\rm d}x)$ (the arguments can
  be easily modified for $\mathcal{H}=\ell^2(\mathbb{Z}^d)$).\\
  The interval $I$ in Theorem~\ref{thmFVE} being compact, we can
  extract from $\cup_{E\in I} I(E)\supset I$ a finite number of
  intervals $I(E)$ that cover $I$. Thus, with no loss, we may assume
  that \eqref{MSAest} is valid for the interval $I$ itself.\\
  We turn to the proof of Theorem~\ref{thmFVE} per se. That
  $(1)\Longleftrightarrow (2)$ is due to \cite{MR2078370}. Note that
  if~\eqref{startMSA} holds, then
  \begin{equation*}
    \limsup_{L\to\infty} \P\left\{\|
      \Gamma_{0,L} R_{\omega,0,L}(E)\chi_{0,{{\frac L 3}}} \| \le
      L^{-\theta-2d+1} \right\} =1. 
  \end{equation*}
  Assume (2), that is \eqref{MSAest} holds for the interval $I$. We
  show that (3) holds. A standard computation (see
  \cite[(EDI)]{MR2078370}, \cite[Eq.(8.10)]{MR2509110}) shows that if
  $\Lambda_L(x)$ is $(\omega,m, E)$-regular, then, if
  $H_{\omega}(\Lambda_L(x) )\varphi_{\omega,\Lambda,j}=\tilde
  E\varphi_{\omega,\Lambda,j}$ and $\|\varphi_{\omega,\Lambda,j}\|=1$,
  then
  \begin{align}
    \| \varphi_{\omega,\Lambda,j} \|_x \lesssim L^{d-1} \|
    \Gamma_{x,L} R_{\omega,x,L}(\tilde E)\chi_{x,{{\frac L 3}}} \| \le
    {\rm e}^{-\frac{m}{2}L} ,
  \end{align}
  whenever $L$ is large enough, depending only on $d,m$.\\
  It follows that for any $(x,y)\in\Lambda^2$, if $k$ is so that
  $L_k\le |x-y| < L_{k+1}$, then for any $\omega\in R(m,L, I,x,y)$,
  for any normalized $\varphi_{\omega,\Lambda,j}$,
  \begin{equation}
    \| \varphi_{\omega,\Lambda,j} \|_x \| \varphi_{\omega,\Lambda,j}
    \|_y \le {\rm e}^{-\frac{m}{2}L_k} \le {\rm
      e}^{-\frac{m}{2}|x-y|^{\frac 1\alpha}} .
  \end{equation}
  As a consequence,
  \begin{equation}
    \E \left\{ \sup_j\| \varphi_{\omega,\Lambda,j} \|_x\|
      \varphi_{\omega,\Lambda,j} \|_y \right\} \le {\rm
      e}^{-\frac{m}{2}|x-y|^{\frac 1\alpha}} + {\rm
      e}^{-|x-y|^{\frac{\zeta}{\alpha}}} \le 2 {\rm
      e}^{-|x-y|^{\frac{\zeta}{\alpha}}} ,
  \end{equation}
  for $|x-y|$ large enough. Thus,~\eqref{FVsudecE} follows for any
  $\xi<\zeta \alpha^{-1}$. Since $\zeta<1<\alpha$ can be chosen
  arbitrarily close to 1,~(3) follows.\\
  (3) clearly implies (4).\\
  To see that (3) implies (5), it is enough to decompose the operator
  $f(H_\omega(\Lambda))$ over the eigenvectors
  $(\varphi_{\omega,\Lambda,j})_j$ and note that $ \| \chi_x
  \Pi_{\varphi_{\omega,\Lambda,j}} \chi_y \|_2= \|
  \varphi_{\omega,\Lambda,j} \|_x\| \varphi_{\omega,\Lambda,j} \|_y$,
  when $ \Pi_{\varphi_{\omega,\Lambda,j}}$ is the orthogonal
  projection on $ \varphi_{\omega,\Lambda,j}$.\\
  That $(5)\Longrightarrow (6) \Longrightarrow (7)$ is immediate.\\
  Assume (3). We prove (8). We have
  \begin{align}
    & \P\left( \exists (x,y)\in\Lambda^2, \, {\rm e}^{|x-y|^\xi}
      \sup_j \| \varphi_{\omega,\Lambda,j} \|_x\|
      \varphi_{\omega,\Lambda,j} \|_y \ge L^{p+d} \right)  \\
    & \quad \le L^d \sup_{y\in \Lambda} \sum_{x\in\Lambda} \P\left(
      {\rm e}^{|x-y|^\xi} \sup_j\| \varphi_{\omega,\Lambda,j} \|_x\|
      \varphi_{\omega,\Lambda,j} \|_y \ge L^{p+d} \right) \\
    & \quad \le L^{-p} \sup_{y\in \Lambda} \sum_{x\in\Lambda} \E\left(
      {\rm e}^{|x-y|^\xi} \sup_j \| \varphi_{\omega,\Lambda,j} \|_x\|
      \varphi_{\omega,\Lambda,j}
      \|_y \right) \\
    & \quad \lesssim L^{-p}.
  \end{align}
  In other terms, with probability at least $1 - L^{-p}$, we have $\|
  \varphi_{\omega,\Lambda,j} \|_x\| \varphi_{\omega,\Lambda,j} \|_y \|
  \le L^{p+d} {\rm e}^{-|x-y|^\xi}$, for any $j$ and
  $(x,y)\in\Lambda^2$, and (8) holds.\\
  We show that (8) and (9) are equivalent. Assume (8) and let
  $x_{\omega,\Lambda,j}$ be a center of localization for
  $\varphi_{\omega,\Lambda,j}$. Since
  $\|\varphi_{\omega,\Lambda,j}\|=1$, note that
  $\|\varphi_{\omega,\Lambda,j} \|_x \gtrsim L^{-d/2} $. We
  write~\eqref{FVsudecP} with $y=x_{\omega,\Lambda,j}$
  and~\eqref{FVsule} follows from the last observation with a constant $q'=q+\frac d2$. Conversely, if
  (9) holds, then it follows from~\eqref{FVsule} that
  \begin{align}
    \| \varphi_{\omega,\Lambda,j} \|_x\| \varphi_{\omega,\Lambda,j}
    \|_y \le L^{2q} \e^{-|x-x_{\omega,\Lambda,j}|^{\xi}
    }\e^{-|y-x_{\omega,\Lambda,j}|^{\xi}}\le L^{2q}\e^{-|x-y|^{\xi}}.
  \end{align}
 We show that $(8)\Longrightarrow
  (2)$. 
  Assume (8). For any given $E\in I$, thanks to (W), we may assume
  that $\mathrm{d}(E,\sigma(H_\omega(\Lambda)))\ge L^{-p}$, with
  probability at least $1-L^{-p+d}$. As a consequence, for any $L$ large
  enough, for any $(x,y)\in\Lambda^2$, $|x-y| \ge L/2$,
  \begin{align}
    \| \chi_x (H_\omega(\Lambda) - E)^{-1} \chi_y \| & \le L^p
    \sum_j\| \varphi_{\omega,\Lambda,j} \|_x\|
    \varphi_{\omega,\Lambda,j} \|_y\\& \le L^{p+q} \e^{-L^\xi }.
  \end{align}
 The estimate (2) follows, regardless of the values of $p,q,\xi$.
 \\
  The proofs related to (10) and (11) are similar.
  \\ 
 It remains to show $(7)\Longrightarrow (2)$.  Assume (7) and pick $E\in I$. Set $\delta=L^{-p}$. By (W), with
  probability at least $1-L^{-p}$, we have
  \begin{equation}
    (H_\omega(\Lambda) - E)^{-1} =
    (H_\omega(\Lambda) - E)^{-1} \car_{\Sigma\setminus
      [E-\delta,E+\delta] } (H_\omega(\Lambda)) . 
  \end{equation}
  Define the function $f_{E,\delta}(\lambda)= \delta (\lambda -
  E)^{-1} \car_{\Sigma\setminus [E-\delta,E+\delta] } (\lambda)$. Note
  that $|f_{E,\delta} | \le 1$. It follows that
  \begin{align}
    & \P \left(\exists (x,y)\in\Lambda^2, \, |x-y|\ge \frac L2, \, \|
      \chi_x (H_\omega(\Lambda) - E)^{-1} \chi_y \| \ge L^{-\theta}
    \right) \\
    & \quad \le L^{-p} + \sum_{\substack{(x,y)\in\Lambda^2 \\ |x-y|\ge
        L/2}} \P \left( \| \chi_x f_{E,\delta}(H_\omega(\Lambda) )
      \chi_y \|
      \ge L^{-\theta-p} \right) \\
    & \quad \le L^{-p } + L^{-p } \sup_{y\in\Lambda} \sum_{\substack{x\in\Lambda
      \\ |x-y|\ge L/2}} L^{\theta+2p+d} \E \left( \| \chi_x
      f_{E,\delta}(H_\omega(\Lambda) ) \chi_y \| \right)\\
    & \quad \le L^{-p } + L^{-p } \sup_{y\in\Lambda}
    \sum_{x\in\Lambda} (2|x-y|)^{\theta+2p+d} \E \left( \| \chi_x
      f_{E,\delta}(H_\omega(\Lambda) ) \chi_y \|_2 \right) \lesssim
    L^{-p} .
  \end{align}
 The last observation stated in the theorem follows from the proof above.
\end{proof}




%
\newcommand{\etalchar}[1]{$^{#1}$}
\def\cprime{$'$} \def\cydot{\leavevmode\raise.4ex\hbox{.}} \def\cprime{$'$}


\end{document}